\documentclass[11pt]{amsart}

\usepackage{amscd}
\usepackage{amssymb}
\usepackage{booktabs}
\usepackage{graphicx}
\PassOptionsToPackage{table}{xcolor}

\usepackage[shortlabels]{enumitem}

\newcommand{\proofsec}[1]{\smallskip\emph{\underline{#1}}}
\newcounter{stepnum}
\newcommand{\thestep}{\arabic{stepnum}}
\newcommand{\step}[1]{\proofsec{Step \refstepcounter{stepnum}\thestep}: \emph{#1}}

\usepackage{times}

\usepackage{color}

 \usepackage[T1]{fontenc}  
 \newcommand{\ditto}{\ \textquotedbl \ }

\DeclareMathOperator{\RE}{Re}
\DeclareMathOperator{\IM}{Im}
\renewcommand{\Re}{\RE}

\usepackage{dynkin-diagrams}

\usepackage[colorlinks=true, pdfstartview=FitV, linkcolor=blue, pagebackref, citecolor=blue, urlcolor=blue]{hyperref}
\usepackage{xurl}

\usepackage{orcidlink}

\usepackage{amsthm}
\newtheorem{thm}[equation]{Theorem}
\newtheorem{lem}[equation]{Lemma}
\newtheorem{cor}[equation]{Corollary}
\newtheorem{prop}[equation]{Proposition}
\newtheorem{conj}[equation]{Conjecture}

\newtheorem*{thm*}{Theorem}
\newtheorem*{prop*}{Proposition}
\newtheorem*{cor*}{Corollary}
\newtheorem*{lem*}{Lemma}
\newtheorem*{MT*}{Main Theorem}

\newtheorem*{ques*}{Question}
\newtheorem*{claim*}{Claim}

\theoremstyle{definition} %
\newtheorem{defn}[equation]{Definition}
\newtheorem*{defn*}{Definition}
\newtheorem{eg}[equation]{Example}

\theoremstyle{remark} %
\newtheorem{rmk}[equation]{Remark}

\newtheorem*{rmk*}{Remark}
\newtheorem{rem}[equation]{Remark}
\newtheorem*{rmks*}{Remarks}

\newcommand{\stacklink}[1]{\href{https://stacks.math.columbia.edu/tag/#1}{#1}}
\newcommand{\stacks}[1]{\cite[Tag \stacklink{#1}]{stacks-project}}

\newcommand{\trans}{\mathsf{T}}   

\newcommand{\Gm}{\mathbb{G}_{\mathrm{m}}}

\newcommand{\A}{\mathbb{A}}
\newcommand{\R}{\mathbb{R}}
\newcommand{\F}{\mathbb{F}}
\newcommand{\Q}{\mathbb{Q}}
\newcommand{\E}{\mathbb{E}}

\newcommand{\C}{\mathbb{C}}
\newcommand{\Z}{\mathbb{Z}}

\newcommand{\mcD}{\mathcal{D}}
\newcommand{\mcO}{\mathcal{O}}
\newcommand{\cO}{\mcO}
\newcommand{\cOh}{\hat{\mcO}}

\newcommand{\ic}{{\circ}}

\newcommand{\kalg}{\bar{k}}
\newcommand{\Qalg}{\bar{\Q}}

\newcommand{\kh}{\hat{k}}
\newcommand{\gh}{\hat{g}}

\newcommand{\Mh}{\hat{M}}

\newcommand{\Tt}{\tilde{T}}

\DeclareMathOperator{\Sym}{Sym}
\DeclareMathOperator{\Hom}{Hom}
\DeclareMathOperator{\Res}{Res}
\DeclareMathOperator{\Der}{Der}
\DeclareMathOperator{\Reg}{Reg}

\newcommand{\la}{\lambda}

\newcommand{\mfa}{\mathfrak{a}}
\newcommand{\mfe}{\mathfrak{e}}
\newcommand{\mfg}{\mathfrak{g}}
\newcommand{\mft}{\mathfrak{t}}

\renewcommand{\:}{\colon}

\DeclareMathOperator{\tr}{Tr}
\DeclareMathOperator{\sgn}{sgn}
\DeclareMathOperator{\diag}{diag}
\DeclareMathOperator{\Diff}{Diff}

\newcommand{\emptyslot}{{-}}

\newcommand{\qform}[1]{{\left\langle{#1}\right\rangle}}                   

\newcommand{\eand}{\quad\text{and}\quad}

\newcommand{\mmin}{\underline{m}}
\newcommand{\mmax}{\overline{m}}

\DeclareMathOperator{\GL}{GL}
\DeclareMathOperator{\PGL}{PGL}

\DeclareMathOperator{\SL}{SL}

\DeclareMathOperator{\Sp}{Sp}
\DeclareMathOperator{\SU}{SU}
\DeclareMathOperator{\SO}{SO}

\DeclareMathOperator{\Spec}{Spec}

\DeclareMathOperator{\class}{Class}
\newcommand{\ClG}{\class(G)}

\DeclareMathOperator{\End}{End}

\DeclareMathOperator{\Tr}{Tr}
\DeclareMathOperator{\Ad}{Ad}
\DeclareMathOperator{\Lie}{Lie}
\DeclareMathOperator{\Aut}{Aut}

\DeclareMathOperator{\Out}{Out}

\DeclareMathOperator{\rank}{rank}

\newcommand{\ov}{\omega^\vee}
\newcommand{\av}{\alpha^\vee}

\newcommand{\sinfty}{c}
\newcommand{\thinfty}{\theta_0}

\def\edgelen{0.75cm}  

\numberwithin{equation}{section}

\begin{document}

\title{Invariant derivations and trace bounds}



\begin{abstract}
About 20 years ago, J-P.~Serre announced a bound on the trace of elements of compact Lie groups under the adjoint representation together with related results,  provided indications of his proofs, and invited a better proof.  This note provides a new, general method for proving such bounds; uses that method to derive Serre's bounds; gives a second proof of Serre's announced results that (we learned) closely follows his original argument; and provides lower bounds for traces of other representations of compact Lie groups and for Brauer characters of finite groups. 
\end{abstract}

\author[S. Garibaldi]{Skip Garibaldi\,\orcidlink{0000-0001-8924-5933}}
\address{Garibaldi: IDA Center for Communications Research-La Jolla, San Diego, CA 92121, USA)}
\email{skip@garibaldibros.com}

\author[R.M. Guralnick]{Robert M. Guralnick\,\orcidlink{0000-0002-9094-857X}}
\address{Guralnick: Department of Mathematics, University of Southern California, Los Angeles, CA 90089-2532, USA}
\email{guralnic@usc.edu}

\author[E.M. Rains]{Eric M. Rains\,\orcidlink{0000-0002-9915-0919}}
\address{Rains: Department of Mathematics, California Institute of Technology, 1200 E California Blvd, Pasadena, CA 91125, USA}
\email{rains@caltech.edu}

\maketitle

\setcounter{tocdepth}{1}
\tableofcontents
\section{Introduction}

During the 2000s, 
J-P.~Serre proved several bounds on the trace of elements in a compact Lie group $G$, see  \ref{devissage.cor}, \ref{simple.thm}, and \ref{short.prop} below.   He announced these results in \cite{Se:trace} but only gave indications to the proofs, writing 
``The classical groups are easy enough, but $F_4$, $E_6$, $E_7$ and $E_8$ are not (especially $E_6$, which I owe to Alain Connes). I hope there is a better proof.''  Part \ref{part.Wders} gives a general approach for proving bounds such as those in \ref{simple.thm} via $G$-invariant derivations on the ring $\C[G]^G$.  We exhibit its use to prove \ref{simple.thm} for groups of type $G_2$ and $E_8$ and provide Magma \cite{Magma} code to verify \ref{simple.thm} for each of the exceptional groups\footnote{Code to replicate the calculations in this paper is provided in the GitHub repo at \url{https://github.com/skipgaribaldi/invariant-derivations}}.  This approach works for optimizing other class functions on a compact $G$, which we illustrate by studying the fundamental irreducible representations of the compact groups of type $F_4$ and $E_8$ in \S\ref{E8.sec}.

Part \ref{part.break} gives a proof of \ref{simple.thm} for the case $s = 1$ and for \ref{short.prop} by a different method involving ad hoc symmetry breaking.  Interestingly, Serre explained to us that this proof is largely the same as his original arguments. It differs in minor details.  For example, when treating the exceptional groups in \S\ref{serre.method}, today we have the luxury to use a few lines of Magma code, whereas 20 years earlier Serre did the same computation by hand.  That is to say an independent event --- the improvement of computer algebra systems over the intervening years --- has improved the proof.  After this paper had been on the arXiv for several months, we learned of the paper \cite{DMMV}, which also contains proofs of \ref{simple.thm} for the case $s=1$ for the groups of types $A$, $D$, and $E$, in a manner substantially similar to what we do here.  We decided to leave this material in this paper in order to give a unified treatment in one place for all the types.  

In contrast to the approach in Part \ref{part.break}, the proof from Part \ref{part.Wders} avoids ad hoc symmetry breaking and stays focused on canonical objects.  It is also most easily adapted to providing bounds on values of characters other than $\Tr \Ad$, as we illustrate.  

Part \ref{completion.part} completes the proofs of Serre's statements \ref{devissage.cor} and \ref{simple.thm}, building on the 
results of Part \ref{part.Wders} or \ref{part.break}.


In Part \ref{min.part}, we expand our view to consider Serre's bounds as a solution to a specific instance of a general problem that we now describe.  Let $\chi := \Tr \rho$ for $\rho \: G \to \GL(V)$ a (real or complex) representation of a compact Lie group $G$.  Since $\rho$ is unitary, the eigenvalues of $\chi(g)$ for $g \in G$ are complex numbers of norm 1, so 
\begin{equation} \label{triv.bounds}
-\chi(1) \le \RE \chi(g) \le \chi(1)
\end{equation}
for all $g \in G$, where $\chi(1) = \dim V$ and $\RE \colon \C \to \R$ gives the real part of a complex number.  As $g$ varies over $G$, the maximum $\chi(1)$ of $\RE \chi(g)$ from \eqref{triv.bounds} is trivially achieved, but the minimum need not be, as Serre's bounds show.  For example, for $G = \SU(n)$ and $\rho = \Ad$, the adjoint representation, $\min_{g \in G} \chi(g)$ is asymptotically equivalent to $-\chi(1)/n^2$, see Example \ref{SUn.eg}.  In \S\ref{what} we observe that $\min_{g \in G} \RE \chi(g)$ is an algebraic number.  In \S\ref{general.sec} we provide some bounds on $\min_{g \in G} \RE \chi(g)$ for all $\chi$ and fixed $G$.

Finally, Appendix \ref{finite.sec} uses results from the main body of the paper to bound Brauer characters of representations of finite groups over fields of prime characteristic.

\subsection*{Notation and conventions} 
In most of this paper, we consider objects over a base field $F$ that is either $\R$ or $\C$.  In \S\ref{what} we consider other fields of characteristic zero, and in Remark \ref{finite1} and the appendix we will consider fields of prime characteristic.

 A real Lie group is a group object in the category of analytic real manifolds \cite[\S{III.1}]{Bou:g1}.  Consistent with the definitions of \cite[\S{IX.7}]{Bou:g7}, a representation of such a group $G$ is a continuous (hence analytic) homomorphism to the (complex) Lie group $\GL_{n,\C}$ of $n$-by-$n$ invertible complex matrices for some $n \ge 1$.

We refer to Bourbaki \cite[Ch.~VI]{Bou:g4} for general information about root systems.  We typically write $\Phi$ for the set of roots and $\Phi^+$ for the positive roots.

We also consider algebraic groups as in Milne's book \cite{Milne}.  An algebraic group over a field $F$ is a group object in the category of schemes of finite type over $F$.  Two examples are $\GL_{n,F}$, the group of $n$-by-$n$ matrices over $F$, and $U(n)_F$, the group of $n$-by-$n$ unitary matrices with entries in $F(\sqrt{-1})$.  When $F = \R$, these are also real Lie groups.  In this paper, we only consider algebraic groups that are affine, and for simplicity we write simply ``$G$ is an algebraic group'' to mean that $G$ is an affine algebraic group in the sense of Milne's book.  Furthermore, apart from the exceptions already mentioned, we will only consider the case where the field $F$ has characteristic zero, so $G$ will be linear (meaning that it is isomorphic to a Zariski-closed $F$-subscheme of $\GL_{n,F}$ in a way that is compatible with the group operations on $G$ and $\GL_{n,F}$) and smooth, see \cite[Prop.~1.41, Th.~3.23, Cor.~4.10]{Milne}.  Because of the assumption of characteristic zero and that $G$ is affine, it is the same to consider these groups from the perspective of \cite{Borel} or \cite{Sp:LAG}.
 For any commutative $F$-algebra $R$, we write $G(R)$ for the set $\Hom(\Spec R, G)$ of $R$-points of $G$.   The algebraic group $\GL_{1,R}$, i.e., the 1-dimensional torus with $R$-points $R^\times$, is typically denoted $\Gm$.  We sometimes write $\Lie(G)$ or $\mfg$ for the Lie algebra of $G$, and similarly for other algebraic groups. 
 
Mostly the symbol $T$ will denote a torus (most often a maximal torus) in a Lie or algebraic group $G$, and $W$ will denote the quotient $N_G(T)/T$.
 
 We say that an algebraic group $G$ over a field of characteristic zero is \emph{reductive} if its unipotent radical is trivial, i.e., if the identity component $G^\circ$ is reductive in the sense of \cite[\S19c]{Milne}.  We say that $G$ is \emph{semisimple} if  it is connected and its radical (the maximal closed normal solvable subgroup of $G$) is the identity.  Every semisimple group is reductive; a connected reductive group is semisimple if and only if its center is finite \cite[Prop.~19.10]{Milne}.
  We say $G$ is \emph{simple} (\cite[Def.~17.97]{Milne} says ``almost simple'') if it is semisimple and non-commutative, and every proper normal algebraic subgroup is finite.
  
A compact (real) Lie group $G$ can also be viewed as a reductive real algebraic group.  The connected components of $G$ are the same, regardless of which viewpoint one takes.  (The fact that the two notions coincide on the identity component is \cite[24.6(ii)]{Borel} and the more general statement follows, see \cite{Conrad:cpt}.)  Of course, for non-compact Lie groups this is typically false.


%
\subsection*{Acknowledgements} We are  grateful to J-P.~Serre for drawing our attention to this problem and for his insightful comments on the subject and on this text, which included suggesting the particularly nice formulations of Lemmas \ref{St.lem} and \ref{D.lem}.  Thank you to the multiple careful referees whose comments improved this text.  SG thanks John Voight for his timely advice on Magma.  RMG was partially supported by NSF grant DMS-1901595 and Simons Foundation Fellowship 609771.


\section{Serre's bounds on the trace} \label{serre.bounds}
 The following is announced as Theorem 3 in \cite{Se:trace}:

\begin{thm} \label{devissage.cor}
For any compact Lie group $K$, and any $x \in K$, one has 
\[
\Tr \Ad(x) \ge -\rank K.
\]
\end{thm}

Here $\Ad \colon K \to \GL(\Lie(K))$ is the adjoint representation, the action of $K$ on its own Lie algebra by conjugation.  This result is proved in \S\ref{devissage.sec} as a consequence of the following more detailed result, which considers both the minima and maxima of $\Tr \Ad$ at the price of restricting to simple compact Lie groups.  See \cite{FreyGriess} and \cite{Frey} for calculations of the value of $\Tr \Ad(x)$ for various $x$, in the case where $K$ is simple of type $E_8$.

Let $G$ be a simple compact Lie group.  Consider $\Out(G)$, the group obtained as the quotient of $\Aut(G)$ by the normal subgroup of inner automorphisms of $G$.  For $s \in \Out(G)$, define 
\[
\mmin(s) := \inf \Tr \Ad(u) \eand \mmax(s) := \sup \Tr \Ad(u)
\]
where $u$ ranges over those $u \in \Aut(G)$ mapping to $s$ in $\Out(G)$ (a connected component of $\Aut(G)$).  As an example, if we take $s = 1$, the elements $u$ range over $g$ in the adjoint group of $G$ and
\begin{equation} \label{max1}
\mmax(1) = \dim G,
\end{equation}
as mentioned in the introduction.

\begin{thm} \label{simple.thm}
Let $G$ be a simple compact Lie group.  For $s \in \Out(G)$:
\begin{enumerate}[(a)]
\item \label{simple.m1} If $s$ corresponds to the automorphism $-1$ of the root system of $G$, we have $\mmin(s) = -(\rank G)$.
\item Otherwise, the values of $\mmin(s)$ and $\mmax(s)$ are given by \eqref{max1} or \autoref{simple.table}.
\end{enumerate}
\end{thm}

\begin{table}[hbtp]
\centering
{\rowcolors{2}{white}{gray!25}
\begin{tabular}{lcrr} \toprule
type of $G$& order of $s$ & $\mmin(s)$&$\mmax(s)$ \\ \midrule
$A_n$ ($n \ge 1$)&1&$-1$ &$\dim G$\\
$A_n$ ($n \ge 2$)&2&$-(\rank G)$&$\begin{cases}
n&\text{if $n$ is even} \\
2+n&\text{if $n$ is odd}
\end{cases}$ \\ 
$D_n$ (even $n \ge 4$)&2&$2-n$&$2n^2 - 5n + 2$\\
$D_4$&3&$2-n = -2$&7\\
$D_n$ (odd $n \ge 5$)&1&$2-n$&$\dim G$\\ 
\ditto&2&$-(\rank G)$&$2n^2 - 5n + 2$\\ 
$E_6$&1&$-3$&$\dim G$ \\ 
\ditto&2&$-(\rank G)$&26\\  \bottomrule
\end{tabular}}
\caption{Values of $\mmin(s)$  and $\mmax(s)$.} \label{simple.table}
\end{table}

\begin{table}[hbtp]
\centering
{\rowcolors{2}{white}{gray!25}
\begin{tabular}{lrrr} \toprule
type of $G$&min&$\dim \rho$ \\ \midrule
$B_n$ ($n \ge 2$)&$1-2n$ &$2n+1$ & \\
$C_n$ ($n \ge 2$)&$\begin{cases}
1-n&\text{if $n$ is odd}\\
-1-n&\text{if $n$ is even}
\end{cases}$
 &$2n^2-n-1$ \\
$F_4$& $-6$&26 \\
$G_2$&$-2$&7 \\ \bottomrule
\end{tabular}}
\caption{The minimum of $\Tr \rho(g)$ for $g \in G$, for each simple compact $G$ with roots of different lengths and $\rho$ the irreducible representation $\rho$ with highest weight the highest short root.} \label{short.root}
\end{table}

Theorem \ref{simple.thm} expands somewhat on what is stated in \cite{Se:trace}, which only reports $\mmin(1)$.  The values of $\mmin(s)$ and $\mmax(s)$ for other $s$ were also computed by Serre in the early 2000s.  In this paper, the value of $\mmax(s)$ is computed in \S\ref{outer.sec} for $s \ne 1$.

Serre also stated, as Theorem 3$'$ in \cite{Se:trace}, the following equality.  In it, $W$ is the Weyl group if $G$ is connected.

\begin{cor} \label{thm.3p}
For $T$ a maximal torus in a compact Lie group $G$ and $W := N_G(T)/T$, 
\[
\inf_{g \in G} \Tr \Ad(g) = \inf_{w \in W} \Tr(w\vert_{\Lie(T)}).
\]
\end{cor}

We provide a proof in \S\ref{finite.proof}.  

In addition to the results we prove here, Serre has recently published in \cite{Se:zeros} his proofs of the other results announced in his 2004 note, so now proofs of all of the announced results are available.

\subsection*{The highest short root representation} 
In the case of a simple group with roots of different lengths, there is an irreducible representation which is almost as interesting as the adjoint one, namely the one whose highest weight is the highest short root.  We bound also $\Tr \rho$ for such a representation $\rho$.  Note that for these groups, there are no outer automorphisms, so it is just a question of bounding $\Tr \rho$ on the simple group itself.  Also, the maximum of $\Tr \rho$ is $\dim \rho$.

\begin{prop} \label{short.prop}
Let $G$ be a simple compact Lie group with roots of different lengths and put $\rho$ for the irreducible representation whose highest weight is the highest short root.  The extreme values of $\Tr \rho(g)$ for $g \in G$ are given in \autoref{short.root}.
\end{prop}

The values for $G_2$ in  \autoref{short.root} were previously observed in \cite[5.5]{Katz:G2} and \cite[eq.~(37)]{Kaiser}.  Still for type $G_2$, see \cite[pp.~270, 271]{Griess:G2} for calculations of the value of $\Tr \rho(g)$ for various $g$ of finite order.

\section{Preliminary observations}

\begin{eg} \label{minus1}
Let $G$ be a simple compact Lie group with maximal torus $T$.  There is a $\phi \in \Aut(G)$ that normalizes $T$ and acts as $\phi(t) = t^{-1}$ for $t \in T$ \cite[\S{IX.4.10}, Prop.~17]{Bou:g7}.  (It is even unique up to conjugation by an element of $T$.)   Then $\phi$ interchanges the positive and negative roots, so the only nonzero diagonal entries of $\Ad(\phi)$ relative to a root space decomposition of $\Lie(G)$ come from its action on the torus and we have $\Tr \Ad(\phi) = -(\rank G)$.    In the notation of Th.~\ref{simple.thm}\ref{simple.m1}, this proves that $\mmin(s) \le -(\rank G)$.
\end{eg}

\begin{eg}[$\SU(n)$] \label{SUn.eg}
Consider the case of $G = \SU(n)$, which is simple compact of type $A_{n-1}$.  The complexification of $G$ is $\SL_n(\C)$ and we may identify the adjoint representation of $G_\C$ with $(\C^n \otimes (\C^n)^*)/ \C$.  For $g \in \SL_n(\C)$ we have $\Tr \Ad(g) = \Tr(g) \Tr(g^{-1}) - 1$.  

An element  $g \in \SL_n(\C)$ coming from $\SU(n)$ has eigenvalues $z_1, \ldots, z_n$ such that $\bar{z}_j = z_j^{-1}$ for all $j$, so setting $z := \sum z_j$ we have $1 + \Tr \Ad(g) = z \bar{z} \ge 0$, so $\mmin(1) \ge -1$.  On the other hand, if $g \in \SU(n)$ is diagonal with entries pairs $i, -i$ or triples $1, \zeta, \zeta^2$ for $\zeta = e^{2\pi i /3}$, then $\Tr(g) = \Tr(g^{-1}) = 0$ and $\Tr \Ad(g) = -1$, which shows the claimed minimum is attained.  This proves that $\mmin(1) = -1$ for $G$ of type $A$ as indicated in \autoref{simple.table}.
\end{eg}

In all our approaches, we will lean on the following observation, which reduces to a classic result by Steinberg:

\begin{lem}  \label{St.lem}
Let $f$ be a polynomial in the fundamental characters $p_1, \ldots,p_r$ of a semisimple complex algebraic group  $G$.
Let $g$ be a semisimple element of $G$ that is a critical point of $f$.   Then at least one of the following holds:
\begin{enumerate}[(i)]
\item \label{St.nreg} $g$ is not regular.
\item \label{St.grad} For all $i$, the partial derivative $\partial_{p_i} f$ of $f$ with respect to $p_i$ is zero at $g$.
\end{enumerate}
\end{lem}

\begin{proof}
The element $g$ is contained in a maximal torus $T$ of $G$.  Put $W$ for the Weyl group of $G$ with respect to $T$.  Write $f$ as a composition
$T \to T/W \to \C$.  Suppose $g$ is regular.  Then the Jacobian of the first map is invertible at $g$ by Steinberg's Theorem as in \cite{St:inv} or \cite[6.1, 6.3]{St:reg}.  Since $g$ is a critical point of $f$, it follows that the gradient of the map $T/W \to \C$ must vanish at the image of $g$, i.e., \ref{St.grad}.
\end{proof}

We also need the following small elaboration of results of Steinberg.  We provide a proof because we do not know a reference for it in the literature.  Every algebraic group $G$ acts on itself by conjugation, and we write $\ClG$ for the (categorical) quotient of $G$ by this action (in order to avoid the notation $G/G$).  As a scheme, for $G$ a reductive algebraic group over $F$, $\ClG$ is $\Spec F[G]^G$ and $F[G]^G$ is finitely generated.

\begin{prop} \label{torus.prod}
Let $T$ be a maximal torus in a connected reductive algebraic group $G$ over a field $F$ of characteristic zero, and put $W := N_G(T)/T$.  
\begin{enumerate}[(a)]
\item \label{St.inv.iso} $\ClG \cong T/W$, i.e., the natural map $\Spec F[G] \to \Spec F[T]$ gives an isomorphism $\Spec F[G]^G \cong \Spec F[T]^W$.
\item \label{St.inv.prod} If the torus $Z(G)^\ic$ is split and the derived subgroup $G'$ of $G$ is simply connected, then 
	\begin{enumerate}[(i)]
	\item \label{St.inv.prod1} $\ClG$ is the product of an affine space and a torus, and
	\item \label{St.inv.prod2} if $G'$ is also split, there is an open
          embedding $\ClG \hookrightarrow \A^{\rank G}$ such that the
          coordinates are irreducible characters.
	\end{enumerate}
\end{enumerate}
\end{prop}

\begin{proof}
To verify that the natural map in \ref{St.inv.iso} is an isomorphism, it suffices to do so when $F$ is algebraically closed.  If $G$ is semisimple, then this is a result of Steinberg, \cite[6.4]{St:reg}.  The proof is a reduction to that case.

Let $T'$ be a maximal torus in the derived subgroup $G'$ of $G$.  The identity component $Z(G)^\ic$ of the center of $G$ is a torus and $T := T' \cdot Z(G)^\ic$ is a maximal torus of $G$.  By conjugacy of maximal tori in $G$, we may assume that the given torus $T$ is this $T$.  Set $H := T' \cap Z(G)^\ic$, $T^+  := T' \times Z(G)^\ic$, and $G^+ := G' \times Z(G)^\ic$.  There is a natural diagram
\[
\begin{CD}
1 @>>> H @>>> T^+ @>>> T @>>> 1\\
@. @| @VVV @VVV @. \\
1 @>>> H @>>> G^+ @>>> G @>>> 1
\end{CD}
\]
which commutes.  

Certainly, $N_{G^+}(T^+)/T^+ \cong N_{G'}(T')/T' \cong W$.  We obtain another commutative diagram
\[
\begin{CD}
T^+/W @>>> T/W \\
@V{\cong}VV @VVV \\
\class(G^+) @>>> \ClG
\end{CD}
\]
For the left vertical arrow, $T^+ / W = T/W \times Z(G)^\ic$ and $\class(G^+) = \class{G'} \times Z(G)^\ic$, and these are isomorphic because $G'$ is semisimple.  That is, the left vertical arrow is an isomorphism.  The top and bottom arrows are a quotient map by $H$, on which $W$ and $G^+$ both act trivially, so the right arrow is also an isomorphism, proving \ref{St.inv.iso}.

For \ref{St.inv.prod}, suppose first that $G'$ is split.  Since $G'$ is
simply connected, $\class(G')$ is an affine space \cite[Th.~6.1(b)]{St:reg}
where the coordinates are fundamental characters.  Thus $\class(G^+)$ is
manifestly a product of an affine space and a torus.

If $G'$ is not split, the previous case descends because the action of
$\Aut(G')$ on $\class(G')$ factors through $\Out(G')$, which permutes the
fundamental characters according to the diagram automorphisms.  In
particular, the induced class in the Galois cohomology set $H^1(F,
\Aut(\class(G')))$ is in the image of $H^1(F, \GL_{\rank G, F})$, which is
a singleton.  So again $\class(G')$ is  an affine space and $\class(G^+)$
is a product of an affine space and a torus, with the affine space
naturally identified with the space of sections of a vector space $V_{G'}$.

It thus remains to consider the action of $H$.  Since
$Z(G)^\ic$ is split, $\Hom(H,\Gm)$ is a finite abelian group with trivial
action by the Galois group, and thus the action of $H$ induces an
$F$-linear grading of $V_{G'}$ by $\Hom(H,\Gm)$.  In addition, the map
$\Hom(Z(G)^\ic,\Gm)\to \Hom(H,\Gm)$ is surjective, and thus we can rescale
each homogeneous component by a character of $Z(G)^\ic$ to get an
$H$-invariant affine factor.  We thus get an isomorphism
\[
\class(G)\cong (\A^{r'}\times Z(G)^\ic)/H
\cong \A^{r'}\times (Z(G)^\ic/H).
\]
Since $Z(G)^\ic$ is a torus, the quotient $Z(G)^\ic/H$ is also a torus (it is also the abelianization $G^\text{ab}$ of $G$),
as required.

When $G'$ is also split, $V_{G'}$ is the span of fundamental characters,
all of which are homogeneous for the grading by $\Hom(H,\Gm)$, and thus the
procedure of the previous paragraph multiplies the fundamental characters
by characters of $Z(G)^\ic$ to obtain new characters that descend to $G$.
We thus get an isomorphism of the form
\[
\class(G)\cong \A^{r'}\times G^{\text{ab}}\cong (\A^1)^{r'}\times
\Gm^{r-r'}
\]
in which the maps to $\A^1$ are irreducible characters of $G$ and the maps
to $\Gm$ are irreducible characters of $G^{\text{ab}}$, thus also of $G$.
Embedding each $\Gm$ in $\A^1$ in the natural way gives \ref{St.inv.prod2}.
\end{proof}

\begin{eg}
  The group $G = U(2)$ has simply connected derived subgroup $G' =
  \SU(2)$, with $\class(G')\cong \A^1$.  However, $Z(G)^\ic$ is a rank $1$
  torus that is not split and we claim that $\ClG$ is not isomorphic to the
  product of an affine space and a torus.

  Over $\C$, the invariant ring is $\C[G]^G \cong \C[e_1,e_2,e_2^{-1}]$
  with $e_1=z_1+z_2$ and $e_2=z_1z_2$, and complex conjugation acts on this
  ring by conjugating scalars and taking $(e_1,e_2)\mapsto
  (e_1/e_2,1/e_2)$.  If the corresponding real subring were the coordinate
  ring of a product $\A^1\times T$, then the map to $T$ would be given by a
  unit in the complex coordinate ring, so is given by (without loss of
  generality) $e_2$, making the torus have real locus $S^1$.  The affine
  coordinate must then correspond to a generator of the complex ring over
  $\C[e_2,e_2^{-1}]$, and thus must have the form $c_1 e_2^k e_1 +
  c_0(e_2)$ where $c_1\in \C^\times$, $c_0\in \C[e_2,e_2^{-1}]$, and $k\in
  \Z$.  The affine coordinate must, of course, also be invariant under the
  conjugation, and thus we find $c_1 e_2^k = \bar{c_1} e_2^{-k-1}$, so that
  $2k+1=0$, giving a contradiction.  (The basic issue in this case is that
  $\Hom(T,\Gm)\to \Hom(K,\Gm)$ is split as a map of sets but not as a map
  of sets with Galois action.)

   One can also use this to produce an example in which $G'$ is split, for
   the simple reason that the central product $\SL_2 \circ\ U(1)$ is an inner
   twist of $U(2)$, and thus there is an isomorphism
   \[
   \class(U(2))\cong \class(\SL_2 \circ\ U(1))
   \]
   of schemes over $\R$.
\end{eg}

\part{Approach via invariant derivations} \label{part.Wders}

\section{Invariant derivations and critical points} \label{Wders.1}

Let $G$ be a compact Lie group, and suppose we wish to optimize a finite
linear combination $f$ of irreducible characters on $G$.  There is an
obvious (and obviously infeasible) approach: the optimum will be attained
at a critical point of $f$ (defined locally on $G$ in the usual way), and
thus one can reduce to finding those critical points.  Of course, as
stated, this is entirely analytic and thus it is unclear whether we can
actually find the critical points; this becomes worse once we realize that
$G$-invariance ensures that the critical locus is a union of $G$-orbits,
almost all of which will be positive-dimensional.

Instead we will
look for critical points on the {\em algebraic} complex group $G(\C)$.
This at least in principle allows one to compute the ideal cutting out the
critical subscheme, at which point standard Gr\"obner basis techniques let
one determine the components of the critical locus.  In practice, however,
this is still not particularly feasible: this rephrases the problem as a
constrained optimization problem with a typically complicated ideal of
constraints.  Moreover, we still must contend with the fact that the output
will be a collection of $G(\C)$-invariant subschemes of $G(\C)$ and each
factor will itself be complicated, thus the Gr\"obner basis
calculations will be quite difficult.

One reduction of the problem is immediate: In order to prove bounds like the ones in \S\ref{serre.bounds}, which apply to compact and connected groups,
it is only necessary to consider critical points in $G(\C)$ that come from $G(\R)$, and in particular are semisimple.
(We remark that one can reduce the general problem of finding values of $f$ at critical points on $G(\C)$ to examining semisimple critical points, but we will omit the details of this reduction because we do not need it for the results in this paper.)

\begin{lem}   \label{torus.reduce}
Let $T$ be a maximal torus of $G$, and let $f$ be a polynomial function on $G$.  Suppose that $T$ is central in $G$ or $f$ is $G$-invariant.  Then $t\in T(\C)$ is a critical
point for $f$ if and only if it is a critical point for $f|_T$.
\end{lem}

\begin{proof} 
The ``only if'' direction is trivial.  So suppose $t$ is a critical point for $f|_T$.  
The identity component $G^\ic$ of $G$ is a union of conjugates of $T$, and thus in particular the conjugates of $\mft$
span $\mfg$.  For $x \in \mfg$, there exists $g_1, \ldots, g_n \in G$ and $y_1, \ldots, y_n \in \mft$ such that $\sum_i (\Ad g_i)y_i = x$ and we have 
\[
xf(t) = \sum_i ((\Ad g_i) y_i)f(t) = \sum_i y_i(g_i^{-1}f)(t) = \sum_i y_i f(g_itg_i^{-1}).
\]
By hypothesis, the summand in the last term equals  $y_i f(t)$, which is zero because $t$ is a critical point for $f|_T$.
\end{proof}

When $G$ is connected, Lemma \ref{torus.reduce} reduces the problem considerably: now we need only work on the
torus $T$.  

\subsection*{Derivations and group actions}
Here is another elementary result on critical points for functions invariant under the action of a finite group.  Let $H$ be a finite group acting on an affine scheme $X$ over a field $F$.  An element $h \in H$ acts on a function $f \in F[X]$ via ${^h f}(x) = f(h^{-1} x)$.  A \emph{derivation} on $X$ is an $F$-linear derivation $F[X] \to F[X]$, specializing a definition in \stacks{00RN}.  An element $h \in H$ acts on a derivation $D$ via $(({^h D})f)(x) = (D({^{h^{-1}}f}))(h^{-1}x)$.

\begin{prop} \label{Winv.crit}
Suppose $H$ and $X$ are as in the preceding paragraph, and $F$ has characteristic zero.  Let $f$ be an $H$-invariant function on $X$ and $x \in X(F)$.  Then $x$ is a critical point for $f$ if and only if $(Df)(x) = 0$ for every $H$-invariant derivation $D$.  
\end{prop}

\begin{proof}
$x$ is a critical point for $f$ if and only if $(Df)(x) = 0$ for every derivation $D$ (not necessarily $H$-invariant), so the ``only if'' direction is trivial and we assume $(Df)(x) = 0$ for every $H$-invariant $D$.  Let $D_0$ be any derivation.  Since $H$ is finite, there is a $\pi \in F[X]$ such that $\pi(x) = 1$ and $\pi$ vanishes on the other elements of the $H$-orbit of $x$.  Put $H_x$ for the centralizer of $x$ in $H$ and define
\[
D_1 := \frac{1}{\left| H_x \right|} \sum_{h \in H} {^h(\pi D_0)}.
\]
Since $f$ is $H$-invariant, we have 
\[
\left[ {^h(\pi D_0)} \right] (f)(x) = \left[ (\pi D_0)(f) \right] (h^{-1} x) = \pi(h^{-1} x) \left[ (D_0 f)(h^{-1}x) \right].
\]
For $h \in H_x$, this is $(D_0 f)(x)$, whereas for $h \not\in H_x$, it is zero.  It follows that $(D_1 f)(x) = (D_0 f)(x)$.

On the other hand, $D_1$ is an $H$-invariant derivation, so $(D_1 f)(x) = 0$ by hypothesis.
\end{proof}

In this section, we have reduced the problem of finding critical points of a $G$-invariant function $f$ on a connected compact Lie group to finding the vanishing locus of $Df$ on a maximal torus $T$ of $G$ as $D$ varies over $W$-invariant derivations on $T$ for $W := N_G(T)/T$.

\section{Derivations on \texorpdfstring{$T/W$}{T/W} versus \texorpdfstring{$W$}{W}-invariant derivations on \texorpdfstring{$T$}{T}} \label{Wders.2}

In this section, we work in the context of a connected reductive algebraic 
group $G$ over $\C$ with a maximal torus $T$ and Weyl group $W := N_G(T)/T$.
We relate the $W$-invariant derivations on $T$ --- which we denote by $(\Der \C[T])^W$ --- to derivations on $T/W$, i.e., $\Der \C[T/W]$.

Any $W$-invariant derivation on $T$ takes $W$-invariant functions to
$W$-invariant functions, so induces a derivation on $T/W$, giving us a
natural map $(\Der \C[T])^W\to \Der \C[T/W]$.  Every derivation on $\C[T/W]$ extends uniquely to its function field
$\C(T/W)$, and hence to that field's finite separable extension $\C(T)$.  In particular, a derivation on $T$ is determined by its values on $\C[T/W]$,
so this map is injective.  This lets us identify $(\Der \C[T])^W$ with a
submodule of $\Der \C[T/W]$.

\begin{eg}[$\SL_2$] \label{A1.eg}
In the special case where $G = \SL_2$, we may describe the relationship in an elementary manner.  Imagine the torus as consisting of diagonal matrices with coordinate functions $s, s^{-1}$.  We have $\C[T] \cong \C[s, s^{-1}]$, where the nontrivial element of the Weyl group $\Z/2$ acts by $s \mapsto s^{-1}$.  In particular, $\C[T]^W \cong \C[t]$ for $t := s + s^{-1}$, the trace of an element of $T$.

It is an exercise to verify that the Weyl-invariant derivations on $T$ form a free $\C[T]^W$-module generated by $(s^2 - 1) \partial_s$.  The derivations of $\C[T]^W$, on the other hand, are evidently generated by $\partial_t$.  The usual chain rule shows the two are related by 
\[
(s^2 - 1) \partial_s = (t^2 - 4) \partial_t.
\]
For computations, we will prefer the version on the right, which is expressed as a derivation on $T/W$.  Note that the displayed equality shows that
\[
(\Der \C[T/W]) / (\Der \C[T])^W \cong \C[t]/(t^2 - 4)
\]
 and in particular that $(\Der \C[T])^W$ is a proper submodule of $\Der \C[T/W]$.
\end{eg}

It turns out that when $T/W$ is sufficiently nice (e.g., when the derived subgroup $G'$ is
simply connected), not only can we show that the relevant submodule of
$\Der \C[T/W]$ is a free submodule, but the proof gives us a reasonably
natural way to construct a basis.  (A similar construction appears in
\S{B.2} of \cite{OrlikTerao}, although that material concerns the linear
action of $W$, the quotient of a euclidean space by the reflection
representation of $W$.)

\subsection*{Interlude: differential operators}

We take a moment to review some background on differential operators on an affine variety $X$ over $\C$.

\begin{defn}
A differential operator of order 0 is $\mu_f \colon \C[X] \to \C[X]$ given by multiplication by an element $f \in \C[X]$.  A differential operator $\mcD$ of order $k + 1$ is a $\C$-linear map $\C[X] \to \C[X]$ such that, for each $g \in \C[X]$, the map $f \mapsto \mcD(gf) - g\mcD(f)$ is a differential operator of order $k$.  
\end{defn}

The collection of differential operators of order $k$ is denoted $\Diff^k_{\C[X]/\C}(\C[X], \C[X])$ in \stacks{09CH}, and it is a $\C$-module.
A differential operator of order $k$ is a differential operator of order $K$ for all $K \ge k$.  A differential operator of order 1 is a derivation if and only if it annihilates 1  \stacks{09CM}.
 
\begin{lem} \label{diffop}
Let $\mcD$ be a second-order differential operator on an affine variety $X$ over $\C$.  Then for $f, g \in \C[X]$:
\begin{enumerate}[(a)]
\item \label{diffop.der} $[\mu_f, \mcD] + \mu_{\mcD(f)} - \mu_{\mcD(1)f}$ is a derivation on $X$.
\item \label{diffop.bider} The pairing $(f, g) \mapsto \mcD(fg) - f \mcD(g) - \mcD(f)g + fg \mcD(1)$ is symmetric and is a derivation on $X$ in each slot.  It is identically zero if $\mcD$ is a first-order differential operator.
\end{enumerate}
\end{lem}

\begin{proof}
For \ref{diffop.der}, we note that the first term is a first-order differential operator (a fact already implicit in \stacks{0G35}).  The other two terms are order 0, so the total expression is a first-order differential operator.  It annihilates 1 because $\mu_f \mcD(1) = \mu_{\mcD(1)f}(1)$ and $\mcD \mu_f(1) = \mcD(f)1$, so it is a derivation.

For \ref{diffop.bider}, the pairing is evidently symmetric and $\C$-linear in each slot.  If $\mcD$ is first order, then there is an $f' \in \C[X]$ such that $\mcD(fg) - f\mcD(g) = f' g$ for all $g$.  Then $\mcD(f)g - fg\mcD(1) = f'g$ as well and the pairing is zero.  Finally, the the image of $(f, g)$ is $-1$ times the derivation in \ref{diffop.der}, so it is a derivation in $g$.
\end{proof}

The pairing in \ref{diffop.bider} is called the \emph{symmetric
biderivation} associated with $\mcD$.  The construction is analogous to the
one for a polynomial function $q$ of degree at most 2 on a vector space
$V$: One defines a pairing $(v, v') \mapsto q(v + v') - q(v) - q(v') +
q(1)$, which is symmetric and linear in each slot and vanishes if and only
if $q$ has degree at most 1.

Note that a symmetric biderivation $\delta$ on a commutative
$\C$-algebra $R$ is
equivalent to a symmetric bilinear form on the module of K\"ahler differentials $\Omega_R$, and thus in
particular induces a map
\[
\Omega_R\to \text{Der}(R)
\]
taking $f\,dg$ to the derivation $h\mapsto f\delta(g,h)$.

\subsection*{Constructing invariant derivations from invariant functions}
  We now return to the situation of $G$, $T$, and $W$ from the beginning of this section.
Let $A$ be a nonsingular $W$-invariant quadratic form on $\mft^*$.
We may view this as an element of $(S^2\mft)^W$, i.e., a Weyl-invariant element of the symmetric square.  Since
elements of $\mft$ are naturally translation-invariant derivations,
$A$ corresponds naturally to a second-order differential operator $\mcD_A$
on $\C[T]$, invariant under $W$ and under translation by $T$.  

The space
$\mft^*$ may be viewed as the space of translation-invariant
differentials, and thus any homomorphism $z\colon T\to \Gm$ induces an element
$dz/z$ of $\mft^*$.  In particular, we get a basis of $\mft^*$ from any isomorphism $T\cong \Gm^n$ by taking the logarithmic
differentials of the coordinates.

Since $\mcD_A$ is $W$-invariant, its value on a $W$-invariant function
is $W$-invariant, and thus the symmetric biderivation associated to $\mcD_A$ induces a $\C[T/W]$-valued biderivation $M_A$ on $T/W$.

\begin{prop}
  If $T/W$ is smooth (e.g., if the derived subgroup of $G$ is simply connected), then the image of the map
\[
\Omega_{\C[T/W]}\to \Der \C[T/W]
\]
induced by $M_A$ is $(\Der \C[T])^W$.
\end{prop}

\begin{proof}
Since $A$ is nonsingular, it induces an isomorphism
\[
\Omega_{\C[T]}\cong \Der \C[T],
\]
which by $W$-equivariance restricts to an isomorphism
\[
\Omega_{\C[T]}^W\cong (\Der \C[T])^W.
\]
The map
\[
\Omega_{\C[T/W]}\to \Der \C[T/W]
\]
induced by $M_A$ is then the composition
\[
\Omega_{\C[T/W]}\to \Omega_{\C[T]}^W\cong (\Der \C[T])^W\hookrightarrow \Der \C[T/W]
\]
in which the maps on either side are natural.  Moreover, by \cite[Cor.\ to
  Th.~II.3]{Knighten}, the first natural map is an isomorphism, and thus
the image of the map induced by $M_A$ is $(\Der \C[T])^W$ as
required.
\end{proof}

\begin{cor} \label{der.free}
  If the derived subgroup $G'$ of $G$ is simply connected, then $(\Der \C[T])^W$ is a free module over $\C[T]^W$.
\end{cor}

\begin{proof}
  Since $G'$ is simply connected, by Prop.~\ref{torus.prod}\ref{St.inv.prod}\ref{St.inv.prod2} we may identify $T/W$ with an open subscheme of $\A^r$ where $r=\dim T$.  Since open embeddings are \'etale, the sheaf of K\"ahler differentials on $T/W$ is the pullback of the sheaf of K\"ahler differentials on $\A^r$. Since the latter is free (with basis $df_1,\ldots,df_r$ where $f_1,\dots,f_r$ are the coordinate functions), so is $\Omega_{\C[T/W]}$.
\end{proof}

From the proof of the corollary, we see that we may associate a basis of
the space of $W$-invariant derivations to any choice of coordinates on
$T/W$, i.e., to any choice of open embedding $T/W\to \A^n$.

Given such coordinates $f_1,\dots,f_n$, the basis consists of the derivations
\begin{equation} \label{Di.def}
D_i\:f\mapsto M_A(f_i,f) = \sum_j M_A(f_i,f_j) \partial_{f_j} f
\end{equation}
and thus is determined by the symmetric matrix $M$ whose $(i,j)$-entry is
$M_A(f_i,f_j)$.  If we also have an isomorphism $(z_1,\dots,z_n)\:T
\xrightarrow{\sim} \Gm^n$, then this matrix factors via the chain rule:
\[
M_A(f_i,f_j) = \sum_{k,l} A_{kl}
(z_k\partial_{z_k}f_i)(z_l\partial_{z_l}f_j)
\]
where 
\[
A_{kl} =
\textstyle\frac{1}{2}\left(A(\textstyle\frac{dz_k}{z_k}+\frac{dz_l}{z_l})-A(\frac{dz_k}{z_k})-A(\frac{dz_l}{z_l}) \right)
\in\C
\]
is the $(k,l)$-entry of the symmetric matrix corresponding to the quadratic form $A$ in the
basis corresponding to $\vec{z}$.  In other words, the matrix $M$ (which depends on the choice of $A$ and the coordinates $f_1, \ldots, f_n$) is 
\[
M = J^\trans A J, 
\]
where $J$ is the Jacobian matrix with $(i, j)$-entry $z_i\partial_{z_i}f_j$.

\begin{eg} \label{simple.A}
Suppose $G$ is simple (and simply connected).  There is a unique Weyl-invariant symmetric bilinear form $b$ on the root lattice $Q$ such that $b(\alpha, \alpha) = 2$ for each short root alpha; it is positive definite and integer-valued.  Put $B$ for the matrix of $b$ with respect to some pinning of $G$ with respect to $T$.  Since $G$ is simply connected, the space $\mft^*$ is identified with the complexification $P \otimes \C$ of the weight lattice $P$ and one can take $f_1, \ldots, f_n$ to be the fundamental characters.  We take $A$ to be some positive integer multiple of $B$; any such choice will determine a real inner product on $P \otimes \R$.  
In the the Magma code for this paper and in the $G_2$, $F_4$, and $E_8$ examples in this part of the paper, we use $A := | P/Q | \, B$.
(Magma's \texttt{CoxeterForm}, when passed a root datum with \texttt{IsogenyClass := "SC"}, returns a matrix for $\frac12 B$ with respect to a basis consisting of fundamental weights.)  Choosing a different scalar factor will scale the matrix $A$, but will leave the critical point computations otherwise unchanged, see Corollaries \ref{fi.crit} and \ref{corner.crit}.

For the case of $\SL_2$ from Example \ref{A1.eg}, $A = (2)$, $J = (z_1)$, and $M = (2z_1^2)$.
\end{eg}

\begin{rmk}
Although it is natural to use the fundamental characters as the generators for the ring of invariants, there are other choices.
For $E_6$, $E_7$, and $E_8$, the matrix $M$ can be significantly simplified by using a basis coming from the identification of $T/W$ with a moduli space of del Pezzo surfaces, see \cite{EOR}.
\end{rmk}

\begin{eg}[$\GL_2$]  \label{GL2.eg2}
  Suppose now that $G = \GL_2$ and take $T$ to be the diagonal matrices,
  with coordinate functions $z_1, z_2$.  Coordinates on $\ClG$ are provided
  by $f_1 := z_1 + z_2$ and $f_2 := z_1 z_2$, so
  \[
J = \left( \begin{smallmatrix} 
z_1 & z_1 z_2 \\
z_2 & z_1 z_2
\end{smallmatrix} \right).
  \]
  There is a 2-dimensional space of choices for $A$, with the simplest one
  corresponding to the pairing $(dz_i/z_i,dz_j/z_j)=\delta_{ij}$ (Kronecker delta).  In this
  form, $M_A$ has matrix
  \[
M =  J^\trans J = \left(
\begin{smallmatrix}
z_1^2 + z_2^2 & z_1^2 z_2 + z_1 z_2^2 \\
z_1^2 z_2 + z_1 z_2^2 & 2z_1^2 z_2^2
\end{smallmatrix} \right) =
\left( \begin{smallmatrix}
f_1^2 - 2f_2 & f_1 f_2 \\
f_1 f_2 & 2f_2^2
\end{smallmatrix} \right),
\]
so that the module of invariant derivations has basis
\[
D_1=(f_1^2-2f_2)\partial_{f_1} + (f_1f_2)\partial_{f_2}, \quad
D_2=(f_1f_2)\partial_{f_1}+2f_2^2\partial_{f_2}.
\]
\end{eg}

For semisimple $g \in G(\C)$, it makes sense to write $M(g)$ for the matrix with complex entries obtained by evaluating each entry of $M$ (which belongs to $\C[T/W]$) at the image of $g$ in $(T/W)(\C)$.

\begin{cor} \label{fi.crit}
Let $t$ be an element of $T(\C)$.  The $i$-th column (equivalently, row) of $M(t)$ vanishes if and only if $t$ is a critical point for  $f_i$ on $G$.
\end{cor}

\begin{proof}
Since $M$ is symmetric, it suffices to consider only the claim involving a
column.  The $i$-th column vanishes if and only if $D_j(f_i)(t) = 0$ for
all $j$; since we have seen that the $D_j$ form a basis of the submodule of
$\Der \C[T/W]$ corresponding to $W$-invariant derivations, this is
equivalent to having $D(f_i)(t)=0$ for every $W$-invariant derivation $D$
on $T$.  This is equivalent by Prop.~\ref{Winv.crit} to $t$ being a
critical point for $f_i\vert_T$, which in turn is equivalent to being a
critical point for $f_i$ on $G$ (Lemma \ref{torus.reduce}).
\end{proof}

\begin{defn}
A \emph{corner} of $G$ is an element of $(T/W)(\C)$ where $M$ vanishes.  Note that since $df_1,\ldots, df_n$ are a basis of the module of differentials, $M$ vanishes at $t$ if and only if $M_A$ vanishes at $t$, and thus this notion is basis independent.  
\end{defn}

Since $(T/W)(\C)$ is in bijection with conjugacy classes of semisimple elements of $G(\C)$, one could equally well refer to the conjugacy classes corresponding to corners as corners themselves.

\begin{cor} \label{corner.crit}
Let $G$ be a semisimple simply connected algebraic group.  The corners are the images of those $t \in T(\C)$ that are critical points for every $G$-invariant function.
\end{cor}

\begin{proof}
If $t$ is a critical point for every $G$-invariant function, it is a critical point for $f_1,\ldots, f_n$ and thus $D_i(f_j)(t) = 0$ for all $i$ and $j$, i.e., $M$ vanishes at $t$.

Conversely, if $M$ vanishes at $t$, then for any $G$-invariant
function $f$,
\[
D_i(f)(t)
=
\sum_j (\partial_{f_j}f)(t) M_A(f_i,f_j)(t)
\]
is zero for each $i$.  Since the $D_i$ are a basis of the module of
$W$-invariant derivations, it follows that $D(f)(t)=0$ for any
$W$-invariant derivation $D$ and thus that $t$ is a critical point for $f$.
\end{proof}

We will say more about corners in \S\ref{corners.sec}.

\begin{rmk} \label{finite1}
Although one no longer has an application to minimization or maximization,
the module of invariant derivations is still a natural object to consider
over a field of prime characteristic.  There are two places in which the above
construction may fail.  For one, it may not be possible to choose a {\em nonsingular}
invariant form $A$.  In that case, $\det(A)=0$, so $\det(M)=0$, and the rows of $M$ do not
even generate the space of \emph{rational} invariant derivations.  

For another,
the module of invariant differentials may differ
from the module of differentials on the quotient.
Per \cite{Knighten}, the latter can only occur when the quotient map is wildly
ramified in codimension 1, or in other words when there are hypersurfaces
fixed pointwise by elements of order equal to the characteristic of the field, which in our case forces characteristic 2.

To illustrate this issue, consider $\SL_2$ over the algebraic closure of $\F_2$, with invariant $f_1=z_1+1/z_1$.  We could take $A$ to be the
matrix $(1)$ (instead of the more natural $(2)$), so nonsingularity is not
an issue, and thus compute $M = (f_1^2)$.  However, the derivation
$f_1^2\partial_{f_1}$ does not generate the module of invariant
derivations: the tangent space to the center of $\SL_2$ induces an
invariant derivation $z_1\partial_{z_1} = f_1\partial_{f_1}$.  (Similarly,
the pullback of $df_1$ does not generate the module of invariant
differentials.)

This leads to a natural question: Under what circumstances is the module of invariant derivations
still free over the ring of invariants in general?   We computed spanning sets\footnote{See the Magma program \texttt{bad-char} in the git repository.} for the modules of invariant derivations for $\SL_5$ in characteristic 5 and $\SL_6$ in characteristics 2 and 3; in each case we have found that the module is not free at the identity.  It thus seems likely that the module is rarely free when the characteristic is not ``very good'', with only sporadic exceptions (e.g., $\SL_2$ and $\SL_4$ in characteristic 2 and $\SL_3$ in characteristic 3).
\end{rmk}

\subsection*{Alternate construction of the invariant derivations}
There is an alternate approach to construct derivations.  There is another
natural way to associate an operator $C_A$ on $\C[T]^W\cong \C[G]^G$ to
$A$, in such a way as to be diagonal on the characters of irreducible representations.
Define
\begin{equation} \label{C.altdef}
C_A := \mcD_A + \frac{1}{2}\Delta^{-1} M_A(\Delta,\emptyslot)
\end{equation}
where
\[
\Delta:= \prod_{\alpha\in \Phi^+} (e^{\alpha}-2+e^{-\alpha}) \quad \in \C[T]^W.
\]
Here we are writing an unadorned symbol such as $\alpha$ for roots and weights (viewed as elements of $\mft^*$, which is canonically isomorphic to $\Hom(T, \Gm) \otimes \C)$ and $e^\alpha$ for the corresponding element of $\Hom(T, \Gm)$.

We now describe $C_A$ in a special case.
Recall the Weyl denominator 
\[
\delta := e^\rho\prod_{\alpha\in \Phi^+} (1-e^{-\alpha})
\]
where $\rho$ is half the sum of the positive roots.  If $\rho$ corresponds to an element of $\Hom(T, \Gm)$, then $\delta$ belongs to $\C[T]$ and is anti-invariant under the Weyl group \cite[\S{VI.3.3}, Prop.~2(i)]{Bou:g4}, and moreover we have
\[
\delta^2=e^{2\rho} \prod_{\alpha\in \Phi^+}(1-e^{-\alpha})^2
= \Delta.
\]

\begin{lem}
If the half sum of the positive roots, $\rho$, corresponds to an element of $\Hom(T, \Gm)$, then 
\begin{equation} \label{C.scdef}
C_A = \delta^{-1} \mcD_A \delta - \delta^{-1}\mcD_A(\delta).
\end{equation}
\end{lem}

\begin{proof}
Since $M_A$ is a derivation in the first argument, we have
\[
\delta^{-k} M_A(\delta^k,\emptyslot)
=
k \delta^{-1} M_A(\delta,\emptyslot)
\]
for any integer $k$.   Combining this with the definition \eqref{C.altdef} of $C_A$ gives
\[
C_A = \mcD_A + \delta^{-1} M_A(\delta, \emptyslot).
\]
Since $\mcD_A$ is a homogeneous quadratic polynomial in derivations (which each kill 1), we have $\mcD_A(1)=0$.  Combining this with the expression for $M_A$ as a biderivation arising from $\mcD_A$ as in Lemma \ref{diffop}\ref{diffop.bider}, we find that $C_A$ equals \eqref{C.scdef}.
\end{proof}

The inspiration for the definition \eqref{C.altdef} came from \eqref{C.scdef}.
By the Weyl character formula, multiplying an irreducible character by $\delta$ gives
the antisymmetrization of a monomial, and thus gives an eigenfunction for
$\mcD_A$.  In other words, $\delta^{-1}\mcD_A\delta$ is diagonal on the
irreducible characters, which suggested  (after subtracting off the eigenvalue at
$1$) \eqref{C.scdef}.

\begin{prop}
  The operator $C_A$ is a second-order differential operator on $\C[T]^W$ with
  associated biderivation $M_A$.  Moreover, it is diagonal in the basis of
  irreducible characters, with eigenvalue $A(\lambda+\rho)-A(\rho)$ on the
  character with highest weight $\lambda$.
\end{prop}

\begin{proof}
  We first show that it is diagonal on the irreducible characters.  If
  $\rho$ is a character of $G$, then $C_A$ is given by \eqref{C.scdef} and the claim follows from the Weyl character formula.  Otherwise,
  there is a double cover $G^+$ for which $\rho$ is a weight, and (since
  isogenies are \'etale) a corresponding induced pairing $A^+$.  Any irreducible
  character of $G$ is an irreducible character for $G^+$, and the action of
  $C_{A^+}$ agrees with the action of $C_A$ where both are defined.

  We have thus seen that $C_A$ preserves $\C[T]^W$.  Since $\mcD_A$ is a second-order differential operator, and $C_A - \mcD_A$ is a derivation by \eqref{C.altdef}, $C_A$ is also a second-order differential operator with the same associated biderivation as $\mcD_A$.
\end{proof}

\begin{rem}
  When $G$ is simply connected, the eigenvalue of $C_A$ on the irreducible
  character of highest weight $\lambda$ agrees with the value of the
  quadratic Casimir element on the corresponding representation.  (See,
  e.g., \cite[Prop.~10.6]{Hall:2nd}; in the reference, one fixes a
  $G$-invariant pairing on $\mfg$ rather than a $W$-invariant pairing
  on $\mft^*$, but these have the same degrees of freedom.)  In
  particular, when $G$ is semisimple, $C_A$ extends to a second-order
  operator on $\C[G]$ which is invariant under both left and right
  multiplication.  The analogous statement for $G$ reductive follows by
  pulling $C_A$ back to $G'\times Z(G)^\ic$.
\end{rem}

\begin{rem}
  If $G$ has maximal compact subgroup $K$ and $T$ contains a maximal torus
  $T_K$ of $K$, then Haar measure on $K$ induces a $W$-invariant measure on $T_K$,
  and the density of that measure relative to Haar measure on $T_K$ is (up
  to an irrelevant scalar multiple) nothing other than $\Delta$.  In
  particular, since $\mcD_A$ is self-adjoint relative to Haar measure on
  $T_K$, $C_A$ is self-adjoint relative to the measure induced from $K$,
  which is why its eigenfunctions form an orthonormal basis relative to
  that measure.
\end{rem}

In particular, since $C_A$ and $\mcD_A$ have the same associated
biderivation, the $G$-invariant derivations obtained from $C_A$ and the
fundamental characters agree with the derivations $D_i$ obtained from
$\mcD_A$ via \eqref{Di.def}.

\subsection*{Computation}
Note that both $\mcD_A$ and $C_A$ are reasonably straightforward to compute
on $\C[G]^G$: $\mcD_A$ is diagonal on monomials and $C_A$ is
diagonal on irreducible characters.  So to apply the operators to an
element of $\C[G]^G$, it suffices to expand it into irreducibles
(or monomials) and apply the relevant diagonal operator.  This gives the
coefficients of the matrix $M$ as an expansion in irreducibles (or
monomials), and thus one reduces to the problem of expressing those
coefficients as polynomials in the fundamental characters.

Having done that, it is then straightforward to write down the $n$
equations cutting out the critical locus on $\ClG$ corresponding to any
(polynomial) class function, and apply the usual Gr\"{o}bner basis
techniques to find the components of the critical locus.  Magma code for
this computation is provided in the GitHub repo for this article.  In the
code, $C_A$ is referred to as the Casimir and $\mcD_A$ is referred to as
the pseudo-Casimir.

\begin{eg}[$G_2$] \label{G2.eg}
Consider the case where $G$ is simple of type $G_2$, where the roots are numbered as in the diagram
\[
\dynkin[label,edge length=\edgelen,backwards,ordering=Carter] G2
\]
The weight $\omega_1$ dual to the simple coroot $\av_1$ is the highest short root, and the weight $\omega_2$ is the highest (long) root.  We take $f_i$ to be the character of the irreducible representation with highest weight $\omega_i$.  With respect to this numbering, the matrix $M$ is
\[
\left( \begin{smallmatrix}
4f_1^2 - 4f_2 - 16  f_1 - 28 & 6f_1f_2 - 14 f_1^2 + 14 f_2 - 16 f_1 + 14 \\
* & -12f_1^3 + 12 f_2^2 + 24 f_1 f_2 - 20 f_1^2 + 8 f_2 + 44 f_1 - 28
\end{smallmatrix} \right)
\]
where the $*$ in the lower left indicates an entry determined by symmetry (and omitted to make the matrix easier to read).

The corners are
\begin{equation} \label{G2.corners}
(f_1, f_2) = 
(-2, 5), \quad (-1, -2), \quad \text{or} \quad (7, 14).
\end{equation}

The critical points of $f_2$, the character of the adjoint representation, are the locus cut out by the second row (or column) of $M$, which has four components, three of which are the corners.  The remaining component is $(f_1, f_2) = (7/9, 10/27)$.  Since no critical point has $f_2$ less than $-2$ and $-2$ occurs as a value of $\Tr \Ad$ on the compact group (because $-1$ belongs to the Weyl group, Example \ref{minus1}), we conclude that $\min \Tr \Ad = -2$, verifying Theorem \ref{simple.thm} in this case.

The critical points of $f_1$, the character of the highest short root representation, are the locus cut out by the first row (or column) of $M$.  This consists solely of the corners \eqref{G2.corners}.
\end{eg}

\section{Critical points from the compact group}  \label{Wders.crit}

Suppose we have computed the
critical points over $\C$, but are only interested in the critical points
lying in the compact group, as in Theorem \ref{simple.thm}.
  If we were given the critical point as an
element of the group, this would be relatively straightforward: we need
simply ask whether the point is in the real locus (of the compact real
form).  Of course, what we are actually given is a point on the {\em
  quotient} $T/W$, and the real locus of the quotient is in general bigger than
the quotient of the real locus.  For instance, for $A_1$, the
coordinate on the quotient is $s+1/s$, which is certainly real
for $\diag(s, 1/s) \in \SU(2)$, but is also real for $\diag(s, 1/s) \in \SL_2(\R)$.
So we are left with the following question: Given a real point on the
quotient, how can we tell if it is in the image of the compact group?

It turns out that if we make a suitable choice of pairing $A$, the matrix
$M$ with entries $M_A(f_i,f_j)$ also answers this question.  We continue
the assumption of \S\ref{Wders.2} that $G$ is a connected reductive group
whose derived subgroup is simply connected.

We have already observed that any choice of isomorphism
$(z_1,\dots,z_n)\:T\xrightarrow{\sim} \Gm^n$ induces a basis of $\mft^*$.  This
corresponds to a basis-independent statement: any homomorphism
$\Hom(T,\Gm)$ induces (by taking the logarithmic differential) an element
of $\mft^*$, and the resulting homomorphism
\[
\Hom(T,\Gm)\to \mft^*
\]
induces an isomorphism
\[
\Hom(T,\Gm)\otimes \C\xrightarrow{\sim} \mft^*.
\]
As a result, we may view $A$ as being a quadratic form on the weight lattice
rather than on the dual of the Cartan.  Therefore it makes sense to require that
\begin{equation}  \label{A.hyp}
\text{$A$ is both real and positive definite}
\end{equation}
which we do for the remainder of the paper.

\begin{defn} \label{Msig.def}
Let $\sigma$ be the automorphism of the maximal torus
$\sigma(t)=t^{-1}$.  Define $M_A^\sigma$ to be the biderivation $(f,g)\mapsto
M_A(\sigma(f),g)$.

Given a choice of coordinates $f_1,\dots,f_n$, we get a
corresponding matrix $M^\sigma$ with entries
\[
(M^\sigma)_{ij} = M^\sigma_A(f_i,f_j) = M_A(\sigma(f_i),f_j)
=
\sum_{1\le k\le n} \partial_{f_k}(\sigma(f_i)) M_{kj};
\]
in other words, $M^\sigma=J(\sigma)^\trans M$ where $J(\sigma)$ is the
Jacobian of the action of $\sigma$ on $T/W$ with $(i,j)$-entry $\partial_{f_i}(\sigma(f_j))$.
\end{defn}

The matrix $M^\sigma$ is no longer symmetric, but is instead hermitian in a suitable sense.

\begin{lem} \label{Msig.herm}
  We have
  \[
  M^\sigma_A(g,f) = \sigma(M^\sigma_A(f,g))
  \]
  for all $f, g \in \C[T/W]$.
\end{lem}

\begin{proof}
  Since $\sigma$ acts as $-1$ on the weight lattice, we find that $A$ is
  $\sigma$-invariant, and thus
  \begin{align*}
  M^\sigma_A(g,f)
 & =
  M_A(\sigma(g),f)
  =
  M_A(f,\sigma(g)) \\
  &=
  \sigma(M_A(\sigma(f),\sigma(\sigma(g))))
  =
  \sigma(M^\sigma_A(f,g)).
  \end{align*}
\end{proof}

In matrix terms, Lemma \ref{Msig.herm} says:
\begin{equation}\label{herm.prop}
  (M^\sigma)^\trans = \sigma(M^\sigma)
\end{equation}

\begin{eg}[$\GL_2$, continued]  \label{GL2.eg3}
We return to the case $G = \GL_2$ from Example \ref{GL2.eg2}, noting that the choice of $A$ satisfies .  Since $\sigma(z_i) = 1/z_i$, we have $f_1 \sigma = 1/z_1 + 1/z_2$, $f_2 \sigma = 1/(z_1 z_2)$, 
\[
J(\sigma) = 
\left( \begin{smallmatrix}
1/f_2 & 0 \\ 
-f_1/f_2^2 & -1/f_2^2
\end{smallmatrix} \right), \eand 
M^\sigma = J(\sigma)^\trans M = 
\left( \begin{smallmatrix}
-2 & -f_1 \\ 
-f_1/f_2 & -2
\end{smallmatrix} \right).
\]
Since $f_1/f_2 = \sigma(f_1)$, we have verified that $M$ satisfies the hermitian property \eqref{herm.prop}.
\end{eg}

\begin{eg} \label{Msig.M}
Continue Example \ref{simple.A}, where $G$ is simple and simply connected. Acting on the root system, $-1$ and $-w_0$ for $w_0$ the longest element of the Weyl group $W$ are in the same coset of $W$.  Since $-w_0$ permutes the simple roots, $\sigma$ permutes the fundamental characters in the same way, and $M^\sigma$ is obtained from $M$ by permuting the rows of $M$ in the same manner.
\end{eg}

Let $T$ be a compact torus.  There is a natural corresponding {\em split}
torus
\[
T_s:=\Hom_\Z(\Hom(T,S^1),\Gm)
\]
over $\Z$, such that $T_s(\C)$ is naturally isomorphic to the
complexification of $T$.  (This follows immediately from the canonical
isomorphism between the complexification of $S^1$ and $\C^\times$.)  The
subgroup $T\subset T_s(\C)$ is then characterized as the elements $t$ such
that $\sigma(t)=\bar{t}$.  (Complex conjugation makes sense because the
base change of $T_s$ to $\C$ factors through the base change to $\R$.)  For any
$f\in \R[T_s]^W$ and $t\in T$, we find
\[
\sigma(f)(t) = f(\overline{t}) = \overline{f(t)},
\]
so that for $t\in T$,
\[
M^\sigma_A(f,g)(t)
=
\overline{M^\sigma_A(g,f)(t)}.
\]
We fix a choice of coordinates $f_1,\dots,f_n\in \R[T_s/W]$, which we may do by viewing $T_s$ as the maximal torus in a split form of $G$ and applying Proposition \ref{torus.prod}.  (In the case of most interest, $G$ will be simple and the coordinates are naturally defined over $\Z$.)  With respect to these coordinates,
the associated matrix $M^\sigma(t)$ is
hermitian for $t\in T$. 

The characterization of $T\subset T_s(\C)$ identifies $T$ as the real locus
of a twisted real form of $T_s$, which we also denote by $T$.  In
particular, there is a quotient scheme $T/W$ defined over $\R$, which is a
twisted real form of $T_s/W$.  The real points $x\in (T/W)(\R)$ are
characterized by the condition that
\[
\sigma(f)(x) = \overline{f(x)}
\]
for all $f\in \R[T_s/W]$, and thus $M^\sigma$ is hermitian at any point
of $(T/W)(\R)$.

We aim to prove:

\begin{thm} \label{neg.sd}
Let $T$ be a maximal torus in a compact Lie group $G$ whose derived subgroup is simply connected.  
 An element $x \in (T/W)(\R)$ is in the image of $T(\R)$ if and only if $M^\sigma(x)$ is negative semidefinite.
\end{thm}

\begin{eg}[$\SL_2$, continued] \label{A1.eg3}
For $G = \SL_2$, $M^\sigma$ is the 1-by-1 matrix with entry $t^2 - 4$.  For $t \in \R$, this is negative semidefinite if and only if $-2 \le t \le 2$.  On the other hand, given a $t$ in this range, the complex number $s := t/2 + i \sqrt{1 - t^2/4}$ has length 1 and $t = s + 1/s$.  This verifies Theorem \ref{neg.sd} for $G = \SL_2$.
\end{eg}

We will prove Theorem \ref{neg.sd} in the following abstracted setting.

Let $V\cong \R^n$ be Euclidean space, and let $\Gamma$
be a discrete subgroup of the ``(real) affine orthogonal group'' $V\rtimes O(V)$ with finite image in $O(V)$.
The action of $\Gamma$ on $V$ extends naturally to an
action on $V\otimes \C$, and the quotient $X:=(V\otimes \C)/\Gamma$ has a
natural structure as an affine scheme over $\R$ where the affine space is 
associated to the orthogonal complement of $\R(\Gamma\cap V)$. Indeed, the quotient
by the lattice $\Gamma \cap V$ is a product of an affine space and an algebraic
torus, with real structure such that the real locus is the product of a
vector space and a topological torus. This is clearly an affine scheme,
and $(V\otimes \C)/\Gamma$ is the quotient of this affine $\R$-scheme by the
real and finite group $\Gamma/(\Gamma \cap V)$.

There is a natural symmetric biderivation on $V\otimes \C$ given by
$(f,g)\mapsto \nabla f\cdot \nabla g$.
This biderivation is clearly $\Gamma$-invariant,
so descends to a symmetric biderivation on $X$. In particular, at any
smooth point $x\in X$, it determines a symmetric pairing on the
cotangent space at $x$. If $x\in X(\R)$, this pairing will be
real-valued, and thus in particular has a well-defined signature.

For $f$ defined on a neighborhood of $x_0\in X(\R)$, we may define
$\bar{f}$ by $\bar{f}(x)=\overline{f(\bar{x})}$. (This is a priori defined on a different neighborhood of $x_0$, but of course both are
defined on the intersection.)

Put $\Reg(X)$ for the regular locus of the affine scheme $X$, equivalently the smooth locus \stacks{0B8X}.

\begin{lem} \label{RegX.lem}
For $x_0\in \Reg(X)(\R)$, the following are equivalent:
\begin{enumerate}[(a)]
\item \label{RegX.im} $x_0$ is in the image of $V$.
\item \label{RegX.pair} For any analytic function $f$ defined on a neighborhood of $x_0$,
\[
(\nabla f\cdot \nabla \bar{f})(x_0)\ge 0.
\]
\item \label{RegX.pos} The pairing on the cotangent space at $x_0$ is positive semidefinite.

\end{enumerate}
\end{lem}

\begin{proof}
Let $f_1,\dots,f_n$ be real analytic local coordinates at $x_0\in
X(\R)$. The gradients give a basis of the cotangent sheaf near $x_0$,
and the pairing relative to that basis has matrix $M = (M_{ij})$ for 
\[
M_{ij} := \nabla f_i\cdot \nabla f_j,
\]
where the gradients are taken on the pullback to $V\otimes \C$. In
particular, the pairing is positive semidefinite, \ref{RegX.pos}, if and only if $M$ is positive
semidefinite, if and only if $\vec\alpha M\vec\alpha^\dagger\ge 0$ for any complex
row vector $\vec\alpha$. So \ref{RegX.pos} is equivalent to the statement that
\[
(\nabla f\cdot \nabla \bar{f})(x_0)\ge 0
\]
where $f=\sum_i \alpha_i f_i$, for every $\vec\alpha \in \C^n$.  Thus \ref{RegX.pair} implies \ref{RegX.pos}.  

Now suppose \ref{RegX.im}, that $x_0$ is in the image of $v_0\in V$. Then any analytic
function $f$ defined on a neighborhood of $x_0$ pulls back to an
analytic function $g$ defined on a neighborhood of $v_0$, and one has
 \[
 (\nabla f\cdot \nabla \bar{f})(x_0)
 =
 (\nabla g\cdot \nabla \bar{g})(v_0)
 =
 \sum_{1\le i\le n} |(\partial_i g)(v_0)|^2
 \ge 0,
 \]
 i.e., \ref{RegX.pair}.

Now suppose the pairing is positive
semidefinite at $x_0$, \ref{RegX.pos}. Since $\Gamma$ is discrete, we may choose a
topologically open neighborhood $U$ of $x_0$ so that each connected
component of its preimage in $V\otimes \C$ contains a unique preimage of
$x_0$. Let $v_0$ be such a preimage, with corresponding component
$U_0$, and consider the polynomial
 \[
 (v-\RE(v_0))\cdot (v - \RE(v_0))
 \]
on $V$. This is entire, so in particular defines an analytic
function on $U_0$, which, since the stabilizer $\Gamma_{v_0}$ fixes
$\RE(v_0)$ and respects the Euclidean inner product, is
$\Gamma_{v_0}$-invariant, so descends to an analytic function $f$ on a
neighborhood of $x_0$.

In particular, the positivity assumption tells us that $\nabla f\cdot
\nabla f$ is nonnegative at $x_0$. But we know the preimage of $f$ near
$v_0$, and may thus directly compute that $\nabla f\cdot \nabla f=4f$.
Evaluating this at $x_0$ gives
 \[
 0\le (\nabla f\cdot \nabla f)(x_0) = 4f(x_0) = 4 (i\IM(v_0))\cdot
(i\IM(v_0))
 = -4 \left|\IM(v_0)\right|^2.
 \]
 We thus see that $\left|\IM(v_0)\right|^2\le 0$, implying $\IM(v_0)=0$, \ref{RegX.im}.
\end{proof}

\begin{rmk*}
    Lemma \ref{RegX.lem} is reminiscent of the main theorem 0.10 in 
\cite{ProcesiSchwarz}, but that result concerns a linear representation of a 
compact group rather than an affine representation of a discrete group. 
When our group is finite, the lemma is a special case of their result, with a 
very similar proof.  (They not only allow the group to be compact, but 
give the appropriate condition at singular points of the quotient.) 
Presumably there is a common generalization to affine representations of 
(topologically) closed subgroups of $V\ltimes O(V)$ with closed image in 
$O(V)$.
\end{rmk*}

We now return to the setting of Theorem \ref{neg.sd}.

\begin{proof}[Proof of Theorem \ref{neg.sd}]
Write $T$ for the maximal torus in $G$, $W$ for the Weyl group $N_G(T)/T$, and $\Tt$ for the universal cover of $T$.
 Then $T(\C)/W$ may be
rewritten as a quotient $\Tt(\C)/\pi_1(T).W$ of the type denoted by $X$ above, so
that Lemma \ref{RegX.lem} applies. Since $T(\C)/W$ is the product of an
algebraic torus and an affine space (Prop.~\ref{torus.prod}\ref{St.inv.prod}), it has a global system of
coordinates $f_1,\dots,f_n$, and thus we find that $x\in (T/W)(\R)$ is
in the image of $\Tt(\R)$, or equivalently in the image of the
compact quotient $T(\R)$, if and only if the hermitian matrix $\Mh$ with $(i,j)$-entry
\[
\Mh_{ij}:=\nabla f_i\cdot \nabla \bar{f}_j.
\]
is positive semidefinite.

Note that here the gradients are computed on $\Tt$.  The natural map $\Tt\to T$ that respects real forms is given by $\vec{x}\mapsto \exp(2\pi i\vec{x})=:\vec{z}$, and thus the ``algebraic'' gradient with coordinates $z_j\partial_{z_j}f$ is $-1/(2\pi i)$ times the ``analytic'' gradient with coordinates $\partial_{x_i}$. It follows that the matrix $M^\sigma$ in Theorem \ref{neg.sd} is $-1/(4\pi^2) \Mh$, so is negative semidefinite if and only if $\Mh$ is positive semidefinite.
\end{proof}

 We remark that the result also applies to the case
of a finite Coxeter group acting on a vector space, since again the
quotient is smooth and has global coordinates (i.e., generators of the
ring of polynomial invariants).

Returning to the main theme of the paper, we have the following corollary of Theorem \ref{neg.sd}.

\begin{cor} \label{neg.G}
Let $G$ be a semisimple and simply connected compact Lie group.  Let $g$ be a semisimple element of $G(\C)$ whose image in $\ClG(\C)$ lies in $\ClG(\R)$.  Then $g$ is conjugate under $G(\C)$ to an element of $G(\R)$ if and only if $M^\sigma(g)$ is negative semidefinite.
\end{cor}

We can restate this in terms of the split real form of the group: the
hypothesis becomes that $g$ and $\bar{g}^{-1}$ are conjugate in $G(\C)$,
while the conclusion becomes that $g$ is contained in a compact subgroup of
$G(\C)$ if and only if $M^\sigma(g)$ is negative semidefinite.

\begin{proof}
Let $T$ be a maximal torus of $G$ containing $g$ and write $W$ for $N_G(T)/T$.  Then by Steinberg \cite[6.4]{St:reg} $\ClG \cong T/W$ over $\R$, so
the hypothesis on $g$ implies that the image of $g$ in $(T/W)(\C)$ belongs to $(T/W)(\R)$.  Theorem \ref{neg.sd} gives: $M^\sigma(g)$ is negative semidefinite if and only if there exists a $t \in T(\R)$ with the same image in $(T/W)(\R)$.  Applying Steinberg again, we see that this condition is equivalent to $g$ being $G(\C)$-conjugate to an element of $T(\R)$.
\end{proof}
  
\section{Corners and the rank of \texorpdfstring{$M$}{M}} \label{corners.sec}

Continue the assumption of sections \ref{Wders.2} and \ref{Wders.crit} that $G$ is a connected
reductive algebraic group over $\C$ whose derived subgroup is simply connected.  We
now investigate the significance of the rank of $M$ at an element $g$.
Note that $M$ and $M^\sigma$ have the same rank at $g$, because $M^\sigma =
J(\sigma)^\trans M$, where $\sigma$ is an automorphism, so $J(\sigma)$ is
invertible.

\begin{lem}
  Let $V$ be a finite-dimensional complex representation of the finite
  group $G$, and let $q\in (S^2(V^*))^G$ be a nonsingular $G$-invariant
  quadratic form on $V$.  Then the restriction of $q$ to $V^G$ is
  nonsingular.
\end{lem}

\begin{proof}
  Since complex representations of finite group are semisimple, there is a
  natural direct sum decomposition $V\cong V^G\oplus W$, in which $W$ is
  the sum of all other isotypic components.  For any $v\in V^G$, $w\in W$,
  $G$-invariance of $q$ implies $q(v+w)=q(v)+q(w)$, so that if $v$ is in
  the radical of $q|_{V^G}$, it is in the radical of $V$, and thus $v=0$.
\end{proof}

\begin{lem} \label{omega.fiber}
Let $t\in T(\C)$, with stabilizer $W_t$.  Then the fiber at $t$ of the natural map
\[
\psi \: \Omega_{\C[T/W]}\otimes_{\C[T/W]}\C[T]\to \Omega_{\C[T]}
\]
has image $(\Omega_{\C[T]}|_t)^{W_t}$, which is canonically identified with $(\mft^*)^{W_t}$.
\end{lem}

Since $T$ is abelian, the map $t \mapsto w(t) t^{-1}$ for $w \in W$ is a homomorphism $T \mapsto T$.  In the proof, we will use the notation $H_0(W; T)$ for the coinvariants, i.e., $T$ modulo the subgroup consisting of elements $w(t) t^{-1}$ for $w \in W$ and $t \in T(\C)$ \cite[Ch.~6]{Weibel:HA}.  Note that the quotient $H_0(W; T)$ of $T$ is itself an algebraic group over $\C$, and therefore is itself a torus.  The vector space $H_0(W; \mft)$ is defined similarly.

\begin{proof}
Consider the factorization
\[
\Omega_{\C[T/W]}\otimes_{\C[T/W]}\C[T]
\to
\Omega_{\C[T/W_t]}\otimes_{\C[T/W_t]}\C[T]
\to
\Omega_{\C[T]}
\]
of $\psi$.
Since $T/W_t\to T/W$ is \'etale at $t$, the first map is an isomorphism near $t$, and thus we may assume without loss of generality that $W_t=W$, i.e., that $t$ is fixed by $W$.  Moreover, since $\psi$ is equivariant under translations by $T^W$, we may assume that $t=1$.

In this case, $\psi$ is $W$-equivariant, and thus the restriction to $1$
is still $W$-equivariant; since $W$ acts trivially on $\Omega_{\C[T/W]}$,
it follows that the image is contained in
$(\Omega_{\C[T]}|_1)^W$.  That space is the same as $(\mft^*)^W$ since $\Omega_{\C[T]}\vert_1$ is canonically identified with $\mft^*$.

To get the other inclusion, consider the natural quotient morphism $T\to 
H_0(W;T)$ of group schemes.  This again induces a pullback map on 
differentials:
\[
\Omega_{\C[H_0(W;T)]}\otimes_{\C[H_0(W;T)]} \C[T]
\to
\Omega_{\C[T]}.
\]
The fiber at 1 of the pullback map is the dual of the corresponding 
morphism of Lie algebras,
\[
H_0(W;{\mft})^*\to {\mft}^*,
\]
which is injective with image $({\mft}^*)^W$.  Since the morphism 
$T\to H_0(W;T)$ factors through the quotient {\em scheme} $T/W$, the map 
on differentials factors through $\psi$, and thus we have a factorization
\[
H_0(W;{\mft})^*\to \Omega_{\C[T/W]}|_1\to {\mft}^*,
\]
so that the image of $\psi|_1$ contains the image $({\mft}^*)^W$ 
of the composite.
\end{proof}

%

\begin{prop} \label{corank}
For each semisimple 
$g \in G(\C)$, the rank of $M(g)$ is equal
to the dimension of the center of the centralizer of $g$.
\end{prop}

There is an analogue of this result in the setting of Lemma \ref{RegX.lem}.  A linear version of \ref{corank} is part of the main theorem of \cite{St:inv}.

\begin{proof} 
We may reduce to the case of a (complex) point of the maximal torus.  We show that each of the numbers in the statement is equal to the dimension of $(\mft^*)^{W_g}$.  

First consider the rank of $M(g)$.  It is the same as the rank of the restriction of $A$ to the image of the fiber at $t$ of the map $\psi$ from Lemma \ref{omega.fiber}, 
which we have just seen is nothing other than $(\mft^*)^{W_g}$.  Since
$A$ is nonsingular and $W_g$-invariant, its restriction to the subspace of
$W_g$-invariants is again nonsingular, proving the desired equality for the rank of $M(g)$.

For the other dimension, we use that the torus $T$ is also a maximal torus in $Z_G(g)$, so the stabilizer $W_g$ of $g$ is the Weyl group $N_{Z_G(g)}(T)/T$ of $Z_G(g)$.  Now $Z_G(g)$ is reductive  \cite[Prop.~13.19]{Borel}, so we may write $T = Z T'$ where $Z$ is the identity component of the center of $Z_G(g)$ and $T'$ is a maximal torus of the semisimple part of $Z_G(g)$.  Each of these tori is stable under $W_g$, so 
\[
\mft^{W_g} = \Lie(Z)^{W_g} \oplus \Lie(T')^{W_g} = \Lie(Z),
\]
where $\Lie(T')^{W_g} = 0$ 
because the Weyl group of the semisimple group $Z_G(g)$ over $\C$ acts without fixed lines on the Lie algebra on the vector space generated by the coroots \cite[VI.1.2, Cor.]{Bou:g4}.  

Since $\mft$ has a non-degenerate $G$-invariant (hence also $W_g$-invariant) bilinear form, as representations of $W_g$, $\mft^*$ is equivalent to $\mft$, and so we have proved
\[
\dim \, (\mft^*)^{W_g} = \dim \mft^{W_g} = \dim \Lie(Z),
\]
as required.
\end{proof}

One special case of the proposition is: $M$ is invertible at $g$ if and only if $g$ is regular.  Here are more.

\begin{cor} \label{corner.ss}
If $G$ has a corner, then $G$ is semisimple (and not merely reductive).
\end{cor}

\begin{proof}
By hypothesis, $G$ is reductive, so the identity component $Z_G(G)^\ic$ of the center of $G$ is a torus $S$.  For every $g \in G(\C)$, the center of $Z_G(g)$ contains $S$.  On the other hand, if the conjugacy class of $g$ is a corner, then by Prop.~\ref{corank} $\dim Z(Z_G(g)) = 0$, so $S = 1$.
\end{proof}

For context, we momentarily leave the domain of algebraic groups over $\C$ to mention the following result, which provides another view on the preceding corollary and Cor.~\ref{corner.crit}:
\begin{prop}
Let $G$ be a connected reductive real Lie group that is not semisimple.  Then for any $g\in G$, there is a real analytic $G$-invariant function on $G$ which is not critical at $g$.
\end{prop}

\begin{proof}
It suffices to find such a function that is pulled back from a nontrivial character $G\to \R_{>0}$ or $G\to S^1$, and thus we may reduce to the case $\dim(G)=1$.  On $\R_{>0}$, the coordinate function has no critical points, while on $S^1$, the only critical points of $\Re(z)$ are $\pm 1$, which are not critical points for $\IM(z)$.
\end{proof}

Return now to the setting where $G$ is an algebraic group over $\C$, and we seek to say more about the corners of $G$. 

\begin{cor} \label{corners.r}
Let $G$ be a semisimple and simply connected algebraic group over $\C$, and suppose $g \in G(\C)$ is semisimple and its conjugacy class is a corner.
\begin{enumerate}[(1)]
\item \label{corners.r.g} $g$ has finite order and the conjugacy class of $g^m$ is also a corner for every integer $m$.
\item \label{corners.r.simple} If $G$ is simple of rank $r$, then $G$ has $r+1$ corners and the order of $g$ divides $a$ times the exponent of $Z(G)$, where $a$ is the coefficient of some simple root in the expression for the highest root of $G$ as a linear combination of simple roots.
\end{enumerate}
\end{cor}

\begin{proof}
Put $H_m := Z_G(g^m)^\ic$.  It is a reductive subgroup of $G$ and contains a maximal torus $T$ containing $g$.  Evidently, 
\begin{equation} \label{corners.r.1}
H_m \supseteq H_1 \supseteq T.
\end{equation}
For every $m$, the center of the reductive group $H_m$ is contained in every maximal torus of $H_m$ \cite[Prop.~21.7]{Milne}, so we find
\[
T \supseteq Z(H_1) \supseteq Z(H_m).
\]
Prop.~\ref{corank} says that $H_1$ is semisimple, so $Z(H_1)$ is finite, ergo $g$ has finite order.  The reductive group $H_m$ also has finite center, so $H_m$ is semisimple (instead of merely reductive)  and the class of $g^m$ is a corner by Prop.~\ref{corank}.  This proves \ref{corners.r.g}.

Now suppose that $G$ is simple of rank $r$.  Since $g$ has finite order, it is conjugate to an element of a maximal compact subtorus $S$ of $T$.  That subtorus is in turn homeomorphic to the closure of an alcove in the coroot space $V$ \cite[\S{IX.5.2}, Cor.~1]{Bou:g7}.  (The coroot space is the real vector space with basis the coroots dual to the simple roots of $G$ with respect to $T$.)  Because $Z_G(g)^\ic$ is semisimple, its root system also has rank $r$, which means that the element $v \in V$ corresponding to $g$ is a vertex of the  closure of the alcove, compare \cite[p.~11]{Reeder:tor}.  The vertices are of the following form.  Put $\ov_i$ for the fundamental dominant weight of the dual root system, so that $\qform{ \ov_i, \alpha_j} = \delta_{ij}$ (Kronecker delta) for $\alpha_1, \ldots, \alpha_r$ the simple roots of $G$ in a pinning involving $T$.  Then $v$ is a vertex if and only if $v = 0$ (in which case $g = 1$) or $v = a_i^{-1} \ov_i$, where $a_i$ is the coefficient of the simple root $\alpha_i$ in the expression $\sum a_j \alpha_j$ for the highest root.  In the second case, the corresponding $g$ has $g^{a_i} \in Z(G)$.  This proves \ref{corners.r.simple}.
\end{proof}

\begin{rmk}[naming the corners] \label{naming}
Since a corner is a conjugacy classes of elements of finite order, when $G$ is simple the corner can be identified uniquely via its Kac coordinates as described in \cite{Kac:aut}, \cite{Se:Kac}, \cite[\S2.2]{Reeder:tor}, or \cite[\S6]{Lehalleur}.  In the notation of the corollary, this is a list $(s_0, s_1, \ldots, s_r)$ of non-negative integers, not all zero, with greatest common divisor 1.  The corners are those conjugacy classes with Kac coordinates that have a 1 in one entry and zeros elsewhere; the corner $a_i^{-1} \ov_i$ has $s_i = 1$.

Alternatively, suppose $G$ is simple of type $G_2$, $F_4$, or $E_8$.  The conjugacy class of some $g \in G(\C)$ is a corner and $p$ is a prime not dividing the order of $g$.  One could use the methods of \autoref{finite.sec} to construct a finite-field analogue $G(\F_p)$ of $G$ over the finite field $\F_p$ and relate $g$ with an element $\bar{g}$ of the same order in $G(\F_p)$.  In this way, it may make sense to refer to the conjugacy class of $g$ itself using a label 2A, 3B, etc., that is  used for the conjugacy class of $\bar{g}$ in the style of the Atlas \cite{atlas}.  We do not pursue this here.
\end{rmk}

\begin{eg}[type $G_2$, continued] \label{G2.contd}
For $G$ of type $G_2$, the 3 corners were found in Example \ref{G2.eg}.  Continue the notation of that example.
The minimum of $f_1$ on the critical points, $-2$, comes from the compact $G_2$ (because it comes from a corner), verifying Proposition \ref{short.prop} in this case.

We remark that the remaining critical point of $f_2$ has
\[
M(7/9, 10/27) = \left( \begin{smallmatrix} -3200/81 & 0 \\ 0 & 0 \end{smallmatrix} \right).
\]
This matrix is negative semidefinite.  Since $M = M^\sigma$ (Example \ref{Msig.M}), this critical point comes from the compact form.
\end{eg}

The values of the fundamental characters of $G_2$, $F_4$, or $E_8$ on the corners, exhibited in \eqref{G2.corners} and Tables \ref{F4.corners} and \ref{E8.corners}, are integers.  This is an example of a more general phenomenon.

\begin{prop}
Suppose $G$ is simple and simply connected.  Then $-1$ is in the Weyl group if and only if every character takes integer values on the corners.
\end{prop}

One could equally well state the criterion as ``$Z(G)$ has exponent at most 2'', see \cite[\S{VI.4}, Exercise 7]{Bou:g4}.  One direction is obvious, because the Weyl group acts trivially on $Z(G)$, whereas $-1$ acts by inverting elements of $Z(G)$.

\begin{proof}
Suppose first that $-1$ is in the Weyl group.  Let $g$ be a representative of the conjugacy class corresponding to a corner, and put $m$ for the order of $g$.  The character $\chi$ restricts to a complex character of the subgroup $\Z/m$ generated by $g$, and therefore $\chi(g)$ is an algebraic integer in the cyclotomic field $\Q(e^{2\pi i/m})$.  Also, $\chi(g^{-1})=\overline{\chi(g)}$.  On the other hand, the assumption that $-1$ is in the Weyl group implies that $g$ and $g^{-1}$ are conjugate, so $\chi(g)$ belongs to the real subfield $\Q(\cos(2\pi/m))$.  The proves the claim if $m$ is 1, 2, 3, 4, or 6.

By Corollary \ref{corners.r}\ref{corners.r.simple}, $m$ divides $a$ times the exponent of center of $G$, for some coefficient $a$ appearing in the expression of the highest root.  For types $G_2$ and $F_4$, the center is trivial and $a$ is 4 or divides 6, which disposes of those cases.
For types $B$, $C$, and $D_{2n}$, $a$ is 1 or 2 while central elements have order dividing 2, so $m$ divides 4.

For $E_7$ and $E_8$, the same argument treats all the corners except one in each case, the one where $\ov$ is dual to the unique simple root in the Dynkin diagram that has three neighbors.  However, the Galois group of $\Q(\cos(2\pi/m))$ permutes the character values $\chi(g), \chi(g^2), \ldots$.  Since every power of $g$ is a corner (Cor.~\ref{corners.r}\ref{corners.r.g}) and $\chi$ takes integer values on all the other corners, $\chi(g)$ must be fixed by the Galois group and so belongs to $\Q$.

For the other direction, we suppose that $-1$ is not in the Weyl group, i.e., $G$ has type $A_n$ for $n \ge 2$, $D_n$ for odd $n \ge 2$, or $E_6$.  The center of $G$ is cyclic of order $m = n+1$, 4, or 3 respectively and is generated by an element $g$ corresponding to the fundamental coweight $\ov$ dual to the simple root $\alpha_x$ for $x = 1, n, 1$ respectively as numbered in \eqref{EDA}, see \cite[p.~298, Table 3]{OV:LGAG}.   In each case the coefficient in the highest root of the simple root corresponding to $\ov$ is 1 so $g$ is a corner.   The fundamental representation with highest weight $\omega_x$ is injective, so $g$ acts as an $m$-th root of unity $\zeta$, ergo $f_x(g) = \zeta f_x(1)$, which is not an integer since $m \ge 3$.
\end{proof}

\section{Fundamental representations of \texorpdfstring{$F_4$}{F4} and \texorpdfstring{$E_8$}{E8}} \label{E8.sec}

In Examples \ref{G2.eg} and \ref{G2.contd}, we applied the method of invariant derivations to determine the range of trace values on the fundamental irreducible representations of the compact simple Lie group of type $G_2$.  In this section, we will do the same for $F_4$ and $E_8$, although for $E_8$ we will only give complete information on the fundamental representations of dimension less than 100 million.  

\subsection*{Method}
We number the fundamental irreducible representations (equivalently, simple roots) as in Bourbaki, i.e., as in the diagrams
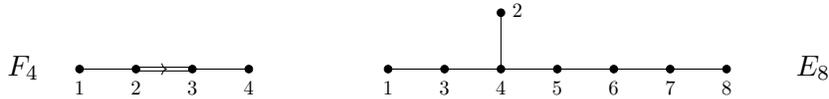
\[
F_4 \quad \dynkin[label,edge length=\edgelen] F4 \qquad \qquad
\dynkin[label,edge length=\edgelen] E8 \qquad E_8
\]
Put $f_i$ for the fundamental character with highest weight $\omega_i$.  Our aim is to compute $\min_{g \in G} f_i(g)$ for $G$ a compact form of $F_4$ or $E_8$, for as many choices of $i$ as possible.  

For $F_4$, we remark that $f_1 = \Tr \Ad$ and $f_4$ is the trace of the highest short root representation from \autoref{short.root}.  

For $E_8$, 
$f_8 = \Tr \Ad$.  The next smallest representation, the case $i = 1$, has dimension 3875 and was a particular object of study in \cite[\S7]{GG:simple}, \cite{ChayetG}, and \cite{DeMedtsVC}.  As a representation, it is an irreducible summand of $S^2(\mfe_8)$. The next smallest representation, with $i = 7$, is isomorphic to $(\wedge^2 \mfe_8)/\mfe_8$ (compare for example \cite{Adams:E8} or \cite{Guillot:lambda}).

\setcounter{stepnum}{-1}
\step{Preliminary step.} A key preliminary step is to compute the matrix $M$, which depends only on $G$.  For $E_8$, 
in contrast to the case of the other exceptional groups, this takes a measurable amount of time on a small computer --- between 7 and 8 minutes on the 2022 M2 MacBook Air being used to compose this text.  However, once $M$ has been calculated, it can re-used for different values of $i$.
 
\step{Compute the value of $f_i$ on the corners}, i.e., the vanishing locus of $M$.  There are $r + 1$ corners for $G$ of rank $r$ (Cor.~\ref{corners.r}).   For $F_4$ and $E_8$, we list the results in Tables \ref{F4.corners} and \ref{E8.corners}, where the header row of the table gives the value of $i$.  In the table, we have sorted the rows so that the first corner, i.e., the first data row, gives the dimension $d := f_i(1)$ of the representation.  In column $i$, we underline the minimum value $m$ of $f_i$ on a corner.  

\step{Compute the critical locus of $f_i$,} i.e., the vanishing locus of column $i$ of $M$.  Despite amounting to a tiny number of lines of Magma code, this is where the larger computations spend the most time.

\step{Examine each closed point.}  For $F_4$ and $E_8$, just as for $G_2$, the collection of critical points is zero-dimensional, so each component of the critical locus defines a number field $K$.  For each real embedding $\tau$ of each $K$, we compute $\tau(f_i) \in \R$.  Since our aim is to find the minimum of $f_i$, we discard $\tau(f_i)$ unless it lies in $[-d, m)$.  Indeed, we are looking for cases where $f_i < m$ (hence the upper bound) and an element with $f_i < -d$ cannot come from the compact group (hence the lower bound).  For those $\tau$ such that $\tau(f_i)$ does lie in the interval, we compute $\tau(M^\sigma)$ and check whether it is negative semi-definite by inspecting the sign of the determinants of the minors, which determines whether this critical point comes from the compact group (Theorem \ref{neg.sd}).

\begin{eg}[Fundamental representations of $F_4$] \label{F4.eg}
For $G$ of type $F_4$ and $i \ne 2$, the minimum of $f_i$ occurs at a corner.  In the remaining case of $f_2$, the minimum on the corners is $-14$, and there is a unique conjugacy class of critical points on the compact group whose trace is less than that.  That minimal value $m$ of the trace is
\[
\frac{98}{27} (1 - 2 \sqrt{7}) \approx -15.58,
\]
which is not an algebraic integer.
\end{eg}
\begin{table}[htbp]
\[
{\rowcolors{2}{white}{gray!25}
\begin{array}{rrrrrrrr}\toprule
1&2&3&4 \\ \midrule
52&1274&273&26 \\ 
20&154&\underline{\mathbf{-15}}&\underline{\mathbf{-6}}  \\
0&-10&5&-2  \\
-2&5&3&-1  \\
\underline{\mathbf{-4}}&\underline{-14}&-7&2  \\ \midrule
-4&\approx -15.58&-15&-6 \\
\bottomrule
\end{array}}
\]
\caption{Value of $f_i$ on each of the 5 corners of $F_4$.  The top row gives $i$, the bottom row gives the actual minimum of $f_i$.} \label{F4.corners}
\end{table}

\begin{table}[htbp]
\[
{\rowcolors{2}{white}{gray!25}
\begin{array}{rrrrrrrr}\toprule
1&2&3&4&5&6&7&8 \\ \midrule
3875& 147250& 6696000& 6899079264& 146325270& 2450240& 30380& 248 \\
35& 50& -960& \underline{-41888}& 3094& 832& \underline{\mathbf{-84}}& \underline{\mathbf{-8}} \\
3& \underline{\mathbf{-494}}& \underline{-2496}& -12320& \underline{-17290}& \underline{\mathbf{-1216}}& 140& 24 \\
5& -35& 0& 924& -231& 44& 14& -4 \\
\underline{\mathbf{-13}}& -8& 135& 411& 390& 71& 5& 5 \\
7& -10& 48& 224& -70& 16& 0& -4 \\
0& 0& 0& 14& 20& -10& 5& -2 \\
3& -8& 15& 19& -10& -1& 5& -3 \\
-5& 10& 0& 0& -50& 0& 4& 0 \\ \midrule
-13&-494&\approx -2524.74&?&?&-1216&-84&-8 \\ \bottomrule
\end{array}}
\]
\caption{Value of $f_i$ on each of the 9 corners of $E_8$.  The top row gives $i$, the bottom row gives the actual minimum of $f_i$.} \label{E8.corners}
\end{table}

\begin{eg}[Fundamental representations of $E_8$] \label{E8.eg}
For $G$ of type $E_8$ and $i = 1$, 2, 6, 7, or 8, the minimum value of $f_i$ on the compact $E_8$ occurs at a corner and is exhibited in  \autoref{E8.corners}.  Magma code to execute this computation is provided in the GitHub repo for this paper.  For any particular $i$, the computation takes longer and more computer memory if the dimension of the corresponding representation is larger.  The case $i =6$ has dimension about 2.5 million and took about 5 hours on an Intel Xeon Platinum 8380HL processor. 

Following the same procedure for the case $i = 3$, which has dimension about 7 million, less than 3 times as large as the $i = 6$ case, took about 5 days.
We found exactly one conjugacy class of critical points for $f_3$ 
on the compact $E_8$ where $f_3$ has value less than $-2496$, the minimum value on the corners.  For the closest root of $27z^3+10z^2-151z-49$ to 0, $z \approx -0.323628$, this minimal value of the trace is
\[
m := \frac{12327541632z^2 - 1318826423376z - 618329293056}{75346115} \approx -2524.74.
\]
The norm of $m$ in $\Q$ is $4894459524756541440/1159171$, so $m$ is not an algebraic integer.

For $E_8$, the cases of $f_4$ and $f_5$ remain to be investigated.
\end{eg}

\part{Approach via symmetry breaking} \label{part.break}

In Part \ref{part.Wders} we introduced a general technique and used it to verify part of Theorem \ref{simple.thm}, namely the values of $\mmin(1)$ for $G$ of exceptional type.  This part of the paper will begin by giving a different proof of those same values (\S\ref{serre.method}) and then move on to compute $\mmin(1)$ for $G$ simple of classical type in \S\ref{D.sec}--\S\ref{B.sec}.  Along the way we will provide a complete proof of Proposition \ref{short.prop}.

\section{Simple groups of exceptional type} \label{serre.method}

We now describe how to compute $\mmin(1)$ as in Theorem \ref{simple.thm} for $G$ compact simple of exceptional type.  The crux point is that $G$ is compact, so every element is semisimple, so it suffices to find the extreme points of $\Tr \Ad$ on a maximal torus (Lemma \ref{torus.reduce}).

\subsection*{Branching to a product of \texorpdfstring{$A_1$}{A1}'s}
Suppose now that $G$ has rank $n$ and $G$ contains a semisimple subgroup $H$ of type $A_1^{\times n}$.  This is the case when $G$ is exceptional and not of type $E_6$.  One can restrict the representation $\Ad$ (or the highest short root representation) of $G$ to $H$ in just a few lines of commands in a computer algebra system like Magma.

The highest weights of the irreducible summands in the restriction of the representation to $H$ give a formula for restriction of $\Tr \Ad$ to a maximal torus of $H$ in terms of the trace (on the natural representation) of each copy of $A_1$.  For example, a toral element in $A_1$ is diagonal with entries $s$, $1/s$ and trace $t = s + 1/s$.  Then the modules with the following highest weights contribute to the total trace via:
\begin{equation} \label{cheb.small}
\begin{array}{c|cccc}  \toprule
\text{weight} & 0&1&2&3 \\ 
\text{restriction} & 1 & t & t^2 - 1 & t^3 - 2t \\ \bottomrule
\end{array}
\end{equation}
And for a weight vector that has nonzero coordinates for more than one copy of $A_1$, we use the usual matrix identity $\tr(A \otimes B) = \tr(A) \tr(B)$.

In this way, we write $\Tr \Ad$ as a restriction to $H$ as a polynomial in $t_1, \ldots, t_n$.  Call this $f$.  For the convenience of the reader, we record in  \autoref{exgroups.table} the polynomials $f$ so obtained.
\begin{table}[hbtp]
\[
{\rowcolors{2}{white}{gray!25}
\begin{array}{cl} \toprule
\text{type of $G$}& \text{objective function $f$} \\ \midrule
G_2& t_1^3 t_2  - 2 t_1 t_2 +  t_1^2+ t_2^2 - 2\\
F_4& t_1 t_2 t_3 t_4 + t_1^2 + t_2^2 + t_3^2 + t_4^2+ t_1 t_2 + t_1 t_3+ t_1 t_4 + t_2 t_3  + t_2 t_4 + t_3 t_4  - 4\\
E_7&t_1t_2t_3t_4 + t_1t_2t_5t_6 + t_3t_4t_5t_6 + t_1t_3t_5t_7 +  t_2t_4t_5t_7 + t_2t_3t_6t_7 \\
&+ t_1t_4t_6t_7 + t_1^2 + t_2^2 + t_3^2  + t_4^2 + t_5^2 + t_6^2 + t_7^2 - 7 \\ 
E_8& t_1t_2t_3t_4 + t_1t_2t_5t_6 + t_3t_4t_5t_6 + t_1t_3t_5t_7 + t_2t_4t_5t_7 + t_2t_3t_6t_7 \\
&+ t_1t_4t_6t_7 + t_2t_3t_5t_8 + t_1t_4t_5t_8 + t_1t_3t_6t_8 + t_2t_4t_6t_8 + t_1t_2t_7t_8 \\
&+ t_3t_4t_7t_8 + t_5t_6t_7t_8 + t_1^2 + t_2^2 + t_3^2 + t_4^2 + t_5^2 + t_6^2 + t_7^2 + t_8^2 - 8\\ \bottomrule
\end{array}}
\]
\caption{The restriction $f$ of $\Tr \Ad$ from $G$ to a copy of $A_1^{\times n}$.} \label{exgroups.table} 
\end{table}

To find the critical points, we use that the critical points of $\Tr \Ad$ are where the Weyl-invariant differentials vanish, which by Example \ref{A1.eg} are the common solutions of
 \begin{equation} \label{inv.diff}
  (t_i^2 - 4) \partial_{t_i} f = 0 \quad \text{for all $i$}.
   \end{equation}
This defines a scheme so we can use the usual commutative algebra tools.

The points of this scheme are the critical points in the complexification of a maximal torus.  Some of these critical points have trace that lie outside the region claimed in Theorem \ref{simple.thm}, but that is no contradiction as long as they do not come from the compact group.  A necessary condition for a toral element to come from the compact group is that its trace on the adjoint representation of $A_1$ is inside the bound, i.e., $| t_i^2 - 2 | \le 2$.  We look at all the critical points that satisfy this necessary condition and find that $\Tr \Ad$ is in the claimed range for all of them.

This quickly handles the cases of $G_2$ and $F_4$.  However, $E_7$ and $E_8$ are a bit more trouble.  When you use a laptop computer to solve for the critical points, it does not finish the computation immediately.  So we make the problem easier in these cases.  One idea is to use interpretations of the adjoint representation of these groups in terms of Hamming codes, but here is a less fancy approach leveraging Lemma \ref{St.lem}.

For $E_7$, we branch to $D_6 A_1$.  Leveraging the lemma, we can assume that the toral element in the $A_1$ component is actually the identity.  Then when we branch $D_6$ to $D_4 A_1 A_1$, we have a similar argument.  So we end up with two of the seven variables fixed before we do the elimination on the computer, and this makes the problem easy for the computer.  The case of $E_8$ is handled the same way except that the argument works for one additional level.

\subsection*{Highest short root representations}
For types $G_2$ and $F_4$, this same procedure can be applied to the highest short root representation  $\rho$ instead of $\Ad$, where we find:
\[
\begin{array}{cl} \toprule
\text{type of $G$}&\text{objective function $f$} \\ \midrule
G_2& t_1^2 + t_1 t_2 - 1\\
F_4& t_1 t_2 + t_1 t_3 + t_2 t_3 + t_1 t_4 + t_2 t_4 + t_3 t_4 + 2 \\ \bottomrule
\end{array} 
\]
Following this procedure again, we obtain a proof of Proposition \ref{short.prop} for these two cases.

\subsection*{Type \texorpdfstring{$E_6$}{E6}}
Suppose now that $G$ is simple compact of type $E_6$.  We aim to compute $\mmin(1)$ for this $G$.

The adjoint character is a fundamental character, so there
are no regular critical points (Lemma~\ref{St.lem}).

There is a unique orbit of codimension 1 tori fixed pointwise by a reflection.  These tori are
naturally maximal tori of subgroups of type $A_5$.  This occurs as
a component of a subgroup $A_5A_1$, and the corresponding
decomposition of $\mfe_6$ is:
\[
\mfa_5+\mfa_1+\wedge^3(V_6)\otimes V_2.
\]
The corresponding polynomial in the invariants $e_1,\ldots,e_5$ is nonlinear
$(e_1e_5+2+2e_3)$, but the derivative in $e_3$ is still a nonzero constant, and
thus a critical point must be a nonregular element of $A_5$.

Any nonregular semisimple element of $A_5$ is contained in a subgroup isomorphic to the double cover of $\GL_4$ obtained by adjoining a square root $f_2$ of the determinant (where the repeated eigenvalue is $f_2^{-1}$).  We implicitly define functions $a_0$, $a_1$, and $b_1$ by writing
the characteristic polynomial of the element of $\GL_4$ as
\[
t^4-a_1t^3+f_2a_0t^2-f_2b_1t+f_2^2.
\]
The decomposition of $\mfe_6$ is linear in $a_0$, so again there
are no regular critical points.

The next case is isogenous to $A_1\times \Gm^2$. We can then restrict $\Tr \Ad$ to the maximal torus in $A_1 \times \Gm^2$ and simply solve for the vanishing
of the gradient.

We provide Magma code for these computations in the GitHub repo for this paper.

\section{Type \texorpdfstring{$D$}{D}} \label{D.sec}

Take $G = \SO(2n)$, a compact connected group of type $D_n$ consisting of $X \in \GL_{2n}(\R)$ such that $X^\trans = X^{-1}$ and $\det X = 1$.  Consider the maximal torus consisting of diagonal matrices.  The entries of such a matrix $g$ come in pairs $s_j$, $s_j^{-1}$, and we consider the coordinate $t_j := s_j + 1/s_j$.  Then 
\begin{equation} \label{D.ad}
\Tr \Ad(g) = n + \sum_{j < k} t_j t_k = n + e_2(s_1, 1/s_1, s_2, 1/s_2, \ldots),
\end{equation}
where $e_2$ denotes the degree 2 elementary symmetric function.

We scale the variables $t_j$ defined above by 1/2, which scales $\sum t_i t_j$ by 4.  Plugging in the result of the following lemma gives the value of $\mmin(1)$ claimed by Theorem \ref{simple.thm}.

\begin{lem} \label{D.lem}
The minimum of the quadratic form $q(t_1, \ldots, t_n) := \sum_{1 \le i < j \le n} t_i t_j$ over $[-1, 1]^{\times n} \subset \R^n$ is $-n/2$ if $n$ is even and $(1-n)/2$ if $n$ is odd.
\end{lem}

\begin{proof}
Let $t := (t_1, \ldots, t_n)$ be a point where the quadratic form $q$ is minimal.  Let $a$, $b$ be the number of the $t_i$ equal to $1$ and $-1$ respectively and set $c := n - a - b$.

In the case $c = 0$, we have $2q(t) = a^2 - a + b^2 - b - 2ab = (a-b)^2 - n \ge -n$.  If $n$ is even, this is what we want.  If $n$ is odd, we have $a \ne b$, hence $(a-b)^2 \ge 1$ and we get $2q(t) \ge 1-n$.  Moreover, in both cases, equality is possible: Take $a = b = n/2$ if $n$ is even, and $a = (n-1)/2$, $b = (n+1)/2$ if $n$ is odd.

In the case $c \ge 1$, let $i$ be an index such that $t_i \ne \pm 1$.  Then $0 = \partial_{t_i} q\vert_t = v - t_i$ for $v = \sum t_i$.  This implies that $|v| \le 1$, $v = a - b + cv$, and
\[
2q(t) = v^2 - \sum_j t_j^2 = v^2 - a - b - cv^2 = (1 - n) + (c-1)(1 - v^2).
\]
Since $v^2 \le 1$, this shows that $2q(t) \ge 1- n$.
\end{proof}

\begin{rmk}
In the case where $n$ is even, we could argue as follows.
For a $g$ coming from the compact real form, we have $\bar{s}_j = 1/s_j$, so $s_j - 1/s_j$ is imaginary and 
\[
-\sum_j (s_j - 1/s_j)^2 \ge 0.
\]
Then
\[
0 \le \sum_j t_j^2 - \sum_j (s_j - 1/s_j)^2 = 2n + 2\Tr\Ad(g),
\]
proving that $\Tr \Ad(g) \ge -(\rank G)$ in this case.
\end{rmk}

\section{Type \texorpdfstring{$C$}{C}} \label{C.sec}

View the group $\Sp_{2n}$ of type $C_n$ as the stabilizer in $\GL_{2n}$ of a non-degenerate skew-symmetric bilinear form, specifically its complex points are the $X \in M_{2n}(\C)$ such that $X^\trans \Omega X = \Omega$ for 
\[
\Omega = \left( \begin{smallmatrix} 0 & I_n \\ -I_n & 0 \end{smallmatrix} \right),
\]
where $I_n$ denotes the $n$-by-$n$ identity matrix.  The compact real form of $\Sp_{2n}$ has real points the intersection with the unitary group.

Consider the maximal torus consisting of diagonal matrices.  The entries of such a matrix $g$ come in pairs $s_j$, $s_j^{-1}$, and we consider the coordinate $t_j := s_j + 1/s_j$.  Then 
\begin{equation} \label{C.ad}
\Tr \Ad(g) = -n + \sum_j t_j^2 + \sum_{j < k} t_j t_k.
\end{equation}

Note that the eigenvalues of an element of $\Sp_{2n}$ coming from the compact form have $\bar{s}_j = 1/s_j$, so each $t_j$ is real.  It follows that the power sums $p_d = \sum_i t_i^d$ have $p_1^2$ and $p_2$ positive.  The expression
$\sum t_j^2 + \sum t_j t_k$ is the degree 2 complete homogeneous symmetric function in the $t_j$'s, so it is $(p_1^2 + p_2)/2 \ge 0$.  Then \eqref{C.ad} gives the required bound $\Tr \Ad(g) \ge -(\rank G)$, proving Theorem \ref{simple.thm} for type $C$.

\subsection*{Highest short root representation}
The highest short root representation $\rho$ of $\Sp_{2n}$ is $(\wedge^2 \C^{2n}) / \C$.  Since this differs by $\C$ from the adjoint representation for type $D$, we see from \eqref{D.ad} that
\begin{equation} \label{C.rho}
\Tr \rho(g) = \sum_{i < j} t_i t_j + n - 1.
\end{equation}
The claim for type $C_n$ in Proposition \ref{short.prop} now follows directly from Lemma \ref{D.lem}.

\section{Type \texorpdfstring{$B$}{B}} \label{B.sec}

Take $G = \SO(2n+1)$, a compact connected group of type $B_n$ consisting of $X \in \GL_{2n+1}(\R)$ such that $X^\trans = X^{-1}$ and $\det X = 1$.  
Consider the maximal torus consisting of diagonal matrices.  The entries of such a matrix $g$ come in pairs $s_j$, $s_j^{-1}$, and we consider the coordinate $t_j := s_j + 1/s_j$.  Then 
\begin{equation} \label{B.ad}
\Tr \Ad(g) = n + \sum_i t_i + \sum_{i<j} t_i t_j.
\end{equation}
We aim to show that $\mmin(1) = -n$.

If the eigenvalue $-1$ occurs in $g$ (as an element of $\SO(2n+1)$), equivalently, $t_i = -2$ for some $i$, then we can pick a non-degenerate subspace $U$ of $\C^{2n+1}$ where $g$ acts as the diagonal matrix with entries $1$, $-1$.  Viewing $g$ as acting on $U$ and on the orthogonal complement of $U$, we can apply induction to deduce the claim.

By reordering the $t_i$, we may assume that $t_1, \ldots, t_{n-c} = 2$ for some $c$ and $t_i \ne \pm 2$ for $i > n-c$.  If $g$ is a critical point for $\Tr \Ad$, then, we have $\partial_{t_i} \Tr \Ad = 1 + v - t_i$ for 
$v := \sum_j t_j$. 
We divide into cases based on the value of $c$.  In the case $c = 0$, $g$ is the identity element and $\Tr \Ad(g) = \dim G$.  In the case $c = 1$, the equation $v = 2(n-c) + (1+v)c$ gives a contradiction.

So we are reduced to consider the case $c \ge 2$.  Rearranging the equation $\Tr \Ad(g) \ge -n$, we aim to show that 
\[
2 v + v^2 - (\sum t_i^2) \ge -4n.
\]
Now, $\sum t_i^2 = 4(n-c) - (1+v)^2 c$, so we aim to show that 
\[
v^2 + 2v + \frac{3c}{1-c} \le 0.
\]
Since $\frac{3c}{1-c} \le -3$, it suffices to verify that 
\[
v^2 + 2v - 3 \le 0.
\]
Which holds because $v = t_i - 1$ is in $[-3, 1]$.  This completes the proof of Theorem \ref{simple.thm} for type $B$.

\subsection*{The highest short root representation}
The highest short root representation $\rho$ of $\SO(2n+1)$ is the natural representation $\R^{2n+1}$, so
\begin{equation} \label{B.rho}
\Tr \rho(g) = 1 + \sum_{i} t_i.
\end{equation}
This is evidently minimized by taking each $t_i = -2$, so $\min \Tr \rho(g) = 1 - 2n$ as claimed in \autoref{short.root}, completing the proof of Prop.~\ref{short.prop}.

\part{Completion of proofs of Serre's claims} \label{completion.part}

We have verified the values of $\mmin(1)$ in Theorem \ref{simple.thm} for $G$ simple, and have done so in two different ways for the crux case of the exceptional groups.  Specifically, the value of $\mmin(1)$ has been computed for type $A$ in Example \ref{SUn.eg}, type $B$ in \S\ref{B.sec}, type $C$ in \S\ref{C.sec}, type $D$ in \S\ref{D.sec}, type $E$ in \cite[App.~B]{DMMV}, type $F$ in \S\ref{E8.sec}, and type $G$ in Example \ref{G2.eg}.  Alternatively, types $E$, $F$, and $G$ were handled in \S\ref{serre.method} and type $E_8$ was also treated in \S\ref{E8.sec}.

We now complete the proofs of the results stated in \S\ref{serre.bounds}.

\section{Components other than the identity} \label{outer.sec}

To complete the proof of Theorem \ref{simple.thm}, it remains to compute $\mmin(s)$ and $\mmax(s)$ when $s \ne 1$.  In the following results, we consider a compact Lie group $G$, and our interest is in the case where $G$ is not connected.

\subsection*{A subquotient \texorpdfstring{$H$}{H}}
Let $T$ be a maximal torus in $G$.
Since $T$ is abelian, for $w \in N_G(T)$ the function $T \to T$ defined by $t \mapsto [w, t] := (w t w^{-1}) t^{-1}$ is a homomorphism.  (This same notion appeared in the proof of Lemma \ref{omega.fiber}.)  We denote its image, a closed subgroup of $T$, by $[w, T]$.

\begin{lem} \label{wT.trace}
Let $w\in N_G(T)$, and let $\rho:G\to \GL(V)$ be a representation of $G$.  Then for all $t\in T$,
\[
\Tr(\rho(wt)) = \Tr(\rho(wt)|_{V^{[w,T]}}).
\]
\end{lem}

\begin{proof}
Note that since $[w,T]=[wT,T]=[wt,T]$, it suffices to prove the claim when $t=1$.

Over $\C$, $V$ is a direct sum of weight spaces $V_\lambda$ for weights $\lambda$.  If $w\cdot\lambda\ne \lambda$, then the action of $w$ permutes the weight spaces for weights in the orbit of $\lambda$, so those summands together contribute 0 to the trace.  We thus find that
\[
\Tr(\rho(w)) = \sum_{\lambda:w\cdot \lambda=\lambda} \Tr(\rho(w)|_{V_\lambda}) = \Tr(\rho(w)|_{V'}),
\]
where $V'$ is the direct sum of the weight spaces corresponding to $w$-invariant weights.  A weight is $w$-invariant if and only if $\lambda(w t w^{-1})=\lambda(t)$ for all $t\in T$, if and only if $\lambda(w t w^{-1}t^{-1}) = 1$ for all $t\in T$, if and only if $\lambda|_{[w,T]}=1$.  In other words, $V'$ is nothing other than the subspace of $V$ of vectors fixed by $[w,T]$, as required.
\end{proof}

In particular, when optimizing a character over a coset $wT$, we may restrict to the subspace $V^{[w,T]}$.  This is naturally a representation of a subquotient of $G$, namely $N_G([w,T])/[w,T]$. Call the identity component of this group $H$; equivalently, $H = N_G([w,T])^\ic / [w,T]$.  (This $H$ will be used in this section and reappear in \S\ref{finite.proof} as part of the proof of Cor.~\ref{thm.3p}.)

Certainly $T$ normalizes its subtorus $[w,T]$.  For $t \in T$ we have 
$w[w,t]w^{-1} = [w, wtw^{-1}] \in [w,T]$, so $wT$ is contained in $N_G([w,T])$ and every element of $wT$ normalizes $H$.

The Lie algebra of $N_G([w,T])$ is the $[w,T]$-invariant subspace of $\mfg$.  Its roots (i.e., the nonzero weights of $\mfg^{[w,T]}$) are the $w$-invariant roots of $G$, while its maximal torus is $T / [w,T] = H_0(\langle w\rangle;T)$.

\begin{lem} \label{H.centralize}
There is an element $w' \in wT$ that centralizes $H$.
\end{lem}

\begin{proof}
Write $\Ad$ for the homomorphisms
\[
H \to \Aut(H) \leftarrow N_G([w,T])
\]
that send an element to conjugation by that element.  
Since $\Ad(w)\in \Aut(H)$ fixes the maximal torus $T/[w,T]$ of $H$ elementwise, $\Ad(w)$ equals conjugation by an element of $T/[w,T]$ by the following lemma.   Since $\Ad(T) = \Ad(T/[w,T])$, we have $\Ad(w) = \Ad(t)$ for some $t \in T$, and $w' := wt^{-1}$ centralizes $H$ as required.
\end{proof}

Here is the promised lemma, which is implicit in standard references such as \cite[\S{IX.4.10}]{Bou:g7}.

\begin{lem}
Suppose the compact Lie group $G$ is connected.  If an automorphism $\phi$ of $G$ fixes the maximal torus $T$ elementwise, then there is a $t \in T$ such that $\phi(g) = tgt^{-1}$ for $g \in G$.
\end{lem}

\begin{proof}
Since $\phi$ normalizes $T$, it acts on $\Hom(T, S^1)$, and this action is trivial since $\phi$ fixes $T$ elementwise. Therefore, the automorphism of the Dynkin diagram of $G$ induced by $\phi$ is trivial and it follows from \cite[\S{IX.4.10}, Eq.~(16)]{Bou:g7} that there is some $x \in G$ such that $\phi(g) = xgx^{-1}$ for all $g \in G$.  Since $x$ centralizes $T$, it belongs to $T$.
\end{proof}

We now return to the setting where $G$ need not be connected.  When $w$ centralizes $H$, Lemma \ref{wT.trace} has a particularly nice consequence.

\begin{lem} \label{H.trace}
With $G$, $T$, $w$ as above and corresponding subquotient $H$, suppose that $w$ centralizes $H$.  Then
\[
\Tr(\Ad_G(wt))
=
\Tr(\Ad_H(t))+\Tr(\Ad_G(w)|_\mft)-\rank(H)
\]
for $t \in T$.
\end{lem}

\begin{proof}
Let $H^+$ be the normalizer of $[w,T]$, so that its Lie algebra is $\mfg^{[w,T]}$.  Then by Lemma \ref{wT.trace}
\[
\Tr(\Ad_G(wt)) = \Tr(\Ad_{H^+}(wt)).
\]
Since $H$ is a characteristic central quotient of the identity component of $H^+$, we may identify its Lie algebra with a direct summand of the Lie algebra of $H^+$, and thus we find that
\begin{align}
\Tr(\Ad_{H^+}(wt))
&=
\Tr(\Ad_H(wt)) + \Tr(\Ad_{H^+}(wt)|_{\Lie([w,T])})\notag\\
&=
\Tr(\Ad_H(wt)) + \Tr(\Ad_G(wt)|_{\Lie([w,T])}).\notag
\end{align}
Since $w$ centralizes $H$, it acts trivially on the Lie algebra, and thus
\[
\Tr(\Ad_H(wt)) = \Tr(\Ad_H(t)).
\]
Similarly,  $t$ centralizes $[w,T]$, so
\[
\Tr(\Ad_G(wt)|_{\Lie([w,T])})
=
\Tr(\Ad_G(w)|_{\Lie([w,T])}).
\]
Since $\mft$ is the direct sum of $\Lie([w,T])$ and the maximal torus of $H$ (on which $w$ acts trivially), we may also write this as
\[
\Tr(\Ad_G(w)|_{\Lie([w,T])})
=
\Tr(\Ad_G(w)|_\mft) - \rank(H). \qedhere
\]
\end{proof}

\begin{cor}\label{GH.reduction}
Let $G$, $T$, $wT$, and $H$ be as above.  Then
  \[
  \inf_{g\in wT} \Tr\Ad_G(g)
  =
\left(  \inf_{h\in H} \Tr\Ad_H(h) \right)
  +
  \Tr(\Ad_G(w)|_{\mft})-\rank(H).
  \]
  The same statement holds with $\inf$ replaced by $\sup$.
\end{cor}

\begin{proof}
We may replace $w$ with $wt$ for $t \in T$, and so (by Lemma \ref{H.centralize}) assume that $w$ centralizes $H$.  Then we have by Lemma \ref{H.trace}:
  \[
  \inf_{t\in T} \Tr\Ad_G(wt)
  =
  \inf_{t\in T} \left(\Tr\Ad_H(t)+ \Tr(\Ad_G(w)|_{\mft})-\rank(H) \right).
  \]
  Since $H$ is compact every $h \in H$ is conjugate to an element in the image of $T$, so 
    \[
  \inf_{t\in T} \Tr\Ad_H(t) = \inf_{h\in H} \Tr\Ad_H(h),
  \]
which  gives the desired result for $\inf$.  The same argument works verbatim with $\inf$ replaced by $\sup$.
\end{proof}

Consider now the case where $G$ is simple and $w$ acts via a diagram automorphism, relative to some pinning of $G^\ic$ with respect to the maximal torus $T$.  Then $\Ad(w)\vert_{\mft}$ is a
permutation matrix in the basis corresponding to simple roots, and thus its
trace equals the number of simple roots fixed by $w$.
To effectively apply Cor.~\ref{GH.reduction} in this case, we need to describe the subquotient $H$, which we now do.

\begin{lem} \label{fixed.system}
Let $G$ and $w$ be as in the preceding paragraph.  If $w \ne 1$, then 
the type of $H$ is described in
\autoref{H.table}.
\end{lem}
\begin{table}[hbtp]
{\rowcolors{2}{white}{gray!25}
\begin{tabular}{lclr} \toprule
type of $G$&order of $c$&type of $H$ & $\Tr(\Ad(w)\vert_{\mft})-\rank(H)$ \\ \midrule
$A_{2m}$ ($m \ge 1$)&2&$A_1^{\times m}$ & $-m$\\
$A_{2m-1}$ ($m \ge 2$)&2&$A_1^{\times m}$  & $1-m$\\
$D_n$ ($n \ge 4$)&2&$D_{n-1}$ & $-1$ \\
$D_4$&3&$A_2$ & $-1$ \\
$E_6$&2&$D_4$ & $-2$ \\ \bottomrule
\end{tabular}}
\caption{Subquotient $H$ of $G$ corresponding to a diagram automorphism
  $w \ne 1$} \label{H.table}
\end{table}

\begin{proof}
We verify this in each case.  For example, when $G$ has type $D_4$ we label the simple roots as in the diagram 
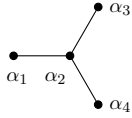
\[
\dynkin[label,label macro/.code={\alpha_{\drlap{#1}}},edge length=\edgelen] D4
\]
We shorten our notation by writing $abcd$ for the root $a \alpha_1 + b \alpha_2 + c \alpha_3 + d \alpha_4$.
if the automorphism $w$ has order 3, then the only positive roots fixed by $w$ are 0100, 1111, and 1211, which generate a root system of type $A_2$.  If the automorphism $w$ has order 2, say by interchanging $\alpha_3$ and $\alpha_4$, then the fixed positive roots are positive linear combinations of 1000, 0100, 0111.  Examining the inner products of these roots we find a root system of type $A_3 = D_3$.

We number the simple roots of types $E_6$, $D_n$, and $A_n$ according to the diagrams
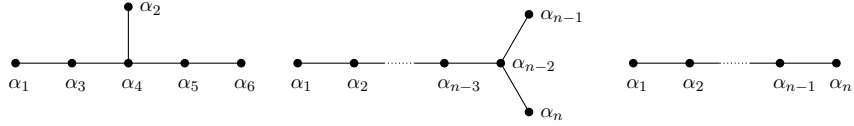
\begin{equation} \label{EDA}
\dynkin[label,label macro/.code={\alpha_{\drlap{#1}}},edge length=\edgelen] E6 
\quad
\dynkin[label,label macro/.code={\alpha_{\drlap{#1}}},edge length=\edgelen,labels={1,2,n-3,n-2,n-1,n},label directions={,,,right,,}] D{}
\quad
\dynkin[label,label macro/.code={\alpha_{\drlap{#1}}},edge length=\edgelen,labels={1,2,n-1,n}] A{}
\end{equation}
Recall that the support of a root $\sum n_i \alpha_i$ for $n_i \in \Z$ is the set of vertices of the Dynkin diagram with $n_i \ne 0$.  We use the facts that the support of root is a connected set of vertices, and that for every connected set of vertices $S$, the sum $\sum_{\alpha \in S} \alpha$ is a positive root \cite[\S{VI.1.6}, Cor.~3]{Bou:g4}. 

For type $E_6$ and $w$ of order 2, the fixed simple roots are spanned by 
\[
0100000, 000100, 001110, 101111,
\]
which span a root system of type $D_4$, with 010000 as the central vertex.  

For type $D_n$ ($n \ge 5$) and $w$ of order 2, the fixed simple roots are spanned by $\alpha_1, \ldots, \alpha_{n-2}$ and $\alpha_{n-2} + \alpha_{n-1} + \alpha_n$, which span a root system of type $D_{n-1}$ with terminal vertices $\alpha_1$, $\alpha_{n-2}$ and $\alpha_{n-2} + \alpha_{n-1} + \alpha_n$.

For type $A_n$ with $n \ge 2$ and $w$ of order 2, we divide into cases based on whether $n$ is even or odd.
If $n$ is even, i.e., $n = 2m$ for $m \ge 1$, the fixed positive roots are spanned by the roots 
\[
\alpha_m + \alpha_{m+1}, \quad \alpha_{m-1} + \alpha_m + \alpha_{m+1} + \alpha_{m+2}, \quad \ldots 
\]
If $n$ is odd, i.e., $n = 2m-1$ for $m \ge 2$, the fixed positive roots are spanned by the roots
\[
\alpha_m, \quad \alpha_{m-1} + \alpha_m + \alpha_{m+1}, \quad \ldots
\]
In either case, there are $m$ listed roots and the listed roots are pairwise orthogonal, so we find a root system of type $A_1^{\times m}$.
\end{proof}

\begin{proof}[Proof of \texorpdfstring{$s \ne 1$}{s not 1} case of Theorem \ref{simple.thm}]
Let $G$ be a compact Lie group that is simple and adjoint.  Pick a maximal torus $T$ of $G$ and a pinning of $G$ with respect to $T$.  The pinning identifies $\Aut(G)$ with $G \rtimes A$ for $A$ the group of automorphisms of the Dynkin diagram of $G$ \cite[\S{IX.4.10}, Cor.]{Bou:g7}, and in particular $A$ is identified with $\Out(G)$.  For $s \ne 1$ in $\Out(G)$, we have a corresponding element $w \in A$.  To prove Theorem \ref{simple.thm}, our task is to find the extreme values of $\Tr \Ad(gw)$ for $g \in G$.

A result of de Siebenthal \cite[Th.~on p.~57]{Sieb} says that the connected component $Gs$ of $\Aut(G)$ is a union of $G$-conjugates of $Tw$, so we are reduced to finding the extreme values of $\Tr \Ad(wt)$ for $t \in T$.  For this, we may apply Cor.~\ref{GH.reduction}, where the extreme values of $\Ad_H(h)$ are provided by the $s = 1$ case of Theorem \ref{simple.thm} (which we have already proved) and $\Tr(\Ad_G(w)\vert_\mft) - \rank(H)$ is provided by \autoref{H.table}.

For example, when $G$ has type $E_6$ and $w$ has order 2, $\Tr \Ad$ on
$D_4$ ranges from $-4$ to 28.  Since $\Tr(\Ad(w)|_{\mft})-\rank(H)=-2$,
we find that the $\Tr \Ad$ on $Gw$ varies from $-6$ to 26.

When $G$ has type $D_4$ and $w$ has order 3, $\Tr \Ad$ on $A_2$ ranges from $-1$ to $8$, so $\Tr \Ad$ on $Gw$ varies from $-2$ to 7.

When $G$ has type $D_n$ and $w$ has order 2, $\Tr \Ad$ on $D_{n-1}$ has
maximum $2n^2 - 5n + 3$,  so $\Tr \Ad$ on $Gw$ has maximum $2n^2 - 5n + 2$.  The minimum of $\Tr \Ad$ on $D_{n-1}$ is $-(n-1)$ if $n$ is odd and $2-(n-1)$ if $n$ is even, whence the minimum of $\Tr \Ad$ on $Gw$ is $-n$ if $n$ is odd and $2-n$ if $n$ is even.

When $G$ has type $A_n$ and $w$ has order 2, $\Tr \Ad$ on $A_1^{\times m}$
ranges from $-m$ to $3m$.  If $n = 2m$, then the correction term is $-m$,
so $\Tr \Ad$ on $Gw$ ranges from $-2m = -n$ to $2m = n$.  If $n = 2m-1$,
then the correction term is $1-m$, so $\Tr \Ad$ on $Gw$ ranges from $-2m+1
= -n$ to $2m+1 = n+2$.
\end{proof}

This completes the proof of Theorem \ref{simple.thm}.


\section{Proof of Theorem \ref{devissage.cor}} \label{devissage.sec}

Serre wrote that Theorem \ref{devissage.cor} follows from Theorem \ref{simple.thm} by ``an easy d\'evissage''.   The key case is where the identity component $G^\ic$ of $G$ is a simple compact Lie group, which is a special case of Theorem \ref{simple.thm}.  

Next suppose that $G^\ic$ is a product of $n$ copies of a simple compact Lie group $H$ for some $n$.  Writing $\Lie(G)$ as a sum of $n$ copies of $\Lie(H)$ suffices to write elements of $\End(\Lie(G))$ as block matrices with blocks in $\End(\Lie(H))$ and $n^2$ blocks.  For $g \in G$, writing $\Ad(g)$ in this way, we find that either $g$ normalizes a copy of $H$, in which case the corresponding block $X$ of $\Ad(g)$ belongs to the image of $\Ad(\Aut(H))$ and so $\Tr X \ge -(\rank H)$, or $g$ permutes some of the factors, in which case those terms contribute nothing to the trace of $\Ad(g)$.  In summary, we find that the conclusion holds for this $G$.

Next suppose that $G^\ic \cong H_1 \times H_2$, where $H_1$ and $H_2$ are semisimple adjoint compact Lie groups that have no simple factors in common.  Then every $g \in G$ normalizes both $H_1$ and $H_2$ and $\Tr \Ad(g) = \Tr \Ad(g\vert_{H_1}) + \Tr \Ad(g\vert_{H_2})$.

If $G$ is a semisimple adjoint compact Lie group, then it is a product of simple adjoint compact Lie groups, and we are done in this case by what has gone before.

Now consider the case where $G^\ic$ is a torus. Since $G$ is compact, every representation is unitary, so the eigenvalues of $\Ad(g)$, for $g \in G$, are all complex numbers of norm 1.  It follows that $\Tr \Ad(g) \ge -(\dim G) = -(\rank G)$, proving the claim in this case.

In the general case, $G$ is merely assumed to be compact, so it is reductive.  In particular, the center $Z$ of $G$ has identity component a torus.  Every $g \in G$ normalizes $Z$, and so $\Tr \Ad(g)\vert_Z \ge -(\rank Z)$.  Furthermore, the identity component of $G/Z$ is semisimple adjoint, and it follows that the action of $g$ on $\Lie(G/Z)$ also satisfies the trace bound.  Because $\rank Z + \rank G/Z = \rank G$, this completes the proof of Theorem \ref{devissage.cor}. \hfill $\qed$

\section{Proof of Corollary \ref{thm.3p}} \label{finite.proof}

Our aim in this section is to prove Corollary \ref{thm.3p}, Serre's Theorem $3'$.  Recall that $\mft$ is the Lie algebra of a maximal torus $T$ in a compact Lie group $G$.

\begin{lem} \label{mean.val}
Let the random element $w' \in G$ be chosen uniformly from a coset $wT$ of
$T$ in $N_G(T)$.  Then $\Tr(\Ad(w'))$ has mean $\Tr(\Ad(w)|_{\mft})$ and
variance the number of roots fixed by $w$.
\end{lem}

\begin{proof}
First suppose that $w = 1$.  Write the Lie algebra $\mfg$ of $G$ as a sum of weight spaces with respect to $T$: $\mfg = \mft \oplus \left( \oplus_{\alpha \in \Phi} \mfg_\alpha \right)$ where $\Phi := \{ \alpha \in T^* \mid \mfg_\alpha \ne 0 \}$ is the set of roots of $G$.  The restriction of each $\alpha$ to a rank 1 subtorus $S^1$ of $T$ is uniformly distributed, so has mean 0.  Since $T$ is isomorphic to a product of copies of $S^1$, it follows that $\Tr(\Ad_G(t) \vert_{\mfg_\alpha})$ has mean 0 and therefore that $\Tr(\Ad_G(t))$ has mean $\rank(G)$ as claimed.
Regarding the variance, $\Tr(\Ad_G(t))$ has Fourier expansion $\rank(G) + \sum_{\alpha \in \Phi} \alpha(t)$, so Parseval's Theorem gives
\[
\E\left[\left|\Tr(\Ad_G(t))-\rank(G)\right|^2\right] = \sum_{\alpha \in \Phi} \left| 1 \right|^2 = \left| \Phi \right|
\]
as claimed.

Recall the subquotient $H$ defined just before Lemma \ref{H.centralize}; it depends on the coset $wT$.  We next treat the case where $w$ centralizes $H$.
By Lemma \ref{H.trace}, we have
\begin{equation} \label{l13.3}
\Tr(\Ad_G(wt)) - \Tr(\Ad_G(w)\vert_\mft) = \Tr \Ad_H(t) - \rank(H)
\end{equation}
for $t \in T$.  Because $T$ surjects on to a maximal torus of $H$, we may apply the $w = 1$ case to the right side of \eqref{l13.3} to see that the left side has mean 0.  Neither $\Tr(\Ad_G(w)\vert_\mft)$ nor $\rank(H)$ depend on $t$, so the variance of $\Tr(\Ad_G(wt))$ is the variance of either side of \eqref{l13.3}, which is in turn the variance of $\Tr(\Ad_H(t))$, which by the $w = 1$ case equals the number of roots of $H$, i.e., the number of roots fixed by $w$, as claimed.

Finally we consider the general case.  By Lemma \ref{H.centralize}, there is some $w''$ in $wT$ such that $w''$ centralizes $H$.  Since $w''$ belongs to the same coset $wT$, we have $\Tr(\Ad(w)\vert_\mft) = \Tr(\Ad(w'')\vert_\mft)$ and $w$ and $w''$ fix the same roots of $G$.
\end{proof}

\begin{cor} \label{mean.cor1}
For any coset $wT$ of $T$ in $N_G(T)$, we have the inequalities
\[
\min_{w'\in wT} \Tr(\Ad(w')) \le \Tr(\Ad(w)\vert_{\mft})\le \max_{w'\in wT} \Tr(\Ad(w')),
\]
with each equality holding if and only if $w$ fixes no roots.
\end{cor}

\begin{proof}
The inequalities follow from connectivity of $T$ and the mean value theorem.  If $w$ fixes no roots, then the trace is constant on $wT$, so equality holds; if $w$ fixes a root, then the real random variable $\Tr\Ad(wt)$ has nonzero variance and thus attains values strictly below and above its mean.
\end{proof}

\begin{cor} \label{mean.cor2}
For any component $C$ of $G$ and $w \in N_G(T)\cap C$,
\[
\min_{g\in C}\Tr\Ad(g)\le \Tr(\Ad(w)\vert_{\mft})\le \max_{g\in C}\Tr\Ad(g).
\]
If either equality holds, then $w$ fixes no roots.
\end{cor}

\begin{proof}
Since $wT\subset C$, we have
\begin{multline*}
\min_{g\in C}\Tr\Ad(g)\le \min_{w'\in wT} \Tr(\Ad(w')) \le \Tr(\Ad(w)\vert_{\mft}) \\
\le \max_{w'\in wT} \Tr(\Ad(w'))\le \max_{g\in C}\Tr\Ad(g),
\end{multline*}
establishing the first claim.  For the second claim, note that if the given equality holds, then so does the corresponding equality in the previous corollary, and thus $w$ fixes no roots as required.
\end{proof}

\begin{lem} \label{lem.3p}
Let $x \in N_G(T)$.  There exists $w \in xN_{G^\ic}(T)$ that fixes no roots of $G$ with respect to $T$ and 
satisfies
\begin{equation} \label{3p.1}
\inf_{g \in xG^\ic} \Tr \Ad(g) = \Tr(w \vert_\mft).
\end{equation} 
If $G^\ic$ is simple, then $w$ can be chosen to additionally satisfy the following equivalent properties:
\begin{enumerate}[(a)]
\item \label{3p.finite} $Z_T(w)$ is finite.
\item \label{3p.comm} $[w, T] = T$.
\item \label{3p.eig} $\mft^{w} = 0$, i.e., $1$ is not an eigenvalue for $w$ acting on $\mft$.
\end{enumerate}
\end{lem}

\begin{proof}
We first prove the equivalence of the \ref{3p.finite}--\ref{3p.eig}.  For \ref{3p.finite}$\Leftrightarrow$\ref{3p.eig}, we have $\dim Z_T(w) = \dim \mft^{w}$.  For \ref{3p.comm}, as in Lemma \ref{wT.trace} the function $T \to T$ defined by $t \mapsto [w, t] = (w t w^{-1}) t^{-1}$ is a homomorphism with kernel $Z_T(w)$.  The image $[w, T]$ has dimension $\dim T - \dim Z_T(w)$, verifying the equivalence \ref{3p.finite}$\Leftrightarrow$\ref{3p.comm}.

Next we prove that any of the three equivalent conditions implies that $w$ fixes no roots of $G$.  Indeed, if $w$ fixes a root $\alpha$, then it fixes the corresponding coroot $\av$ in the dual root system, then the corresponding cocharacter $\Gm \to T$, which would imply $\mft^{w} \ne 0$.

\proofsec{Simple case}:
Suppose first that the identity component $G^\ic$ of $G$ is simple; this is the crux case.
If $-1 \in xN_{G^\ic}(T)/T$, we can take $w = -1$; it fixes no roots and $\Tr(w\vert_\mft) = -(\rank G)$ is the minimum of $\Tr \Ad(g)$ by Theorem \ref{simple.thm}\ref{simple.m1}.  So suppose $-1 \notin xN_{G^\ic}(T)/T$, and in particular $G^\ic$ has type $A_n$ for $n \ge 2$, $D_n$ for $n \ge 4$, or $E_6$.  As in the statement of Theorem \ref{simple.thm}, conjugation by $x$ induces a graph automorphism $s$ of the Dynkin diagram, and we refer to Theorem \ref{simple.thm} for the value of $\mmin(s)$.

If $G^\ic$ has type $A$, then it is isogenous to $\SU(n)$ for some $n \ge 3$ and $s = 1$.  Take 
$w$ to be any $n$-cycle 
in the Weyl group $S_n$.  The action of the Weyl group on the coroot space is equivalent to the subspace of the permutation representation of $S_n$ consisting of vectors whose coordinates sum to 1.  Any fixed vector in $\mft$ then must have all coordinates equal and therefore must be 0.  Since $w$ acts on the permutation representation with trace zero, it has  trace $-1 = \mmin(s)$ on $\mft$.

If $G^\ic$ has type $D_n$ for odd $n \ge 5$, then $s$ has order 1 and $\mmin(s) = 2 - n$.  We view $G$ as $\SO_{2n}$, which contains a copy of $\SO_6 \times \SO_{2n-6}$ of type $A_3 \times D_{n-3}$.   Take $w$ to be a 4-cycle in the Weyl group of $\SO_6 = \SL_4/\mu_2$ and $-1$ in the Weyl group of $\SO_{2n-6}$.  
The first part has trace $-1$ on its 3-dimensional subalgebra of $\mft$ and the second has trace $3-n$ on its $(n-3)$-dimensional subalgebra of $\mft$,
for a total trace of  $2 - n = \mmin(s)$.  By the same argument as in the previous paragraph, $w$ does not fix any vector in $\mft$. 

If $G^\ic$ has type $D_n$ for even $n \ge 4$ and $s$ has order 2, then after multiplying $x$ by an element of $G^\ic$ we may assume that conjugation by $x$ normalizes the set of simple roots and so is a graph automorphism of the Dynkin diagram order 2.  Composing conjugation by $x$ with the Weyl group element that is $-1$ on the torus gives an element $w$ that evidently fixes no roots and has $\Tr(w\vert_\mft) = 2 - n = \mmin(s)$.

If $G^\ic$ has type $D_4$ and $s$ has order 3, then the subgroup of $N_G(T)/T$ generated by $xT$ and the Weyl group has order 576 and contains two conjugacy classes of elements of order 3 that have exactly 8 elements.  One of these two classes is in $xN_{G^\ic(T)}/T$.  Any representative of that class provides an element $w$ that acts on the root space $\mft$ with trace $-2 = \mmin(s)$ and does not have 1 as an eigenvalue. 

If $G^\ic$ has type $E_6$, take $w$ to be an element of order 3.  There are three conjugacy classes of such elements in the Weyl group, see for example \cite[Table 9]{Carter:weyl}.  We take $w$ to be in the unique conjugacy class of size 80.  In the action on the root space $\mft$, $w$ has trace $-3 = \mmin(s)$ and does not have 1 as an eigenvalue.

The computations in the last two paragraphs can be easily verified using Magma's \texttt{ReflectionGroup} and \texttt{ConjugacyClass} functions.

\proofsec{Remaining cases}:
We have now handled all cases where $G^\ic$ is simple.  If the identity component of $G^\ic$ is a torus, then there is nothing to prove.

For the general case, conjugation by $x$ permutes the simple factors of $G^\ic$.  Put $G^\ic \cong G_1 \times G_2$, where $G_1$ is a product of the central torus and the simple factors that are normalized by $x$ and $G_2$ is the a product of the simple factors that are not normalized by $x$.  Now every $g \in xG^\ic$ has $\Tr \Ad(g)$ equal to the trace of the projection of $\Ad(g)$ on $\Lie(G_1)$ and every $w \in xN_{G^\ic}(T)/T$ moves every root in $G_2$ and $\Tr(w \vert_{\Lie(T)})$ equals the trace of the restriction of $w$ to $\Lie(T \cap G_1)$.  So it suffices to apply the lemma for the case $G^\ic = G_1$ and extend the resulting $w$ in any way to an element of $xN_{G^\ic}(T)/T$.
\end{proof}

\begin{proof}[Proof of Corollary \ref{thm.3p}]
Take $x \in G$ minimizing $\Tr \Ad(x)$.  Then $x$ normalizes some maximal torus $T$ of $G$ \cite[\S{IX.5.3}, Th.~1(b)]{Bou:g7}.  Apply Lemma \ref{lem.3p}
to produce a $w \in xG^\ic$ with $\Tr \Ad(x) = \Tr(w \vert_\mft)$, so
\[
\min_{g \in G} \Tr \Ad(g) \ge \min_{w \in W} \Tr(w\vert_{\mft}).
\]

On the other hand, for every component $C$ of $G$, 
\[
\Tr \Ad(x) \le \min_{g \in C} \Tr \Ad(g) \le \min_{w \in N_G(T) \cap C} \Tr \Ad(w\vert_\mft),
\]
where the last inequality is by Corollary \ref{mean.cor2}.  So
\[
\Tr \Ad(x) \le \min_{w \in N_G(T)} \Tr \Ad(w\vert_\mft). \qedhere
\]
\end{proof}

\part{Minimum values of characters} \label{min.part}

\section{Minimum values of characters: \texorpdfstring{$\SU(2)$}{SU(2)}} \label{SU2.sec}

Let's now discuss the minimum value of $\RE \chi(g)$ for $\chi := \Tr \rho$ the character of a representation $\rho$ of a compact group $G$.  As explained in the introduction, for elementary reasons we have $-\chi(1) \le \RE \chi(g)$, and Serre's results show that this bound can be far from sharp.  So far we have considered the cases $\rho = \Ad$ or the highest short root representation.

Perhaps the most natural collection of representations to consider next are the absolutely irreducible representations of $\SU(2)$.  Put $\chi_d$ for the character of the irreducible (complex) representation $\Sym^d(\C^2)$ with highest weight $d$, i.e., $s^d + s^{d-2} + \cdots + s^{-(d-2)} + s^{-d}$.  In the notation of Example \ref{A1.eg} and \eqref{cheb.small}, we have $t = \chi_1 = s + s^{-1}$ and $\chi_d \in \C[T]^W = \C[t]$.  That is, $\chi_d$ can be written as a polynomial in $t$, which we denote by $\chi_d^t$, where we have already written the first few cases in \eqref{cheb.small}.  It is easy to see that 
\begin{equation} \label{chi.ind}
\chi^t_d = t \chi^t_{d-1} - \chi^t_{d-2}, 
\end{equation}
so $\chi^t_d = U_d(t/2)$ for $U_d$ the Chebyshev poynomial of the second kind of degree $d$.
Therefore, $\chi^t_d$ has $d$ real zeros in $(-2, 2)$, namely $2\cos(2\pi k/(d+1))$ for $1\le k\le d$.  It follows that its derivative 
has $d-1$ real zeros in $(-2,2)$ and thus the critical points are those 
zeros together with $\pm 2$, so all lie in $[-2,2]$.

The inductive formula \eqref{chi.ind} quickly gives $\chi^t_d(\pm 2) = (\pm 1)^d (d+1)$, so the minimum of $\chi_d$ is $-(d+1)$ if $d$ is odd.  Alternatively, the fact that the representation maps $-1 \in \SU(2)$ to $-1 \in \GL(\Sym^d(\C^2))$ for $d$ odd shows that the minimum of $\chi_d$ is $-\chi_d(1) = -(d+1)$.

The case of even $d$ is more interesting.  The critical points and extreme values of $\chi_d^t$ are symmetric around $t = 0$, and the minimum (for $d$ even) occurs at the critical points closest to $\pm 2$ \cite{Cheb}.  

\begin{eg} \label{SU2.chi}
For $d = 6$, the critical points of $\chi^t_d$ closest to $\pm 2$ are at 
\[
t_0 = \pm \sqrt{(5 + \sqrt{7})/3} \approx \pm 1.596
\]
where
\[
\chi^t_d(t_0) = -\frac{7}{27} \left(1+2 \sqrt{7}\right) \approx  -1.631.
\]
In particular, the minimum of $\chi_d$ is neither rational nor an algebraic integer.
\end{eg}

We can provide a uniform bound.  

\begin{prop} \label{SU2.neg}
For $d > 0$, 
\[
\min_{g \in \SU(2)} \chi_d(g) \le -\sinfty \chi_d(1) \quad \text{and} \quad \limsup_{d\to\infty} \min_{g\in\SU(2)} \frac{\chi_d(g)}{\chi_d(1)}
= -\sinfty,
\]
where 
\[
\sinfty := -\min_{\theta \in \R} \frac{\sin(\theta)}{\theta} \approx 0.2172.
\]
\end{prop}

\begin{proof}
We first prove the inequality.  Put $\chi^s_d$ for $\chi_d$ written as a function of $s$, so
\[
\chi^s_d(s) = \frac{s^{d+1}-s^{-d-1}}{s-1/s}
\]
and
\[
\chi^s_d(e^{i\theta}) = \frac{\sin((d+1)\theta)}{\sin(\theta)}.
\]
Trivially
\begin{equation} \label{SU2.neg.1}
\min_{g} \frac{\chi_d(g)}{\chi_d(1)} = \min_{\theta\in \R} \frac{\chi_d^s(e^{i\theta})}{\chi_d^s(1)} \le \frac{\sin(\thinfty)}{(d+1)\sin(\thinfty/(d+1))}.
\end{equation}
where $\thinfty\approx 4.493$ is the positive argument that minimizes $\sin(\theta)/\theta$.

Multiplying by $\theta/\theta$,
\[
\frac{\sin((d+1)\theta)}{(d+1)\sin\theta} = \frac{\sin((d+1)\theta)/((d+1)\theta)}{\sin(\theta)/\theta}.
\]
Since $\left| \sin (\theta) / \theta \right| \le 1$, it follows that 
\[
\left| \frac{\sin((d+1)\theta)}{(d+1)\sin(\theta)} \right| \ge \left| \frac{\sin((d+1) \theta)}{(d+1)\theta} \right|.
\]
The expression on the right side of \eqref{SU2.neg.1} is negative, so plugging in $\thinfty/(d+1)$ for $\theta$ we find
\[
\min_{g} \frac{\chi_d(g)}{\chi_d(1)} \le  \frac{\sin(\thinfty)}{\thinfty}.
\]

As for the claimed limit, we have
\[
\lim_{d\to\infty}\frac{\chi^s_d(e^{it/d})}{\chi_d(1)} = \frac{\sin(t)}{t}.
\]
Since this is a limit of analytic functions that converges uniformly on
compacta, the derivatives also converge, and so do the zeros of the
derivatives.  We thus conclude that
\[
\limsup_{d\to\infty} \min_g \frac{\chi_d(g)}{\chi_d(1)} = -\sinfty,
\]
as claimed.
\end{proof}

\section{What kind of numbers are these?} \label{what}

For $G$ a compact  Lie group, what sort of real number might $\min_{g \in G} \RE \chi(g)$ be?  When $G$ is simple and $\chi = \Tr \Ad$ or the trace of the highest short root representation, we showed that it is an integer (Theorem \ref{simple.thm} and Proposition \ref{short.prop}).  For the same $G$ and general $\chi$, since elements of finite order are dense in $G$, it follows that the minimum is a limit of algebraic integers.  Examples \ref{F4.eg}, \ref{E8.eg}, and \ref{SU2.chi} exhibit cases where $\min_{g \in G} \RE \chi(g)$ is  an algebraic number, but not an algebraic integer.  We have the following general result.

\begin{prop} \label{which.prop}
Let $\rho$ be a representation of a compact Lie group $G$.  For each connected component $C$ of $G$, the extreme values of $\RE \Tr \rho$ on $C$ are algebraic numbers.
\end{prop}

\begin{proof}
We first set up some notation.  Put $\Qalg$ for the algebraic closure of $\Q$ and $F := \Qalg \cap \R$, the algebraic closure of $\Q$ in $\R$.
The unit circle $U(1)$, which we have denoted elsewhere by $S^1$ and by which we mean the group scheme over $\R$ whose $\C$-points are the norm 1 complex numbers, is obtained by base change from a group scheme over $F$ we denote by $U(1)_F$, whose $F$-points are the norm 1 elements of $\Qalg$.

Pick an element $c \in C$.  By the result of de Siebenthal from \cite{Sieb} that we used before in the proof of Theorem \ref{simple.thm} in \S\ref{outer.sec}, there is a torus $T$ in $G$ such that $C$ is a union of $G^\ic$-conjugates of $cT$.  Therefore, 
 the extreme values of $\Re \Tr \rho$ are obtained on $cT$.  Replacing $G$ by the subgroup generated by $c$ and $T$, we may assume that $G$ is abelian and $G^\ic$ is a torus, i.e., a product $U(1)^r$ of $r$ copies of $U(1)$ for some $r$.

We claim that $G$, viewed as an algebraic group, is defined over the algebraic closure $F$ of $\Q$ in $\R$.  
For $m$ the order of $c$ in $G$, we have an exact sequence
\[
\begin{CD}
1 @>>> K @>>> U(1)^r \times \Z/m @>{(1, \ldots, 1, 1) \mapsto c}>> G @>>> 1
\end{CD}
\]
where $U(1)^r$ is identified with $G^\ic$.  The elements of the kernel $K$ have coordinates in $U(1)$ that are elements of order dividing $m$, and in particular these elements already belong to $U(1)_F$, verifying the claim.
This argument also shows that the connected component $C$ of $G$ is defined over $F$.

Next we prove that $\rho$ is also defined over $F$.
Since $G$ is abelian, $\rho$ is a sum of one-dimensional representations,  i.e., it amounts to a homomorphism $G \to U(1)^n$ for some $n$.  Its restriction to $T$ is a homomorphism $U(1)^r \to U(1)^n$, such that composing this map with projection on any of the $n$ coordinates gives a map of the form $(z_1, z_2, \ldots, z_r) \mapsto \prod_{j=1}^n z_j^{e_j}$ for integers $e_j$ \cite[\S{IX.7}]{Bou:g7}. 
If we re-write this product, expanding $z_j = x_j + iy_j$ for $x_j, y_j$ real, we find that the real and imaginary parts of $\prod z_j^{e_j}$ are each given by a polynomial in the $x_j$ and $y_j$ with integer coefficients.  That is, $\rho\vert_T$ arises by base change from a polynomial function $U(1)_F^r \to U(1)_F^n$.  The projection of $\rho(c)$ in each of the $n$ coordinates is an $m$-th root of unity, so $\rho(c)$ is an element of $U(1)_F^n$.

The trace provides a morphism of schemes from $U(1)_F^n$ (the norm 1 elements of $\Qalg$) to the Weil restriction $\Res_{\Qalg/F}(\A^1_F)$, whose $F$-points equal the algebraic closure $\Qalg$ of $\Q$.  (The trace does not respect the group operation so this is a morphism of schemes but not group schemes.)  This variety is naturally identified with pairs of elements $(z, \bar{z})$ for $z \in \Qalg$, and the map 
\[
\pi(z, \bar{z}) := \tfrac12(z + \bar{z})
\]
is a morphism of schemes $\pi \: \Res_{\Qalg/F}(\A^1_F) \to \A^1_F$.  The composition $\pi \Tr \rho$ provides a morphism of schemes $G \to \A^1_F$ whose base change to $\R$ is the map we have denoted elsewhere by $\Re \Tr \rho$.  We restrict this map to the connected component $C$ of $G$.  Its graph 
is a closed subscheme of $C \times \A^1_F$.
The projection of the graph on the second component --- i.e., the image of $\Re \Tr \rho$ --- is a semi-algebraic subset of $\A^1_F$ by the Tarski-Seidenberg Projection Theorem \cite[Th.~1.5.8]{Scheiderer}. 

Let us now consider the image of $\Re \Tr \rho$ over $\R$.  Because $C$ is compact, its image is a finite union of closed intervals in $\R$.  In particular the boundary of the image is a finite set of points including the extreme values.  The Tarski Transfer Principle implies then that the boundary of the image over the subfield $F$ is a semi-algebraic set over $F$ \cite[Cor.~1.6.18]{Scheiderer}.  Therefore the extreme values are algebraic over $F$, proving the claim.
\end{proof}

The second and third paragraphs of the proof, showing that $G$ and $\rho$ are defined over $F$, could be replaced by various other arguments.  For example, one could use the classification of compact Lie groups to express $G$ as obtained by Galois descent from a reductive complex algebraic group whose identity component is split.  Repeating this same descent for the quadratic extension $\Qalg$ of $F$ shows that $G$ is defined over $F$, and a similar argument would show the same for $\rho$.

\section{Minimum values of characters: general case} \label{general.sec}

In the examples we have considered so far, $\min_{g \in G} \RE \chi(g)$ has been negative, which is easily seen to be a general phenomenon:

\begin{prop} \label{neg.val}
Let $\chi$ be the character of an irreducible (complex) representation of a
compact Lie group $G$.  
\begin{enumerate}[(a)]
\item \label{neg.neg} If $\chi$ is nontrivial, then there exists $g\in G$ such that $\RE(\chi(g))<0$.
\item \label{neg.0} If $\chi(1) > 1$, then there exists $g \in G$ such that $\RE \chi(g) = 0$.
\end{enumerate}
\end{prop}

\begin{proof}
Since $\chi$ is nontrivial, orthogonality implies that the 
integral of $\chi$ relative to Haar measure is 0, and thus so is the 
integral of $\RE\chi$.  Since the integral of a nonzero nonnegative 
continuous function is positive, \ref{neg.neg} follows.

Now consider \ref{neg.0}.    The identity component $G^\ic$ is a normal subgroup of $G$, so its fixed subspace is an invariant submodule.  It follows that $G^\ic$ either acts trivially or without fixed points.
If $G^\ic$ acts trivially, then the action of $G$ factors through the finite group $G/G^\ic$, which is not cyclic (because $\chi$ is irreducible and $\chi(1) > 1$), so \ref{neg.0} holds by Burnside \cite[p.~115]{Burnside:tr}.  If $G^\ic$ acts without fixed points, then the average value of $\RE \chi$ on $G^\ic$ is zero; since $\RE \chi(G^\ic)$ is connected and takes a positive value $\chi(1)$, it must also take a negative value, and we conclude \ref{neg.0}.
\end{proof}

For self-dual representations, the statement of \ref{neg.val}\ref{neg.0} says that there exists a $g \in G$ such that $\chi(g) = 0$, which
sounds similar to Theorem 1 in \cite{Se:trace}.  Serre's statement adds the extra conclusion that $g$ can be chosen to be of finite order and does not require the self-dual hypothesis.  See \cite{Se:zeros} for a proof of that result.

The aim of this section is to strengthen \ref{neg.val}\ref{neg.neg} by proving the following uniform 
bound on $\RE\chi / \chi(1)$:
\begin{thm} \label{neg.thm}
  For any connected compact Lie group $G$, there is a constant $c_G>0$ such that for
  any nontrivial, irreducible character $\chi$, there exists $g\in G$ such that
  $\RE\chi(g) \le -c_G\chi(1)$.
\end{thm}

We treated the case $G = \SU(2)$ in Prop.~\ref{SU2.neg} by explicit calculation and found the optimal value of $c_{\SU(2)}$.  (That calculation suggests that there is little hope in making the optimal $c_G$ explicit for general $G$.)
As we saw in that proof, the claim in the theorem
is closely related to the behavior of the characters themselves as the
weight becomes large, suggesting that we should consider limits of the form
\begin{equation} \label{X.lim}
\lim_{n\to\infty}
\frac{\chi_{\lambda^{(n)}}(e^{it/d_n})}{\chi_{\lambda^{(n)}}(1)}
\end{equation}
in which the sequence $\lambda^{(n)}/d_n$ converges, where $\lambda^{(n)}$ is a dominant weight and $d_n \in \R$ is positive.   
When the limit of
$\lambda^{(n)}/d_n$ is in the interior of the fundamental chamber, this
limit is straightforward to compute, and leads us to define the function
\begin{equation} \label{Xdef}
X(s,t):=
\left( \prod_{r\in \Phi^+} \frac{r\cdot \rho}{(r\cdot s)(r\cdot t)} \right)
\left( \sum_{w\in W} \sgn(w) e^{s\cdot wt} \right)
\end{equation}
for $\Phi^+$ the positive roots and $s,t$ in the Cartan algebra $\mft(\C)$ of the {\em complex} form of $G$.  By convention, we take the
split form, so that the Lie algebra of $T$ is $i{\mft}(\R)$.  In the right side of \eqref{Xdef}, the products $r \cdot s$ and $r \cdot t$ are canonically determined by the root system.  Since $G$ is compact, there is a $G$-invariant inner product on $\mfg$ \cite[\S{IX.1.3}, Prop.~1]{Bou:g7}; pick one and use it to calculate $s \cdot wt$.  The inner product defines an isomorphism of $\mft(\C)$ with its dual, and therefore also an inner product among elements of the dual which we use to evaluate $r \cdot \rho$.  Up to scaling one of the arguments by $i$, $X$ is the Fourier transform of the Duistermaat-Heckman measure \cite{EtingofRains}.

Proposition \ref{167} below makes explicit the connection between \eqref{X.lim} and \eqref{Xdef}.

We will now take the analytic continuation of $X$,
 which has the 
advantage that it lets us extend the
function to the reflection hyperplanes.

\begin{prop} \label{Xprop}
The function $X(s,t)$ extends to an entire function of $\mft(\C)^2$
that satisfies the symmetries 
\begin{enumerate}[(a)]
\item $X(s,t)=X(t,s)$, 
\item \label{Xprop.scale} $X(s,ut)=X(us,t)$ for $u\in \C$, 
\item $X(s,wt)=X(s,t)$ for $w\in W$, and 
\item $X(0,0)=1$.
\end{enumerate}
\end{prop}

The proof will use the following lemma, which is implicit already in the development of the Weyl character formula.  It is a special case of a more general statement where $\C^n$ is replaced by a connected complex manifold $X$ and the reflection $r$ is replaced by an automorphism of $X$ of finite order.

\begin{lem} \label{value}
Let $r\in \GL_n(\C)$ be a reflection, and let $f(z)$ be an $r$-invariant meromorphic function on an $r$-invariant region in $\C^n$
containing points fixed by $r$.  Then the valuation of $f$ along the hyperplane $r\cdot z=0$ is even.
\end{lem}

\begin{proof}
This claim is clearly basis-independent, so we may assume that $r$ fixes the first $n-1$ coordinates and negates the last coordinate, so that we are claiming that $f$ has even order along the hyperplane $z_n=0$.
Let $z'$ be a fixed point of $r$ in the region where $f$ is defined, such that no other polar divisor of $f$ contains $z'$.  Then $f$ has a Laurent series expansion in some polydisc around $z'$ of the form
\[
f(z) = \sum_{m_1,\dots,m_{n-1}\ge 0,m_n\ge v} c_{\vec{m}} \prod_{1\le i\le n-1} (z_i-z'_i)^{m_i} z_n^{m_n},
\]
where $v$ is the valuation, so that there is at least one nonzero term with $m_n=v$.  But this must be invariant under negating $z_n$, so that any monomial in the expansion must have $m_n$ even, implying that $v$ is even as well.
\end{proof}

\begin{proof}[Proof of Proposition \ref{Xprop}]
  The function $X(s,t)$ is certainly analytic on the complement of the
  reflection hyperplanes $r\cdot s = 0$ and $r\cdot t = 0$, where it clearly
  satisfies the symmetries $X(s,t)=X(t,s)$; $X(s,ut)=X(us,t)$ for $u\ne 0$;
  and $X(s,wt)=X(s,t)$.

In particular, these same symmetries hold for the function $X(s,t)$, which
by inspection is a meromorphic function with at most a simple pole along
any reflection hyperplane.  But the Lemma \ref{value} tells us that the order of such a pole must be even, and thus $X(s,t)$ actually has nonnegative valuation along any reflection hyperplane.  Thus $X(s,t)$ is a meromorphic function with trivial polar divisor, so is holomorphic.

The symmetries continue to hold on the entire continuation, and thus in
particular we may take the limit $u\to 0$ in $X(s,ut)=X(us,t)$ to obtain
$X(s,0)=X(0,t)$ for all $s$ and $t$, and thus in particular both are equal to
$X(0,0)$.

  It remains only to verify that $X(0,0)=1$.  From the degenerate scaling
  symmetry, this is equal to $X(\rho,0)$ for $\rho$ the half-sum of the positive roots, and thus to $\lim_{t\to 0}
  X(\rho,t)$.  We moreover have the product formula
  \[
  X(\rho,t) = 
   \prod_{r\in \Phi^+} \frac{e^{r\cdot t/2}-e^{-r\cdot t/2}}{r\cdot t}
  \]
  in which each factor converges to $1$, so that $X(\rho,0)=1$ as required.
\end{proof}

\begin{rem}
  Note that $X$ only depends on the Weyl group, even though the root system appears in the definition of $X$.
  Rescaling the roots (including
  swapping short and long roots in non-simply-laced cases)  only
  multiplies $X$ by a constant, which must be 1 by virtue of $X(0,0)=1$.
\end{rem}

Since $X(s,t)$ differs by simple factors from the Weyl numerator, we can
express the normalized characters in terms of $X(s,t)$.

\begin{lem} \label{X.dom}
  For any $t\in \mft(\R)$ and any dominant weight $\lambda$, one has
  \[
  \frac{\chi_{\lambda}(e^{it})}{\chi_{\lambda}(1)}
  =
  \frac{X(\lambda+\rho,it)}{X(\rho,it)}.
  \]
\end{lem}

\begin{proof}
Writing out $X(\lambda+\rho, it)/X(\rho, it)$ using the definition \eqref{Xdef}, gives $\chi_\lambda(e^{it})$ divided by
\[
  \prod_{r\in \Phi^+} \frac{r\cdot (\lambda+\rho)}{r\cdot \rho} = \chi_\la(1),
  \]
where the equality is the Weyl dimension formula.
\end{proof}

This lets us show that $X(s,t)$ is indeed the limiting character even on
the boundary of the fundamental chamber.

\begin{prop} \label{167}
 Let $(\lambda^{(n)},d_n)$ be a sequence of pairs where each $\lambda^{(n)}$ is a dominant weight and $d_n \in \R$ is such that
  $\lambda^{(n)}/d_n$ converges and $d_n \to \infty$.  Then
  \[
  \lim_{n\to\infty}
  \frac{\chi_{\lambda^{(n)}}(e^{it/d_n})}{\chi_{\lambda^{(n)}}(1)}
  =
  X\left(\lim_{n\to\infty} \lambda^{(n)}/d_n,it \right)
  \]
  for all $t\in \mft(\R)$.
\end{prop}

\begin{proof}
  From Lemma \ref{X.dom}, one has
  \begin{align*}
  \lim_{n\to\infty}
  \frac{\chi_{\lambda^{(n)}}(e^{it/d_n})}{\chi_{\lambda^{(n)}}(1)}
& =
  \lim_{n\to\infty}
  \frac{X(\lambda^{(n)}+\rho,it/d_n)}{X(\rho,it/d_n)} \\
&=
  \lim_{n\to\infty}
  \frac{X(\lambda^{(n)}/d_n+\rho/d_n,it)}{X(\rho/d_n,it)}
  \end{align*}
  with the second equality following from the scaling symmetry of $X$.
  Since $X$ is continuous, we can pull the limit inside, at which point
  the claim follows by noting that $X(0,it)=1$.
\end{proof}

To establish the desired bound on $\min\RE\chi/\chi(1)$, we will
show after some work that $X(s,it)$ takes negative values for any nonzero $s\in \mft(\R)$ (Lemma \ref{t.neg}).
We write $\lambda$ for the Haar measure on $\mft(\R)$, since it is essentially just Lebesgue measure.
To avoid issues with the normalization of the Haar measure, the root
system, and the $G$-invariant inner product, define a scaling constant
\[
Z := \int_{t\in \mft(\R)} \left( \prod_{r\in \Phi^+}(r\cdot t)^2 \right) e^{-t\cdot
t/2}
\lambda(dt).
\]

We also use the following mild restatement of the well-known formula for the Fourier/Laplace transform of a Gaussian:
\begin{lem} \label{fourier}
For any $s\in \mft(\C)$, one has
\[
\int_{t \in \mft(\R)} e^{s\cdot t} e^{-t\cdot t/2} \lambda(dt)
=
e^{s\cdot s/2} \int_{t \in \mft(\R)} e^{-t\cdot t/2} \lambda(dt).
\]
\end{lem}

\begin{proof}
The integral on the left is absolutely convergent for any $s\in \mft(\C)$, and uniformly so given any bound on $\left|\RE(s)\right|$, so describes an entire function of $s$.  Since the right-hand side is manifestly entire, it suffices to verify that the two sides agree when $s$ is real.  But for $s\in \mft(\R)$, we may complete the square on the left by translating $t$ by $s$ to obtain
\begin{align*}
\int_{t\in \mft(\R)} e^{s\cdot t} e^{-t\cdot t/2} \lambda(dt)
&=  
\int_{t\in -s+\mft(\R)} e^{s\cdot s/2} e^{-t\cdot t/2} \lambda(dt)\\
&= 
e^{s\cdot s/2} \int_{t\in \mft(\R)} e^{-t\cdot t/2} \lambda(dt)
\end{align*}
as required.
\end{proof}

\begin{prop} \label{158}
  For $s,s'\in \mft(\C)$, one has
  \begin{multline*}
  \frac{1}{Z}
  \int_{t\in \mft(\R)} X(s,it)X(s',it)
  \left( \prod_{r\in \Phi^+} (r\cdot t)^2 \right)
  e^{-t\cdot t/2}
  \lambda(dt)
  = \\
  e^{-s\cdot s/2-s'\cdot s'/2}
  X(s,s').
  \end{multline*}
\end{prop}

\begin{proof}
  Since the right-hand side is analytic in $s$, $s'$, it suffices to show
  the equality when neither of $s$ or $s'$ lies in any reflection
  hyperplane.  We can then write
  \begin{multline*}
  \int_{t\in \mft(\R)} X(s,it)X(s',it)
 \left( \prod_{r\in \Phi^+} (r\cdot t)^2\right)
  e^{-t\cdot t/2}
  \lambda(dt)
  = \\
  \left(\prod_{r\in \Phi^+} \frac{(r\cdot \rho)^2}{(r\cdot s)(r\cdot s')}\right)
  \int_{t\in \mft(\R)}
  \sum_{w,w'\in W}
  \sgn(ww')
  e^{i(ws+w's')\cdot t}
  e^{-t\cdot t/2}
  \lambda(dt).
  \end{multline*}
  Applying Lemma \ref{fourier} gives:
  \begin{multline*}
  \int_{t\in \mft(\R)}
  \sum_{w,w'\in W}
  \sgn(ww')
  e^{i(ws+w's')\cdot t}
  e^{-t\cdot t/2}
  \lambda(dt)
  = \\
  \sum_{w,w'\in W}
  \sgn(ww')
  e^{-(ws+w's')\cdot(ws+w's')/2}
  \int_{t\in \mft(\R)}
  e^{-t\cdot t/2}
  \lambda(dt).
  \end{multline*}
  We then observe that
  \begin{multline*}
  \sum_{w,w'\in W}
  \sgn(ww')
  e^{-(ws+w's')\cdot (ws+w's')/2}
  \\
  = e^{-s\cdot s/2}e^{-s'\cdot s'/2}
  \sum_{w,w'\in W} \sgn(ww') e^{ws\cdot w's'}
  \\
  = |W| \,
  e^{-s\cdot s/2} e^{-s'\cdot s'/2}
  \left( \prod_{r\in \Phi^+}\frac{(r\cdot s)(r\cdot s')}{r\cdot \rho} \right)
  X(s,s').
  \end{multline*}
  We thus conclude that
  \[
  \int_{t\in \mft(\R)} X(s,it)X(s',it)
 \left( \prod_{r\in \Phi^+} (r\cdot t)^2 \right)
  e^{-t\cdot t/2}
  \lambda(dt)
  \]
  is proportional to
  \[
    e^{-s\cdot s/2}e^{-s'\cdot s'/2}
  X(s,s'),
  \]
  and can recover the constant by setting $s=s'=0$.
\end{proof}

The proposition is already enough to give us a weak form of orthogonality sufficient
to prove the desired negativity result.

\begin{cor} \label{limit.0}
  If $\RE(s\cdot s)>0$, then
  \[
  \lim_{u\to\infty}
  \int_{t\in \mft(\R)} X(s,it)
 \left( \prod_{r\in \Phi^+} (r\cdot t)^2 \right)
  e^{-t\cdot t/2u}
  \lambda(dt)
  =
  0.
  \]
\end{cor}

\begin{proof}
  Setting $s'=0$ and abbreviating $p(t) := \prod_{r\in \Phi^+} (r\cdot t)^2$ in Proposition \ref{158} gives
  \[
  \frac{1}{Z}
  \int_{t\in \mft(\R)} X(s,it)
 \,p(t)\,
  e^{-t\cdot t/2}
  \lambda(dt)
  =
  e^{-s\cdot s/2},
  \]
  and thus
  \[
  \frac{1}{Z}
  \int_{t\in \mft(\R)} X(u^{1/2}s,it)
\,  p(t)\,
  e^{-t\cdot t/2}
  \lambda(dt)
  =
  e^{-u(s\cdot s)/2}.
  \]
  We can use the scaling symmetry of $X$ to move the $u^{1/2}$ to the other
  argument and then do a change of variables to move it to the other parts
  of the integrand.  We thus find
  \begin{multline*}
  \frac{1}{Z}
  \int_{t\in \mft(\R)} X(u^{1/2}s,it)
\,p(t)\,
  e^{-t\cdot t/2}
  \lambda(dt)
  = \\
  u^{-|\Phi^+|/2-\dim\mft/2}
  \frac{1}{Z}
  \int_{t\in \mft(\R)} X(s,it)
\, p(t) \,
  e^{-t\cdot t/2u}
  \lambda(dt),
  \end{multline*}
  and thus
  \[
  \frac{1}{Z}
  \int_{t\in \mft(\R)} X(s,it)
 \, p(t) \,
  e^{-t\cdot t/2u}
  \lambda(dt)
  = 
  u^{(|\Phi^+|+\dim\mft)/2}
  e^{-u(s\cdot s)/2},
\]
  which converges to $0$ as required.
\end{proof}

\begin{lem} \label{t.neg}
  For any nonzero $s\in \mft(\R)$, there exists $t\in \mft(\R)$
  such that $\RE X(s,it)<0$.
\end{lem}

\begin{proof}
  Suppose otherwise, and consider the integral
  \[
  \frac{1}{Z}
  \int_{t\in \mft(\R)}
  \RE X(s,it) \,
  \Big(\prod_{r\in \Phi^+} (r\cdot t)^2 \Big)
  e^{-t\cdot t/2u}
  \lambda(dt)
  \]
  as a function of $u$.  On the one hand, the integrand is nonnegative,
  not identically 0, and monotonically increasing in $u$, and thus the
  integral itself must be a positive, monotonically increasing, function of
  $u$.  On the other hand, Corollary \ref{limit.0} tells us that the
  integral converges to 0, giving a contradiction.
\end{proof}

\begin{proof}[Proof of Theorem \ref{neg.thm}]
  If $G$ has positive-dimensional center $Z$, then for any character $\chi$
  which is nontrivial on $Z^\ic$, there is an element $z\in Z^\ic$ such that
  $\chi(z)=-\chi(1)$.  Since one has the universal bound $\RE(\chi(g))\ge
  -\chi(1)$, it follows that any constant that works for $G/Z^\ic$ will work
  for $G$ as well.  We thus reduce to the case that $G$ is semisimple.

In the semisimple case, replacing $G$ by its universal cover only makes the conclusion harder to satisfy, and thus we may assume $G$ simply connected.  Moreover, if $G=G_1\times G_2$, then any irreducible character of $G$ factors as $\chi_1\boxtimes \chi_2$; it follows that if the theorem holds for both factors with constants $c_1$, $c_2$, then it holds for $G$ with constant at least $\min(c_1,c_2)$.  (Take $g$ to be the identity on one of the two factors.)  We thus reduce to the case that $G$ is simply connected and simple.

  Now, if the claim fails for $G$, then there exists a sequence $\lambda^{(n)}$
  of dominant weights such that
  \[
  \limsup_{n\to\infty} \inf_{g\in G}
  \frac{\RE\chi_{\lambda^{(n)}}(g)}{\chi_{\lambda^{(n)}}(1)}\ge 0.
  \]
  We may replace this by a subsequence to replace the limsup by a limit.
  Moreover, since spheres are compact, we may pass to a further subsequence
  to ensure that the sequence $\lambda^{(n)}/|\lambda^{(n)}|$ converges,
  say to $s$.  We then have by Prop.~\ref{167}, for any $t\in \mft(\R)$
  \[
  \lim_{n\to\infty}
  \frac{\chi_{\lambda^{(n)}}(e^{it/|\lambda^{(n)}|})}{\chi_{\lambda^{(n)}}(1)}
  =
  X(s,t),
  \]
  and by Lemma \ref{t.neg} may choose $t$ so that $\RE X(s,t)<0$.  But this gives a
  contradiction:
  \[
  \inf_{g\in G}
  \frac{\RE\chi_{\lambda^{(n)}}(g)}{\chi_{\lambda^{(n)}}(1)}
  \le
  \frac{\RE \chi_{\lambda^{(n)}}(e^{it/|\lambda^{(n)}|})}{\chi_{\lambda^{(n)}}(1)}
  \]
  but the right-hand side has negative limit while the left-hand side has
  nonnegative limit by assumption.
\end{proof}

\section{Extensions}

It is natural to wonder whether there is a more universal version of Theorem \ref{neg.thm} in which the dependence of the constant on $G$ is removed.  This cannot be done, as we have already seen that for $\chi=\Tr\Ad$ on $\SU(n)$, the minimum of $\RE\chi/\chi(1)$ is $-1/(n^2-1)$ (Example \ref{SUn.eg}) and thus $c_{\SU(n)}$ approaches $0$ as $n\to\infty$.  This suggests a more refined question to produce an explicit function $F$ on the natural numbers such that $c_G$ in Theorem \ref{neg.thm} is asymptotically $F(\rank(G))$ as $\rank(G) \to \infty$.  Note that this is really a question about the classical groups.

\begin{eg}
Pavel Etingof explained to us that for $G = \SU(n)$, 
\[
X(\rho, it) \ge -(4/\pi^2)^{n-2},
\]
see \cite[Prop.~6.1]{EtingofRains}.
Defining a sequence $\lambda^{(j)} := j \rho$ and $d_j = j$ in $j$ and applying Proposition \ref{167}, we find that $c_{\SU(n)}$ decays at least exponentially in $n$.
\end{eg}

\begin{conj}
For any nontrivial irreducible character $\chi$ of $G$, one has
\[
\min_{g\in G} \frac{\RE\chi(g)}{\chi(1)}
\le
\sup_{s\in \mft(\R)}
\inf_{t\in \mft(\R)}
X(s,it),
\]
and thus the optimal constant above is
\[
c_G = \sup_{s\in \mft(\R)}
\inf_{t\in \mft(\R)}
X(s,it),
\]
and in particular depends only on the Weyl group of $G$.
\end{conj}

A slightly weaker, but more explicit version has recently been shown in \cite[Cor.~1.4]{EtingofRains}, namely that
\[
\min_{g\in G} \frac{\RE\chi_\lambda(g)}{\chi_\lambda(1)}
\le
\inf_{t\in \mft(\R)} \RE X(\lambda+\rho,it)
\]
for all but finitely many dominant weights $\lambda$. 

\appendix
\section{Application to finite groups and fields of prime characteristic} \label{finite.sec}

Serre’s bounds on the trace of elements of a compact Lie group can also be used to give bounds on the values of the Brauer character for some representations of finite groups over fields of prime characteristic.  Here is one such result:

\begin{prop} \label{Brauer.adj}
Let $G$ be an adjoint simple algebraic group over a finite field $k$ of characteristic $p \ne 2$.  Then the Brauer character of $G(k)$ acting on $\Lie(G)$ is real-valued and its minimum value is $\mmin(1)$ as given by Theorem \ref{simple.thm}, unless perhaps $p = 3$ and $G$ is of type $A_2$ or $E_6$.
\end{prop}

Here is how we define the Brauer character, compare \cite[\S18.1]{Se:rep} or \cite[\S5.3]{Benson1}.  Let $\Gamma$ be a finite group and $k$ a finite field of characteristic $p \ne 0$.  Write the exponent of $\Gamma$ as $\gamma p^r$ for $\gamma$ not divisible by $p$.  Pick a ring $\cOh \subset \C$ containing the primitive root of unity $\exp(2\pi i/\gamma)$ and with a surjection $\cOh \to \kh$ where $\kh$ is an algebraic extension of $k$.  The surjection gives an isomorphism between the group of $\gamma$-th roots of unity in $\C$ and those in $\kh$ (because $p$ does not divide $\gamma$), and for each such root of unity $\zeta \in \kh$, we write $\tilde{\zeta}$ for the corresponding root of unity in $\C$.

Now, given a representation $\rho \: \Gamma \to \GL_{n,k}$, every $p$-regular element $g \in \Gamma$ (i.e., element $g \in \Gamma$ whose order divides $\gamma$) has $\rho(g)$ similar over $\kh$ to a diagonal matrix with diagonal entries roots of unity $\zeta_1, \ldots, \zeta_n$.  The \emph{Brauer character} $b_\rho$ is defined to be $b_{\rho}(g) := \tilde{\zeta}_1 + \cdots + \tilde{\zeta}_n \in \C$.

\subsection*{Brauer characters and group schemes}
This notion applies also to group schemes through their $k$-points.  Let $\cO \subset \C$ be a ring with a surjection $\cO \to k$ and let $G$ be an affine group scheme of finite type over $\cO$.  For every homomorphism of group schemes $\rho \: G \to \GL_{n,\cO}$, we obtain by base change a representation $\rho_k \: G(k) \to \GL_{n,k}$ of the finite group $G(k)$ and so it makes sense to speak of the Brauer character of $\rho_k$ using the definitions from the preceding two paragraphs where $\cOh$ is chosen to contain $\cO$.

In the special case where $G$ is a split torus, the map $G(\cOh) \to G(\kh)$ gives an isomorphism on the subgroups of elements of order dividing $\gamma$.  For such a $g \in G(\cO)$, the Brauer character $b_{\rho_k}(g)$ equals $\Tr \rho_{\C}(g)$.  

In the next result, we compare the Brauer character $b_{\rho_k}$ for $g \in G(k)$ versus the usual character of the representation $\rho_\C \: G(\C) \to \GL_{n,\C}$.  

\begin{lem}  \label{Brauer.ineq}
Suppose $G$ is a split connected reductive group over $\cO$ and put $K$ for a maximal compact subgroup of $G(\C)$.  Then
\[
\min \{\Re b_{\rho_k}(g) \mid \text{$p$-regular $g \in G(k)$})\} \ge \min \{\Re \Tr \rho_{\C}(g) \mid g \in K \}.
\]
If $\Tr((\rho_\C)\vert_K)$ is real-valued, then so is $b_{\rho_k}$.
\end{lem} 

\begin{proof}
Suppose $g \in G(k)$ is $p$-regular.  We aim to show that the image of $g \in G(\kh)$ is the image of an element $\gh \in G(\cOh)$ of the same order.  Then $b_{\rho_k}(g) = \Tr \rho_{\C}(\gh)$.  Since elements of finite order in $G(\C)$ are conjugate to elements of $K$, this will suffice to prove the claim.

The group $G$ is split over $\cO$ and the splitting data for $G$ includes a choice of split maximal torus $T$ defined over $\cO$.  If $g$ is in $T(\kh)$, then as explained before the statement of the lemma, $g$ is the image of an element  of $T(\cOh)$ of the same order and we are done.

For every representation of the base change $G \times k$ of $G$ to $k$, the image of $g$ has minimal polynomial dividing $x^\gamma - 1$ and in particular the minimal polynomial has distinct roots, so $g$ is semisimple.  It follows that $g$ is contained in a maximal torus $T'$ of $G$ defined over $k$.  Enlarging $\cOh$ and $\kh$ if necessary, we may assume that $\kh$ splits $T'$.  Therefore, $T' \times \kh$ is conjugate under $G(\kh)$ to $T \times \kh$.  The value of $b_{\rho_k}(g)$ only depends on the conjugacy class of $g$ in $G(\kh)$, so we may replace $g \in G(\kh)$ by its conjugate lying in $T(\kh)$ and the proof is complete.
\end{proof}

In the definition of Brauer character for a finite group, the representation of $\Gamma$ over $k$ is not required to have an analogous representation over $\C$, whereas such is required in Lemma \ref{Brauer.ineq}.  Such analogous representations exist in many cases of interest, as we now explain.  Given a root datum in the sense of \cite[\S{C.d}]{Milne} or \cite[\S{XXI.1.1}]{SGA3.3:new}, there is up to isomorphism a unique split connected reductive group $G_\Z$ over $\Z$ with that root datum \cite[XXV.1.2]{SGA3.3:new}.  Given a field $k$, one may pick an $\cO$ and define the group scheme $G$ over $\cO$ to be the base change $G_\Z \times \cO$.   For each dominant weight $\lambda$ that is a character on the maximal torus $T$ of $G$, there are various representations of $G$ that are defined over $\Z$ with highest weight $\lambda$ such as the Weyl module $V(\lambda)$, see \cite[II.8.3]{Jantzen}; the base change of $V(\lambda)$ to $\C$ is the irreducible representation of $G(\C)$ with highest weight $\lambda$.

\subsection*{The case of an algebraically closed field}  Continue the setup where $k$ is finite and $G$ is defined over $\cO$, but now choose $\kh$ to be an algebraic closure of $k$.  Since the group scheme $G$ is of finite type, for each $g \in G(\kh)$, there is a finite extension $k_g$ of $k$ contained in $\kh$ such that $g$ comes from $G(k_g)$, and in particular $g$ has finite order.  If $g$ is $p$-regular, then $\kh$ contains a primitive root of unity of the same order as $g$, and the preceding definitions provide a well-defined Brauer character of $g$.  Lemma \ref{Brauer.ineq} also holds in this situation, where the inequality is strengthened to an equality
\begin{equation} \label{Brauer.eq}
\inf \{\Re b_{\rho_{\kh}}(g) \mid \text{$p$-regular $g \in G(\kh)$})\} = \inf \{\Re \Tr \rho_{\C}(g) \mid g \in K \}.
\end{equation}
To see this, note that the elements of $K$ of finite order not divisible by $p$ are dense in $K$.  So it suffices to show that each such element is conjugate to an element coming from $G(\cOh)$, which follows by the same sorts of argument as in the proof of the lemma.

The ``inf'' on the right side of \eqref{Brauer.eq} can be replaced by ``min'' if and only if the minimum is obtained at an element of $K$ of finite order not divisible by $p$.

Let us now return to the adjoint representation and the proof of Prop.~\ref{Brauer.adj}.

\begin{proof}[Proof of Proposition \ref{Brauer.adj}]
That the Brauer character is real-valued is contained in Lemma \ref{Brauer.ineq}.

\proofsec{The split case}: Suppose first that $G$ is split over $k$, and view it as obtained by base change from the split group over $\cO$ with the same root datum, which we also denote by $G$.
 In view of Lemma \ref{Brauer.ineq}, it suffices to exhibit a $p$-regular element $g \in G(k)$ whose Brauer character equals $\mmin(1)$.

One approach to this is via the Weyl group.  
Since 2 is invertible in $\cO$, there is a finite subgroup $W^* \subseteq N_G(T)(\cO)$ that is an extension of the Weyl group (as a finite group) by a product of copies of $\Z/2$ and there is a section $\nu \: W \to W^*$, see for example \cite{Ti:etendu}, \cite[\S{IX.4}, Exercise 12]{Bou:g7}, or \cite{DemTits}. 
If $-1$ is in the Weyl group, take $w = -1$.  If $G$ has type $D_n$ for odd $n \ge 5$ or type $E_6$, then we take $w \in W$ to be the element of order 4 or 3 respectively from the proof of Lemma \ref{lem.3p}.  In each of these cases, $w$ fixes no root, $\nu(w)$ has order not divisible by $p$, and $\Tr(w\vert_{\mft}) = \mmin(1)$.  Thus $\nu(w) \in G(\cO)$ maps to a $p$-regular element of $G(k)$ with Brauer character $\mmin(1)$.

Now suppose $G$ has type $A_{n-1}$ with $n \ge 3$, i.e., $G = \PGL_n$.  We construct an element $g$ as the image in $G(k)$ of a block diagonal matrix $B$ with diagonal blocks of size 2, 3, or 5, with minimum polynomial $t^2 + 1$, $t^3 - 1$, or $t^5 - 1$ respectively (e.g., the companion matrix of one of these polynomials).  If $n = 2m$ is even, we take $B$ to have all 2-by-2 blocks.  If $n = 2m+1$ and $p \ne 3$, then we take $B$ to have a single 3-by-3 block and the rest 2-by-2.  If $n = 2m+1$ and $p = 3$, then $n \ge 5$, and we take $B$ to have a single 5-by-5 block and the rest 2-by-2.  Note that in each case $B$ has trace zero so $\Tr \Ad(B) = -1$ as in Example \ref{SUn.eg}.

\proofsec{Non-split case}: Suppose now that $G$ is not split. Since $k$ is finite, $G$ is quasi-split \cite[Prop.~17.99]{Milne}; 
has type $^2A_n$ for $n \ge 2$, $^2D_n$ for $n \ge 4$, $^3D_4$, or $^2E_6$; 
and is split by an extension $K$ of $k$ of degree 2, 2, 3, or 2 respectively.  
Phrased differently, $G$ is obtained by twisting the split adjoint group $G_0$ over $k$ of the same type as $G$ by a semilinear automorphism we denote by $\sigma$ that combines a graph automorphism of order $[K:k]$ of $G_0$ with a non-trivial $k$-algebra automorphism of $K$.  Since the minimum value of the Brauer character on $G(k)$ is no smaller than that on $G(K)$, it suffices to verify that an element of $G(K)$ whose Brauer character equals $\mmin(1)$ belongs to $G(k)$, i.e., is stable under the Galois action.

For type $A_{n-1}$ with $n \ge 3$, we imitate the split case, where we take a $d$-by-$d$ diagonal block of the matrix $B$ to belong to $\SU(d)$.  View $\SU(d)$ as the group of norm 1 isometries of the $K/k$-hermitian form on $K^d$ defined by 
$(x, y) \mapsto \sum_{i=1}^d \sigma(x) y$.  The companion matrix of $t^d + (-1)^d$ belongs to this interpretation of $\SU(d)$ and realizes a $d$-cycle in the Weyl group.

If $-1$ is in the Weyl group, then (as the longest element of the Weyl group relative to the choice of simple roots) it is fixed by the graph automorphism.  Since it fixes no roots and $G$ is adjoint, the image $g$ of $\nu(-1)$ in $G(K)$ does not depend on the choice of $\nu$, ergo $\sigma$ fixes $g$.  (We will repeat this argument below, but in this special case we could have alternatively referred to \cite[Prop.~XXIV.3.16.2(iii)]{SGA3.3:new}.)

For $G$ of type $E_6$, let $x$ be the Weyl group element from the proof of Lemma \ref{lem.3p}.  Since the conjugacy class of $x$ is the unique one of that size, there is a $y \in W$ such that $\sigma(x) = yxy^{-1}$.  Put $g$, $h$ for the images of $\nu(x)$, $\nu(y)$ in $G(k)$.  Since $g$ normalizes $T$, so does $\sigma(g)$, which acts as $yxy^{-1}$ on $T$ and preserves no roots, so $\sigma(g) = hgh^{-1}$.  By Lang's Theorem \cite[Cor.~17.97]{Milne}, there is an element $h_0 \in G(\kalg)$ for $\kalg$ the algebraic closure of $\kh$, such that $h = \sigma(h_0) h_0^{-1}$.  It follows that
$h_0^{-1} g h_0$ is fixed by $\sigma$, i.e., it belongs to $G(k)$.  Since it acts as $x$ on $T$, this case is complete.

In the remaining case  where $G$ has type $^2D_n$ for odd $n \ge 5$, we imitate the construction in the split case but noting that the 4-cycle in the Weyl group of $\SU(4)$ exists as described three paragraphs above.
\end{proof}

The hypothesis that $p \ne 2$ in Proposition \ref{Brauer.ineq} is not totally necessary.  For example, the proof provided for types $A_2$ and $E_6$ uses an element of order 3 in each case, so it works without change when $p = 2$.

\subsection*{Ree and Suzuki groups}
One could also ask for an analogue of Proposition \ref{Brauer.ineq} for the Ree and Suzuki groups.  Only one of those occurs in characteristic different from 2, namely $^2G_2(k)$ as the subgroup of the $k$-points $G_2(k)$ of the split algebraic group of type $G_2$ when $k$ is a finite field of order $3^n$ for odd $n \ge 1$, see \cite[Ch.~12--14]{Carter:simple}.  As a subgroup, $^2G_2(k)$ acts on the adjoint representation of the algebraic group $G_2$ and by the proposition the Brauer character of 3-regular elements of $^2G_2(k)$ will be real and $\ge -2$.  Conversely, the element $-1$ in the Weyl group of $G_2(k)$ (which can also be viewed as an element of $G_2(k)$ and belongs to the unique conjugacy class of elements of order 2 in $G_2(k)$ \cite[\S4]{Jac:comp}) is in $^2G_2(k)$, so this minimum is obtained on $^2G_2(k)$.
 

\begin{thebibliography}{DMMV21}

\bibitem[Ada85]{Adams:E8}
J.F. Adams, \emph{The fundamental representations of {$E\sb 8$}}, Conference on
  algebraic topology in honor of Peter Hilton (Saint John's, Nfld., 1983),
  Contemp. Math., vol.~37, Amer. Math. Soc., Providence, RI, 1985, pp.~1--10.

\bibitem[BCP90]{Magma}
W.~Bosma, J.~Cannon, and C.~Playoust, \emph{The {M}agma algebra system {I}: the
  user language}, J. Symbolic Computation \textbf{9} (1990), 677--698.

\bibitem[Ben91]{Benson1}
D.J. Benson, \emph{Representations and cohomology {I}: Basic representation
  theory of finite groups and associative algebras}, Cambridge Studies in
  Advanced Mathematics, no.~30, Cambridge University Press, 1991.

\bibitem[Bor91]{Borel}
A.~Borel, \emph{Linear algebraic groups}, second ed., Graduate Texts in
  Mathematics, vol. 126, Springer-Verlag, New York, 1991.

\bibitem[Bou89]{Bou:g1}
N.~Bourbaki, \emph{Lie groups and {L}ie algebras: Chapters 1--3},
  Springer-Verlag, Berlin, 1989.

\bibitem[Bou02]{Bou:g4}
\bysame, \emph{Lie groups and {L}ie algebras: Chapters 4--6}, Springer-Verlag,
  Berlin, 2002.

\bibitem[Bou05]{Bou:g7}
\bysame, \emph{Lie groups and {L}ie algebras: Chapters 7--9}, Springer-Verlag,
  Berlin, 2005.

\bibitem[Bur03]{Burnside:tr}
W.~Burnside, \emph{On an arithmetical theorem connected with roots of unity and
  its application to group characteristics}, Proc. London Math. Soc. \textbf{1}
  (1903), 112--116.

\bibitem[Car72]{Carter:weyl}
R.W. Carter, \emph{Conjugacy classes in the {W}eyl group}, Compositio Math.
  \textbf{25} (1972), no.~1, 1--59.

\bibitem[Car89]{Carter:simple}
\bysame, \emph{Simple groups of {L}ie type}, Wiley, 1989.

\bibitem[CCN{\etalchar{+}}85]{atlas}
J.H. Conway, R.T. Curtis, S.P. Norton, R.A. Parker, and R.A. Wilson,
  \emph{Atlas of finite groups}, Oxford University Press, Eynsham, 1985.

\bibitem[CG21]{ChayetG}
M.~Chayet and S.~Garibaldi, \emph{A class of continuous non-associative
  algebras arising from algebraic groups including ${E}_8$}, Forum of
  Mathematics: Sigma \textbf{9} (2021), e6.

\bibitem[Con10]{Conrad:cpt}
B.~Conrad, \emph{Answer to ``{C}lassification of (compact) {L}ie groups''},
  available at \url{https://mathoverflow.net/questions/6079}, February 2010.

\bibitem[DG11]{SGA3.3:new}
M.~Demazure and A.~Grothendieck, \emph{Sch{\'{e}}mas en groupes {III}:
  {S}tructure des sch\'{e}mas en groupes r\'{e}ductifs}, Soci\'et\'e
  Math\'ematique de France, 2011, re-edition edited by P. Gille and P. Polo.

\bibitem[DMMV21]{DMMV}
M.~{Dutour Sikiri\'c}, D.A. Madore, P.~Moustrou, and F.~Vallentin,
  \emph{Coloring the {V}oronoi tessellation of lattices}, J. London Math. Soc.
  (2) \textbf{104} (2021), 1135--1171.

\bibitem[dS56]{Sieb}
J.~de~Siebenthal, \emph{Sur les groupes de {L}ie compacts non connexes},
  Commentarii Mathematici Helvetici \textbf{31} (1956), 41--89.

\bibitem[DV21]{DeMedtsVC}
T.~{De Medts} and M.~{Van Couwenberghe}, \emph{Non-associative {F}robenius
  algebras for simply laced {C}hevalley groups}, Trans. Amer. Math. Soc.
  \textbf{374} (2021), no.~12, 8715--8774.

\bibitem[EOR07]{EOR}
P.~Etingof, A.~Oblomkov, and E.~Rains, \emph{Generalized double affine {H}ecke
  algebras of rank 1 and quantized del {P}ezzo surfaces}, Adv. Math.
  \textbf{212} (2007), no.~2, 749--796.

\bibitem[ER24]{EtingofRains}
P.~Etingof and E.~Rains, \emph{Bounds for the asymptotic characters of simple
  {L}ie groups}, arXiv:2405.10341, May 2024.

\bibitem[FG98]{FreyGriess}
D.D. Frey and R.L. {Griess, Jr.}, \emph{Conjugacy classes of elements in the
  {B}orovik group}, J. Algebra \textbf{203} (1998), 226--243.

\bibitem[Fre01]{Frey}
D.D. Frey, \emph{Conjugacy of $\mathrm{Alt}_5$ and $\mathrm{SL}(2,5)$ subgroups
  of ${E}_7(\mathbb{C})$}, J. Group Theory \textbf{4} (2001), 277--323.

\bibitem[GG15]{GG:simple}
S.~Garibaldi and R.M. Guralnick, \emph{Simple groups stabilizing polynomials},
  Forum of Mathematics: Pi \textbf{3} (2015), e3 (41 pages).

\bibitem[{Gri}95]{Griess:G2}
R.L. {Griess, Jr.}, \emph{Basic conjugacy theorems for ${G}_2$}, Invent. Math.
  \textbf{121} (1995), 257--277.

\bibitem[Gui07]{Guillot:lambda}
P.~Guillot, \emph{The representation ring of a simply connected {L}ie group as
  a {$\lambda$}-ring}, Comm. Algebra \textbf{35} (2007), no.~3, 875--883.

\bibitem[Hal15]{Hall:2nd}
B.C. Hall, \emph{{L}ie groups, {L}ie algebras, and representations: an
  elementary introduction}, 2nd ed., Graduate Texts in Mathematics, no. 222,
  Springer, 2015.

\bibitem[Jac58]{Jac:comp}
N.~Jacobson, \emph{Composition algebras and their automorphisms}, Rendiconti
  del Circolo Matematico di Palermo. Serie II \textbf{7} (1958), 55--80, [\#60
  in \emph{Coll.\ Math.\ Papers}, vol.~2].

\bibitem[Jan03]{Jantzen}
J.C. Jantzen, \emph{Representations of algebraic groups}, second ed., Math.
  Surveys and Monographs, vol. 107, Amer. Math. Soc., 2003.

\bibitem[Kac69]{Kac:aut}
V.G. Kac, \emph{Automorphisms of finite order of semisimple {L}ie algebras},
  Funct. Anal. Appl. \textbf{3} (1969), 252--254.

\bibitem[Kai06]{Kaiser}
N.~Kaiser, \emph{Mean eigenvalues for simple, simply connected, compact {L}ie
  groups}, J. Phys. A: Math. Gen. \textbf{39} (2006), 15287--15298.

\bibitem[Kat04]{Katz:G2}
N.M. Katz, \emph{Notes on ${G}_2$, determinants, and equidistribution}, Finite
  Fields and their Applications \textbf{10} (2004), 221--269.

\bibitem[Kni73]{Knighten}
C.M. Knighten, \emph{Differentials on quotients of algebraic varieties}, Trans.
  Amer. Math. Soc. \textbf{177} (1973), 65--89.

\bibitem[Leh16]{Lehalleur}
S.P. Lehalleur, \emph{Subgroups of maximal rank of reductive groups}, Autour
  des Sch\'emas en Groupes III, vol.~47, Soci\'et\'e Math\'ematique de France,
  2016, pp.~147--172.

\bibitem[Mil17]{Milne}
J.S. Milne, \emph{Algebraic groups: the theory of group schemes of finite type
  over a field}, Cambridge studies in Advanced Math., vol. 170, Cambridge
  University Press, 2017.

\bibitem[MPS88]{DemTits}
L.~Michel, J.~Patera, and R.T. Sharp, \emph{The {D}emazure-{T}its subgroup of a
  simple {L}ie group}, J. Math. Phys. \textbf{29} (1988), no.~4, 777--796.

\bibitem[OT13]{OrlikTerao}
P.~Orlik and H.~Terao, \emph{Arrangements of hyperplanes}, Die Grundlehren der
  mathematischen Wissenschaften, vol. 300, Springer Science \& Business Media,
  2013.

\bibitem[OV90]{OV:LGAG}
A.~Onishchik and E.B. Vinberg, \emph{Lie groups and algebraic groups}, Springer
  Series in Soviet Mathematics, Springer, 1990.

\bibitem[PS85]{ProcesiSchwarz}
C.~Procesi and G.~Schwarz, \emph{Inequalities defining orbit spaces}, Invent.
  Math. \textbf{81} (1985), 539--554.

\bibitem[Ree10]{Reeder:tor}
M.~Reeder, \emph{Torsion automorphisms of simple {L}ie algebras}, Enseign.
  Math. (2) \textbf{56} (2010), no.~1--2, 3--47.

\bibitem[Sch24]{Scheiderer}
C.~Scheiderer, \emph{A course in real algebraic geometry}, Graduate Texts in
  Mathematics, vol. 303, Springer Nature, 2024.

\bibitem[Ser77]{Se:rep}
J-P. Serre, \emph{Linear representations of finite groups}, Graduate Texts in
  Mathematics, vol.~42, Springer, 1977.

\bibitem[Ser04]{Se:trace}
\bysame, \emph{On the values of the characters of compact {L}ie groups},
  Oberwolfach Reports \textbf{1} (2004), no.~1, 666--667.

\bibitem[Ser06]{Se:Kac}
\bysame, \emph{Coordonn{\'e}es de {K}ac}, Oberwolfach Reports \textbf{3}
  (2006), 1787--1790.

\bibitem[Ser23]{Se:zeros}
\bysame, \emph{Z\'eros de caract\`eres}, arXiv:2312.17551, December 2023.

\bibitem[Spr98]{Sp:LAG}
T.A. Springer, \emph{Linear algebraic groups}, second ed., Birkh\"{a}user,
  1998.

\bibitem[{Sta}]{Cheb}
{Stackoverflow contributors}, \emph{On the extrema of {C}hebyshev polynomials
  of the second kind}, available at
  \url{https://math.stackexchange.com/questions/2665518}.

\bibitem[{Sta}18]{stacks-project}
The {Stacks Project Authors}, \emph{Stacks project},
  \url{https://stacks.math.columbia.edu}, 2018.

\bibitem[Ste60]{St:inv}
R.~Steinberg, \emph{Invariants of finite reflection groups}, Canadian J. Math.
  \textbf{12} (1960), 616--618, [\#11 in \emph{Collected Papers}].

\bibitem[Ste65]{St:reg}
\bysame, \emph{Regular elements of semisimple algebraic groups}, Publ. Math.
  IHES (1965), no.~25, 49--80, [\#20 in \emph{Collected Papers}].

\bibitem[Tit66]{Ti:etendu}
J.~Tits, \emph{Normalisateurs de tores. {I}. {G}roupes de {C}oxeter \'etendus},
  J. Algebra \textbf{4} (1966), 96--116, [\#66 in \emph{{\OE}uvres--Collected
  Works}, vol.~{II}].

\bibitem[Wei94]{Weibel:HA}
C.A. Weibel, \emph{An introduction to homological algebra}, Cambridge Studies
  in Advanced Mathematics, no.~38, Cambridge University Press, 1994.

\end{thebibliography}

\newcommand{\etalchar}[1]{$^{#1}$}
\providecommand{\bysame}{\leavevmode\hbox to3em{\hrulefill}\thinspace}
\providecommand{\MR}{\relax\ifhmode\unskip\space\fi MR }
\providecommand{\MRhref}[2]{%
  \href{http://www.ams.org/mathscinet-getitem?mr=#1}{#2}
}
\providecommand{\href}[2]{#2}

\end{document}